\documentclass[aop,MSNbibl,citesort,dvips]{arximspdf}

\doi{10.1214/11-AOP739} 
\volume{41}
\issue{3B}
\pubyear{2013}
\firstpage{2103}
\lastpage{2181}

\makeatletter
\newcommand{\eqref}[1]{(\ref{#1})}
\newcommand{\Link}{{\mathrm{Link}}}
\newcommand{\salt}{{\mathrm{salt}}}
\newcommand{\alt}{{\mathrm{alt}}}
\newcommand{\Var}{\operatorname{Var}}
\newcommand{\Afrak}{{\mathfrak{A}}}
\newcommand{\Corr}{\operatorname{Corr}}
\newcommand{\Rub}{{\mathrm{Rub}}}
\newcommand{\Error}{{\mathrm{Err}}}
\newcommand{\Lfrak}{{\mathfrak{L}}}
\newcommand{\Mfrak}{{\mathfrak{M}}}
\newcommand{\TTT}{{\mathcal{T}}}
\newcommand{\tbold}{{\mathbf{t}}}
\newcommand{\ybold}{{\mathbf{y}}}
\newcommand{\zbold}{{\mathbf{z}}}
\newcommand{\MMM}{{\mathcal{M}}}
\newcommand{\BF}{{\mathrm{BF}}}
\newcommand{\Cbold}{{\mathbf{C}}}
\newcommand{\Mbold}{{\mathbf{M}}}
\newcommand{\Cfrak}{{\mathfrak{C}}}
\newcommand{\Bias}{{\mathrm{Bias}}}
\newcommand{\Seq}{{\mathrm{Seq}}}
\newcommand{\Part}{{\mathrm{Part}}}
\newcommand{\Gfrak}{{\mathfrak{G}}}
\newcommand{\supp}{\operatorname{supp}}
\newcommand{\sa}{{\mathrm{sa}}}
\newcommand{\trunc}{{\mathrm{trunc}}}
\newcommand{\Spec}{\operatorname{Spec}}
\newcommand{\Tfrak}{{\mathfrak{T}}}
\newcommand{\Ffrak}{{\mathfrak{F}}}
\newcommand{\Efrak}{{\mathfrak{E}}}
\newcommand{\SSS}{{\mathcal{S}}}
\newcommand{\T}{{\mathrm{T}}}
\newcommand{\VVV}{{\mathcal{V}}}
\newcommand{\WWW}{{\mathcal{W}}}
\newcommand{\DD}{{\mathbf{D}}}
\newcommand{\DDD}{{\mathcal{D}}}
\newcommand{\EEE}{{\mathcal{E}}}
\newcommand{\FFF}{{\mathcal{F}}}
\newcommand{\GGG}{{\mathcal{G}}}
\newcommand{\III}{{\mathcal{I}}}
\newcommand{\fbold}{{\mathbf{f}}}
\newcommand{\Ibold}{{\mathbf{I}}}
\newcommand{\Xbold}{{\mathbf{X}}}
\newcommand{\trace}{\operatorname{tr}}
\newcommand{\ii}{{\mathrm{i}}}
\newcommand{\AAA}{{\mathcal{A}}}
\newcommand{\BBB}{{\mathcal{B}}}
\newcommand{\one}{{\mathbf{1}}}
\newcommand{\Dbold}{{\mathbf{D}}}
\newcommand{\Ebold}{{\mathbb{E}}}
\newcommand{\CC}{{\mathbb{C}}}
\newcommand{\RR}{{\mathbb{R}}}
\newcommand{\Mat}{{\mathrm{Mat}}}
\newcommand{\GL}{{\mathrm{GL}}}
\newcommand{\HHH}{{\mathcal{H}}}
\newtheorem{Theorem}{Theorem} 
\newtheorem{Proposition}{Proposition} 
\newtheorem{Lemma}{Lemma} 
\newtheorem{Corollary}{Corollary} 
\newproclaim{Remark}{Remark} 
\newproclaim{Definition}{Definition} 
\makeatother

\begin{document}
\begin{frontmatter}

\title{Convergence of the largest singular value of a~polynomial in
independent Wigner matrices}
\runtitle{Polynomials in Wigner matrices}

\begin{aug}
\author[A]{\fnms{Greg W.} \snm{Anderson}\corref{}\ead[label=e1]{gwanders@umn.edu}}
\runauthor{G. W. Anderson}
\affiliation{University of Minnesota}
\address[A]{School of Mathematics\\
University of Minnesota\\
127 Vincent Hall\\
206 Church St. SE\\
Minneapolis, Minnesota 55455\\
USA\\
\printead{e1}} 
\end{aug}

\received{\smonth{3} \syear{2011}}
\revised{\smonth{12} \syear{2011}}

%
\begin{abstract}
For polynomials in independent Wigner matrices, we prove convergence of
the largest singular value to the operator norm of the corresponding
polynomial in free semicircular variables,
under fourth moment hypotheses.
We actually prove a more general result
of the form ``no eigenvalues outside the support of
the limiting eigenvalue distribution.''
We build on ideas of Haagerup--Schultz--Thorbj\o rnsen on the one hand
and Bai--Silverstein on the other.
We refine the linearization trick so as to preserve self-adjointness
and we develop a secondary trick bearing on the calculation of
correction terms.
Instead of Poincar\'{e}-type inequalities, we use a variety of matrix identities
and $L^p$ estimates.
The Schwinger--Dyson equation controls much of the analysis.
\end{abstract}

%

\end{frontmatter}

\maxcsec{10}
\tableofcontents
\eject
\section{Introduction and statement of the main result}
\subsection{Background and general remarks}
As part of a larger operator-theo\-retic investigation,
it was shown in~\cite{HST} (refining earlier work of~\cite{HT})
that there are for large $N$ almost surely no eigenvalues outside the
$\varepsilon$-neighborhood of the support of the
limiting spectral distribution
of a self-adjoint polynomial in independent GUE matrices.
(See~\cite{AGZ}, Chapter 5, Section 5, for another account of that result.)
It is natural to ask if the same is true for Wigner matrices.
We answer that question here in the affirmative.
To a large extent, this is a matter of learning to get by without the
Poincar\'{e} inequality.
Now the template for results of the form
``no eigenvalues outside the support\ldots'' was established a number of
years earlier
in the pioneering work of~\cite{BaiSil}, and moreover
the authors of that paper got along without the Poincar\'{e} inequality
quite well---erasure of rows and columns, classical $L^p$ estimates and
truncation arguments sufficed.
Moreover, they got their results under stringent fourth moment hypotheses.
In this paper, we channel the separately flowing streams of ideas of
\cite{BaiSil} and~\cite{HST} into one river, encountering a few
perhaps unexpected bends.

Any discussion of largest eigenvalues of Wigner matrices must mention the
classical work~\cite{BaiYin}. In that paper, convergence of the
largest eigenvalue of a Wigner matrix
to the spectrum edge was established under fourth moment hypotheses and
it was furthermore established
that in a certain sense fourth moments are optimal.

Our main result (see Theorem~\ref{TheoremMainResult} immediately below)
is both a ``polynomialization'' of the main result of~\cite{BaiYin}
and a generalization of the random matrix result of~\cite{HST}.
Roughly speaking, to prove our main result,
we let the results of~\cite{BaiYin} do the hard work of attracting the
eigenvalues to a compact neighborhood of the spectrum and then we draw them
the rest of the way in by using refinements of tools from both \cite
{BaiSil} and~\cite{HST},
among them matrix identities, $L^p$ estimates for quadratic forms in
independent random variables,
and the Schwinger--Dyson equation.

Generalizations of the ``no eigenvalues outside the support\ldots'' result
of~\cite{HST}
were quick to appear and continue to do so. In~\cite{Schultz},
following up
on the earlier results of~\cite{HT},
results in the GOE and GSE cases were obtained, and they
revealed a key role for ``correction terms'' of the sort we expend much effort
in this paper to control. In~\cite{CapDon}, a generalization
to non-Gaussian distributions satisfying Poincar\'{e}-type inequalities
was obtained.
In~\cite{Male}, a generalization was given
involving polynomials in GUE matrices
and deterministic matrices with convergent joint law which, in
particular, established
various rectangular analogues.

All the works following upon~\cite{HST} including this one build on
two extraordinarily powerful ideas from that paper:
(i) a counterintuitively ``backwards'' way of estimating the
error of approximate solutions of the Schwinger--Dyson equation
and (ii) the famous linearization trick.
We refine both ideas in this paper.
The refinements are closely intertwined and involve a gadget we call a
\textit{SALT block design}.

We have been significantly influenced by the paper~\cite{HRS07} which
explored geometry
and numerical analysis of the Schwinger--Dyson equation, and which could
serve uninitiated readers as an introduction to the use of matricial
semicircular elements. We were influenced also by~\cite{HLN07}
and~\cite{HLN08} which developed and applied Girko's notion of
deterministic equivalent.
The notion of deterministic equivalent is in effect exploited here as well,
but, more or less following~\cite{HST}, we simply harvest the needed
solutions of the Schwinger--Dyson equation from Boltzmann--Fock space
fully formed, thus avoiding iterative schemes for producing solutions.

There has been a lot of progress recently on universality in the bulk
and at the edge for
single Wigner matrices and sample covariance matrices.
Edge-universality results in the single matrix case
greatly refine and indeed render obsolete results of ``no eigenvalues
outside the support\ldots'' type,
albeit usually under more generous moment assumptions.
We mention for example~\cite{FeldheimSodin} which proves convergence
of the law of the suitably rescaled smallest eigenvalue of a sample
covariance matrix
with nonunity aspect ratio
to the Tracy--Widom distribution. Of course, many other papers could be
mentioned---the area is profoundly active at the moment.
It seems likely that similar edge-universality results are true in the
polynomial case. From this aspirational point of view, our results are crude.
But we hope they could serve as a point of departure.

\subsection{The main result}
We now formulate our main result, which at once generalizes the main
result of~\cite{BaiYin}
and the random matrix result of~\cite{HST}. Notation cursorily
introduced here
is explained in greater detail in Section \ref
{subsectionInitialNotation} below.

\subsubsection{Matrices with noncommutative polynomial entries}
Let ${\CC\langle\Xbold\rangle}$ be the noncommutative polynomial
algebra generated
over $\CC$ by a sequence
$\Xbold=\{\Xbold_\ell\}_{\ell=1}^\infty$ of independent
noncommuting algebraic variables.
We equip ${\CC\langle\Xbold\rangle}$ with an involution by
declaring all the variables
$\Xbold_\ell$ to be self-adjoint.
Given a sequence $a=\{a_\ell\}_{\ell=1}^\infty$ in an algebra $\AAA
$ and $f\in\Mat_n({\CC\langle\Xbold\rangle})$,
we define $f(a)\in\Mat_n(\AAA)$ by evaluating each entry of $f$ at
$\Xbold_\ell=a_\ell$
for all $\ell$. 

\subsubsection{Free semicircular variables}
Let
$\Xi=\{\Xi_\ell\}_{\ell=1}^\infty$ be a sequence of free
semicircular noncommutative random variables
in a faithful $C^*$-probability space. Given self-adjoint $f\in\Mat
_n({\CC\langle\Xbold\rangle})$,
let the law of the noncommutative random variable $f(\Xi)$ be denoted
by $\mu_f$.
The latter is a compactly supported probability measure on the real
line which depends only on the joint law of $\Xi$.
See Section~\ref{subsectionPositivityAndStates} below for a quick
review of
$C^*$-probability spaces and laws of single operators.
See Section~\ref{subsectionBoltzmannFock} for the Boltzmann--Fock
space construction which
yields a sequence $\Xi$ of the required type embedded in an algebraic
setup with further structures
useful to us.
For extensive discussion of noncommutative probability laws, including
joint laws,
and background on free probability,
see~\cite{AGZ}, Chapter~5.

\subsubsection{Random matrices}
Let
\[
\{\{x_\ell(i,j)\}_{1\leq i\leq j<\infty}\}_{\ell=1}^\infty
\]
be an array of independent $\CC$-valued random variables.
We assume the following for all $i$, $j$ and $\ell$:
%
%
\begin{eqnarray}
\label{equationBaiPoly1}
&\mbox{The law of $x_\ell(i,j)$ depends only on $\ell$ and $\one
_{i<j}$},&\\ \label{equationBaiPoly2}
&\Ebold|x_\ell(1,1)|^4<\infty\quad \mbox{and}\quad \Ebold|x_\ell
(1,2)|^4<\infty,&\\
\label{equationBaiPoly3}
&\Ebold x_\ell(1,1)=\Ebold x_\ell(1,2)=0 \quad\mbox{and}\quad \Ebold
|x_\ell(1,2)|^2=1,&\\
\label{equationBaiPoly4}
&\mbox{$x_\ell(1,1)$ is real-valued,}&\\
\label{equationBaiPoly5}
&\mbox{the real and imaginary parts of $x_\ell(1,2)$ are independent.}&
\end{eqnarray}
For positive integers $N$ and $\ell$, we then construct an
$N$-by-$N$ random hermitian matrix $X^{N}_\ell$ with entries
\[
X^{N}_\ell(i,j)=\cases{
x_\ell(i,j),&\quad$\mbox{if $i<j$,}$\vspace*{2pt}\cr
x_\ell(i,i),&\quad$\mbox{if $i=j$,}$\vspace*{2pt}\cr
x_\ell(j,i)^*,&$\quad\mbox{if $i>j$}$}
\]
and for each fixed $N$, we assemble these matrices into a sequence
$X^{N}=\{X^{N}_\ell\}_{\ell=1}^\infty$.
In turn, given self-adjoint $f\in\Mat_n({\CC\langle\Xbold\rangle
})$, let
$\nu^N_f$ be the empirical distribution\vspace*{1pt} of eigenvalues of the random
hermitian matrix $f(\frac{X^{N}}{\sqrt{N}})$.
We use the notation $\nu^N_f$ rather than, say, $\mu^N_f$ because we
are saving the latter
for use in our main technical result, namely Theorem \ref
{TheoremMainResultBis} below.

The next result is essentially
well known and provides context for our main result.
\begin{Theorem}\label{TheoremPreMainResult}
For all self-adjoint $f\in\Mat_n({\CC\langle\Xbold\rangle})$, the
empirical distribution
$\nu^N_f$ converges weakly to $\mu_f$ as $N\rightarrow\infty$,
almost surely.
\end{Theorem}

See~\cite{VDN} or~\cite{AGZ}, Chapter~5, for background, similar
results, and many references.
See~\cite{MeckesSzarek} for an interesting recent approach to the
proof of a similar result.
For the reader's convenience we give in Section~\ref
{subsubsectionByproduct} below a quick derivation of Theorem~\ref
{TheoremPreMainResult} from our main technical result.

Now we can state our main result.

\begin{Theorem}\label{TheoremMainResult}
For every self-adjoint $f\in\Mat_n({\CC\langle\Xbold\rangle})$
and $\varepsilon>0$,
$\supp\nu_f^N$ is contained in
the $\varepsilon$-neighborhood of $\supp\mu_f$ for $N\gg0$, almost
surely.\vadjust{\goodbreak}
\end{Theorem}

See Section~\ref{subsubsectionSidestep} below for the derivation of this result
from our main technical result.

The next corollary justifies the title of this paper.

\begin{Corollary}\label{CorollaryPostMainResult}
For every $f\in\Mat_n({\CC\langle\Xbold\rangle})$,
\[
\lim_{N\rightarrow\infty} \biggl[\!\biggl [f\biggl(\frac{X^{N}}{\sqrt {N}}\biggr)\biggr]\! \biggr]=[\![f(\Xi)]\!]\qquad\mbox{a.s.}
\]
\end{Corollary}

\begin{pf}
After replacing $f$ by $ff^*$, we may assume that $f$ is self-adjoint,
and furthermore that $f(\frac{X^N}{\sqrt{N}})$ and $f(\Xi)$ are positive.
We then need only show that the largest eigenvalue of
$f(\frac{X^N}{\sqrt{N}})$ converges as $N\rightarrow
\infty$ to the
largest element of the spectrum of $f(\Xi)$, almost surely.
In any case, $\Spec(f(\Xi))=\supp\mu_f$ by the very important
Lemma~\ref{LemmaFaithfulSupport} below
and thus
\[[\![f(\Xi)]\!]=\sup\Spec(f(\Xi))=\sup\supp\mu_f.
\]
%
Finally, we have
\[[\![f(\Xi)]\!]\leq\liminf_{N\rightarrow\infty}\biggl[\!\biggl[f \biggl(\frac
{X^N}{\sqrt{N}}\biggr) \biggr]\!\biggr]
\leq\limsup_{N\rightarrow\infty}\biggl[\!\biggl[f\biggl(\frac{X^N}{\sqrt {N}}\biggr)\biggr]\!\biggr]\leq[\![f(\Xi)]\!]\qquad \mbox{a.s.}
\]
by Theorem~\ref{TheoremPreMainResult} on the left
and Theorem~\ref{TheoremMainResult} on the right.
\end{pf}

\subsection{Outline of the paper and plan of proof}
\subsubsection{Truncation and reduction steps}
In Section~\ref{sectionModel}, after introducing general notation and
terminology,
we make the truncation step common to the proofs of Theorems~\ref
{TheoremPreMainResult} and~\ref{TheoremMainResult}.
Then we formulate our main technical result, namely Theorem~\ref
{TheoremMainResultBis} below,
which concerns $L^p$-norms of ``randomized and corrected'' Stieltjes transforms,
and from it we derive Theorems~\ref{TheoremPreMainResult} and~\ref
{TheoremMainResult}.
The proof of the main technical result then
takes up the remainder of the paper.

\subsubsection{The self-adjoint linearization trick}
To launch the ``blocky'' approach taken in the rest of the paper,
at the end of Section~\ref{sectionModel} we present
a very simple self-adjointness-preserving variant of the linearization trick.
In Remark~\ref{RemarkPrePreCobbling} below,
we explain how the trick gives access to the Stieltjes transforms
considered in Theorem~\ref{TheoremMainResultBis}.

\subsubsection{General tools}
In Section~\ref{sectionOperatorTheoryTools}, we review elementary
topics concerning $C^*$-algebras
and in particular we recall the Boltzmann--Fock space construction.
We also consider an ad hoc version of the notion of Schur
complement in a $C^*$-algebra
and use it to solve an abstract version of the Schwinger--Dyson equation.
In Section~\ref{sectionConcentration}, after introducing tensor
products and norming rules,
we write down an ensemble of mostly familiar estimates
that we will use in place of the Poincar\'{e} inequality.
\subsubsection{$\SSS$-(bi)linear machinery}
In Section~\ref{sectionUpgrade}, we introduce a collection of
algebraic tools needed
to take advantage of the fine structure possessed by the random
matrices and operators
described in Remark~\ref{RemarkPrePreCobbling}. In particular, we
introduce the notion of SALT block design
to streamline the self-adjoint linearization trick
and we develop a secondary trick for making new SALT block designs from
old. The secondary trick is a ``bootstrapping'' technique indispensable
to our study of corrections.

\subsubsection{Study of the Schwinger--Dyson equation}
In Section~\ref{sectionOpTheoSchwingerDyson}, we recall the
Schwinger--Dyson (SD) equation
and we construct solutions of it by using
the apparatus of Section~\ref{sectionOperatorTheoryTools}.
Following~\cite{HST}, we use certain of these solutions to represent
the Stieltjes transform $S_{\mu_f}(z)$
figuring in Theorems~\ref{TheoremPreMainResult},~\ref{TheoremMainResult}
and~\ref{TheoremMainResultBis}. See Remark~\ref{RemarkCobblingBis} below.
We next introduce a secondary version of
the SD equation involving notions introduced in Section~\ref{sectionUpgrade}
and we show how a solution of it can be extracted from the ``upper right
corner'' of a solution of a suitably chosen
(larger and more complicated)
instance of the SD equation itself.
We then construct our candidate for the correction ${\mathrm
{bias}}_f^N(z)$ figuring in
Theorem~\ref{TheoremMainResultBis}.
See Remark~\ref{RemarkCobblingTer} below.
In Section~\ref{sectionMagSetup}, working in a relatively simple geometry,
we refine the idea of~\cite{HST} for controlling errors of approximate
solutions to the SD equation.
By means of the secondary trick, we will be able to use the estimates
of Section~\ref{sectionMagSetup} not only to
study the convergence of empirical distributions of eigenvalues to
their limits, but also to study the limiting behavior of corrections.

\subsubsection{Matrix identities and $L^p$ estimates}
In Section~\ref{sectionUrMatrixIdentities}, we present a carefully
edited catalog of identities
satisfied by objects built out of finite-sized chunks of an infinite
matrix with entries which are themselves
matrices of some fixed finite size.
One among these objects via the self-adjoint linearization trick
specializes to the randomized Stieltjes transform\vspace*{1pt} $S_{\mu^N_f}(\zbold
)$ figuring in Theorem~\ref{TheoremMainResultBis}.
See Remark~\ref{RemarkCobbling} below. We note also that identity
\eqref{equationBias} of Section~\ref{subsectionBias}
is the ultimate source of all the correction terms studied here.
In Section~\ref{sectionMatrixIdentities}, we introduce the block
Wigner model and work through a long series of
$L^p$ estimates culminating in Theorem~\ref{TheoremBias} below which
converts identity \eqref{equationBias} to a crucial approximation. We
emphasize that all the arguments and calculations presented
in Sections~\ref{sectionUrMatrixIdentities} and~\ref
{sectionMatrixIdentities}
make sense for Wigner matrices
when specialized to the case in which the constituent blocks are copies
of $\CC$.
In many cases, the calculations so specialized then run along familiar lines.
A reader who has already developed some intuition about Wigner matrices
should, we hope, be able to build on that base in order to understand
our work.

\subsubsection{Concluding arguments}
Finally,
in Section~\ref{sectionEndgame},
we combine the tools collected above to complete the proof of Theorem
\ref{TheoremMainResultBis}
in relatively short order.

\section{The truncation step and the main technical result}
\label{sectionModel}
In Section~\ref{subsectionInitialNotation},
we introduce general notation in force throughout the paper.
In Section~\ref{subsectionTruncationStep}, we carry out the
truncation step for proving Theorems
\ref{TheoremPreMainResult} and~\ref{TheoremMainResult}.
In Section~\ref{subsectionStieltjesTransforms}, we recall a method
for reconstructing a probability measure from its Stieltjes transform.
In Section~\ref{subsectionModel},
we formulate our main technical result,
namely Theorem~\ref{TheoremMainResultBis} below,
and we explain how to check its most important hypotheses efficiently using
the classical estimate of~\cite{FurKom}.
In Section~\ref{subsectionFromMainToBis},
we recover both Theorems~\ref{TheoremPreMainResult}
and~\ref{TheoremMainResult} from Theorem~\ref{TheoremMainResultBis}.
The proof of the latter result will then take up the rest of the paper.
Finally, in Section~\ref{subsectionTheLinearizationTrick}, we
introduce a simple self-adjointness-preserving variant
of the famous linearization trick of~\cite{HST}, thereby banishing
nonlinear noncommutative polynomials from
further consideration in the main body of the paper.

\subsection{Notation and terminology} \label{subsectionInitialNotation}
\subsubsection{General notation}
We use $\vee$ and $\wedge$ for maximum and minimum, respectively.
Given a complex number $z\in\CC$, let
$\Re z=\frac{z+z^*}{2}$ and $\Im z=\frac{z-z^*}{2\ii}$,
and put ${\mathfrak{h}}= \{z\in\CC\vert\Im z>0\}$,
which is the classical upper half-plane.
Let $\Ebold$ denote expectation and let $\Pr$ denote probability.
(We save the letters $E$ and $P$ for other purposes.)
We write $\one_A$ for the indicator of an event $A$.
Let $\supp\nu$ denote the support of a probability measure $\nu$,
and similarly, let $\supp\varphi$ denote the support of a function
$\varphi$.
(Recall that supports are always closed sets.)
For any $\CC$-valued random variable $Z$
and exponent $p\in[1,\infty]$, let ${\Vert Z \Vert}_p$ denote the
$L^p$-norm of $Z$,
that is, let ${\Vert Z \Vert}_p=(\Ebold|Z|^p)^{1/p}$ for $p\in
[1,\infty)$
and otherwise
let ${\Vert Z \Vert}_\infty$ denote the essential supremum of $|Z|$.
For a matrix $A$ with complex entries, let $A^*$ denote the transpose
conjugate, $A^\T$ the transpose
and $[\![A]\!]$ the largest singular value of $A$. More generally, we use
$[\![\cdot]\!]$ to denote the norm on a $C^*$-algebra.
We denote the spectrum of an element $x$
of a $C^*$-algebra $\AAA$ by $\Spec(x)$; Proposition \ref
{PropositionSubStable} below justifies
omission of reference to $\AAA$ in this notation.
\begin{Remark}
We use the not-so-standard notation $[\![\cdot]\!]$ for $C^*$-norms
in order not to collide with
the notation ${\Vert\cdot \Vert}_p$ for $L^p$-norms of random variables.
We will in fact have to consider expressions of the form ${\Vert[\![A]\!] \Vert}_p$ rather frequently.
\end{Remark}

\subsubsection{Algebras and matrices}
An \textit{algebra} $\AAA$ always has $\CC$ as scalar field, is associative,
and possesses a unit denoted by $1_\AAA$. (Other notation for the unit
may also be used, e.g., simply $1$.)
Let $\Mat_n(\AAA)$ denote the algebra of $n$-by-$n$ matrices with
entries in $\AAA$.
More generally, let $\Mat_{k\times\ell}(\AAA)$ denote the space of
$k$-by-$\ell$ matrices with entries in $\AAA$.
The $(i,j)$-entry of a matrix $A$ is invariably denoted $A(i,j)$
(never $A_{ij}$). Let $\AAA^\times$ denote the group of invertible
elements of an algebra $\AAA$,
put $\GL_n(\AAA)=\Mat_n(\AAA)^\times$ ($\GL$ for \textit{general
linear} group)\vadjust{\goodbreak}
and for $A\in\Mat_n(\AAA)$, let $\trace_\AAA A=\sum_{i=1}^n A(i,i)$.
In the special case \mbox{$\AAA=\CC$}, we write $\trace=\trace_\CC$.
Let $\Ibold_n\in\Mat_n(\CC)$ denote the $n$-by-$n$ identity matrix
and more generally, given an element $a\in\AAA$,
let $\Ibold_n\otimes a\in\Mat_n(\AAA)$ denote the diagonal matrix
with entries $a$ on the diagonal. Given a \mbox{$*$-algebra $\AAA$}, that is,
an algebra endowed with an involution denoted $*$,
and an element $a\in\AAA$,
we say that $a$ is \textit{self-adjoint}
if $a^*=a$ and we denote the set of such elements by
$\AAA_{\sa}$. Given a matrix $A\in\Mat_{k\times\ell}(\AAA)$ with entries
in a \mbox{$*$-algebra~$\AAA$},
we define $A^*\in\Mat_{\ell\times k}(\AAA)$
by $A^*(i,j)=A(j,i)^*$. In particular, by this rule $\Mat_n(\AAA)$
becomes a $*$-algebra whenever
$\AAA$ is.\vspace*{-2pt}

\subsubsection{\texorpdfstring{The noncommutative polynomial ring ${\CC\langle\Xbold\rangle}$}
{The noncommutative polynomial ring C<X>}}

Let ${\CC\langle\Xbold\rangle}$ be the noncommutative polynomial
ring generated
over $\CC$ by a sequence
$\Xbold=\{\Xbold_\ell\}_{\ell=1}^\infty$ of independent
noncommuting variables.
By definition, the family of all \textit{monomials}
\[
\bigcup_{m=0}^\infty\{\Xbold_{i_1}\cdots \Xbold_{i_m}\vert i_1,\ldots
,i_m=1,2,3,\ldots\}
\]
(including the empty monomial, which is identified to $1_{\CC\langle
\Xbold\rangle}$)
forms a Hamel basis for the vector space underlying $\CC\langle\Xbold
\rangle$.
In particular, $\CC\langle\Xbold\rangle=\bigcup_{m=1}^\infty\CC
\langle\Xbold_1,\break\ldots,\Xbold_m\rangle$.
We equip ${\CC\langle\Xbold\rangle}$ with $*$-algebra structure by
the rule
$\Xbold_\ell^*=\Xbold_\ell$ for all $\ell$.
Let $S^\infty$ denote the space of sequences in a set $S$.
Given an algebra $\AAA$, a sequence $a\in\AAA^\infty$
and matrix $f\in\Mat_n({\CC\langle\Xbold\rangle})$,
let $f(a)\in\Mat_n(\AAA)$ denote the matrix obtained by evaluating
each entry
at $\Xbold=a$ (and evaluating $1_{{\CC\langle\Xbold\rangle}}$ to
$1_\AAA$).
Note that if $\AAA$ is a $*$-algebra and $a\in\AAA^\infty_\sa$,
then $f(a)^*=f^*(a)$,
that is, the evaluation map $f\mapsto f(a)$ is a $*$-algebra homomorphism.
If $\AAA=\Mat_N(\CC)$,
then we view $f(a)$ as an $n$-by-$n$ array of
$N$-by-$N$ blocks, thus identifying it with an element of $\Mat
_{nN}(\CC)$.

\subsubsection{Empirical distributions of eigenvalues}
Given an $N$-by-$N$ hermitian matrix $A$, the
\textit{empirical distribution} of its eigenvalues $\lambda_1\leq\cdots
\leq\lambda_N$
is by definition the probability
measure
$\frac{1}{N}\sum_{i=1}^N \delta_{\lambda_i}$ on the real line.

\subsection{The truncation step}\label{subsectionTruncationStep}
We begin the proofs of Theorems~\ref{TheoremPreMainResult} and \ref
{TheoremMainResult}
by proving the following proposition.
In doing this, we are imitating the initial truncation step taken in
\cite{BaiSil}.
\begin{Proposition}\label{PropositionTruncationStep}
To prove Theorems~\ref{TheoremPreMainResult} and \ref
{TheoremMainResult}, we may
augment assumptions~\eqref{equationBaiPoly1}--\eqref{equationBaiPoly5}
without loss of generality
by the following assumptions holding for every index $\ell$:
%
%
\begin{eqnarray}
\label{equationBaiPoly6}
&\ii^{-\ell}x_\ell(1,2)\qquad \mbox{is real-valued,}&\\
\label{equationBaiPoly7}
&x_\ell(1,1)=0\qquad \mbox{for odd $\ell$,}&\\
\label{equationBaiPoly8}
&{\Vert x_\ell(1,1) \Vert}_\infty<\infty \quad\mbox{and}\quad
{\Vert x_\ell(1,2) \Vert}_\infty<\infty.&\vspace*{-2pt}
\end{eqnarray}
\end{Proposition}

We collect several tools before completing the proof in
Section~\ref{subsubsectionTruncationEnd} below.
The first tool is simply a couple of standard eigenvalue inequalities.\vadjust{\goodbreak}

\begin{Lemma}\label{LemmaEigenvalueSwarm}
For $A,B\in\Mat_N(\CC)_{\sa}$ let
$\lambda_i(A)$ and $\lambda_i(B)$ denote the $i${th} largest
eigenvalue, respectively.
Then we have \textup{(i)} $\bigvee_{i=1}^N|\lambda_i(A)-\lambda_i(B)|\leq
[\![A-B]\!]$ and \textup{(ii)} the corresponding empirical distributions
are within distance
$[\![A-B]\!]$ as measured in the Lipschitz bounded metric.
\end{Lemma}

Recall that the distance of probability measures $\mu$ and $\nu$
on the real line in the \textit{Lipschitz bounded metric}
is the supremum of $|\int\varphi \,d\mu-\int\varphi \,d\nu|$ where
$\varphi\dvtx \RR\rightarrow\RR$ ranges over functions
with supremum norm and Lipschitz constant both $\leq1$.
Recall also that the Lipschitz-bounded metric is compatible with weak
convergence.
\begin{pf*}{Proof of Lemma~\ref{LemmaEigenvalueSwarm}} (i) This is well-known. See~\cite{HoJo} or~\cite{Simon}.
(ii) For any test function $\varphi\dvtx \RR\rightarrow\RR$ with sup
norm and Lipschitz constant both $\leq1$, since $|\varphi(x)-\varphi
(y)|\leq|x-y|$,
we have $|\int\varphi \,d\mu_A-\int\varphi \,d\mu_B|\leq
[\![A-B]\!]$
by part (i) of the lemma.
\end{pf*}

\subsubsection{The Bai--Yin model}
Let
\[
\{w(i,j)\}_{1\leq i\leq j<\infty}
\]
be an independent family of real random variables
such that the law of $w(i,j)$ depends only on $\one_{i<j}$.
Assume furthermore that $w(1,1)$ and $w(1,2)$ have finite fourth
moments and zero means.
Let $\sigma={\Vert w(1,2) \Vert}_2$.
Given a positive integer $N$, let
$W^N$ be the $N$-by-$N$ random real symmetric matrix with entries
\[
W^N(i,j)=\cases{
w(i,j),&\quad$\mbox{if $i\leq j$,}$\vspace*{2pt}\cr
w(j,i),&\quad$\mbox{if $i>j$.}$}
\]
To have a convenient catchphrase, let us call
\[
\{W^N\}_{N=1}^\infty
\]
the \textit{Bai--Yin model} for Wigner matrices.
We have the following fundamental result.
\begin{Theorem}[(\cite{BaiYin}, Theorem~C)]\label{TheoremBaiYin}
In the Bai--Yin model $\{W^N\}_{N=1}^\infty$, the largest eigenvalue
of $\frac{W^N}{\sqrt{N}}$ converges to $2\sigma$ as $N\rightarrow
\infty$, almost surely.
\end{Theorem}

\begin{Remark}
By~\cite{BaiYin}, Theorem~A, the fourth moment hypothesis
in Theorem~\ref{TheoremBaiYin} cannot be improved while maintaining strong
overall assumptions concerning the form of the joint law of the family
$\{W^N\}$
and in particular enforcing the identification of $W^N$ with the upper
left $N$-by-$N$ block of~$W^{N+1}$.
\end{Remark}

\begin{Remark}
It is trivial but useful to observe that
Theorem~\ref{TheoremBaiYin} continues to hold in the case $\sigma
=0$, that is,
in the case in which $W^N$ is diagonal. In this case, the proof is just
an exercise in applying the Borel--Cantelli lemma.\vadjust{\goodbreak}
\end{Remark}

\begin{Remark}\label{RemarkOfCourse}
Only real symmetric matrices were treated in~\cite{BaiYin} but
all the arguments carry over to the hermitian case.
In particular, Theorem~\ref{TheoremBaiYin}
continues to hold if we replace $\{W^N\}_{N=1}^\infty$ by the slightly
altered family
\[
\{\{(\ii\one_{i<j}-\ii\one_{i>j})W^N(i,j)\}_{i,j=1}^N\}
_{N=1}^\infty.
\]
Just to have a convenient catchphrase (and to avoid introducing yet
more notation) let us call the latter
the \textit{twisted} Bai--Yin model.
\end{Remark}

\subsubsection{$C$-truncation}
Given a $\CC$-valued random variable $Z$ such that\break ${\Vert Z \Vert
}_2=1$ and
$\Ebold Z=0$,
along with a constant $C>0$, put
\[
\rho_C(Z)={\bigl\Vert Z\one_{|Z|\leq C}-\Ebold Z\one_{|Z|\leq C}
\bigr\Vert}_2,\qquad
\theta_C(Z)={\bigl\Vert Z\one_{|Z|>C}-\Ebold Z\one_{|Z|>C} \bigr\Vert}_2,
\]
and if $\rho_C(Z)>0$ put
\[
\trunc_C(Z)=
\bigl(Z\one_{|Z|\leq C}-\Ebold Z\one_{|Z|\leq C}\bigr)/\rho_C(Z).
\]
Note that
%
%
\begin{equation}\label{equationDCapp}
\rho_C(Z)\rightarrow_{C\rightarrow\infty} 1\quad
\mbox{and}\quad \theta_C(Z) \rightarrow_{C\rightarrow\infty} 0
\end{equation}
by dominated convergence.

\begin{Lemma}\label{LemmaTruncApprox}
Consider again the Bai--Yin model $\{W^N\}_{N=1}^\infty$. But now
assume that
$\sigma=1$. Let $C>0$ be large enough so that $ \rho_C(w(1,2))>0$.
Let $\widehat{W}^N$ be the result of applying the truncation procedure
$\trunc_C$ to the off-diagonal
entries of $W^N$ and putting the diagonal entries of $W^N$ to $0$.
\textup{(i)} We have
%
%
\begin{equation}\label{equationTruncApprox}
\limsup_{N\rightarrow\infty}\biggl[\!\biggl[\frac{W^N-\widehat {W}^N}{\sqrt
{N}}\biggr]\!\biggr]\leq2(\theta+1-\rho)\qquad
\mbox{a.s.},
\end{equation}
where $\theta=\theta_C(w(1,2))$ and $\rho=\rho_C(w(1,2))$.
\textup{(ii)} The analogous statement holds for the twisted Bai--Yin model.
\end{Lemma}

\begin{pf}
Let $D^N$ be the result of putting the off-diagonal entries of $W^N$ to zero.
We have in any case a bound
\[\biggl[\!\biggl[\frac{ W^N-\widehat{W}^N }{\sqrt{N}}\biggr]\!\biggr]\leq
\biggl[\!\biggl[\frac{D^N}{\sqrt{N}}\biggr]\!\biggr]+
\biggl[\!\biggl[\frac{W^N-D^N-\rho\widehat{W}^N}{\sqrt{N}}\biggr]\!\biggr]+\biggl[\!\biggl[\frac {(1-\rho
)\widehat{W}^N}{\sqrt{N}}\biggr]\!\biggr]\quad
\mbox{a.s.}
\]
The terms on the right-hand side almost surely tend as $N\rightarrow\infty$ to $0$,
$2\theta$ and $2(1-\rho)$, respectively, by Theorem~\ref{TheoremBaiYin}.
Thus, (i) is proved and (ii) is proved similarly.
\end{pf}

\subsubsection{\texorpdfstring{Proof of Proposition \protect\ref{PropositionTruncationStep}}
{Proof of Proposition 1}}\label{subsubsectionTruncationEnd}

The permissibility of assuming \eqref{equationBaiPoly6}
and \eqref{equationBaiPoly7} is clear---one has only to break the
originally given system of
matrices into symmetric and antisymmetric pieces, rescale and relabel.
We may assume \eqref{equationBaiPoly6}
and \eqref{equationBaiPoly7} henceforth.
Now fix $f\in\Mat_n({\CC\langle\Xbold\rangle})_\sa$ and
$\varepsilon>0$ arbitrarily.
With a large constant $C_\ell>0$
depending on $\ell$, to be aptly chosen presently, let $\widehat
{X}_\ell^N$ be the result of applying
the truncation operation $\trunc_{C_\ell}$ to the off-diagonal\vadjust{\goodbreak}
entries of $X_\ell^N$
and putting the diagonal entries to $0$.
Let $\hat{\nu}_f^N$ denote the empirical distribution of eigenvalues
of $f(\frac{\widehat{X}^N}{\sqrt{N}})$.
By Theorem~\ref{TheoremBaiYin}, Remark~\ref{RemarkOfCourse}, Lemma
\ref{LemmaTruncApprox}
and \eqref{equationDCapp}
we can choose constants $C_\ell$ large enough, depending on $f$, so that
\[
\limsup_{N\rightarrow\infty}
\biggl[\!\biggl[f\biggl(\frac{X^N}{\sqrt{N}}\biggr)-f\biggl(\frac{\widehat {X}^N}{\sqrt{N}}\biggr)\biggr]\!\biggr]<\frac{\varepsilon}{2}
\qquad\mbox{a.s.}
\]
By Lemma~\ref{LemmaEigenvalueSwarm}(i),\vspace*{-1pt} almost surely for $N\gg0$,
we have that $\supp\nu^N_f$ is contained
in the $\frac{\varepsilon}{2}$-neighborhood
of $\supp\hat{\nu}^N_f$, and in turn,
by Theorem~\ref{TheoremMainResult} applied under the additional
assumptions of Proposition~\ref{PropositionTruncationStep},
almost surely for $N\gg0$,
we have that $\hat{\nu}^N_f$ is contained
in the $\frac{\varepsilon}{2}$-neighborhood
of $\supp\mu_f$.
Thus, Theorem~\ref{TheoremMainResult} in the general case
follows from the special case considered in Proposition \ref
{PropositionTruncationStep}.
A~similar argument using Lemma~\ref{LemmaEigenvalueSwarm}(ii)
derives Theorem~\ref{TheoremPreMainResult} in the general case
from the special case considered in Proposition \ref
{PropositionTruncationStep}.

\subsection{Stieltjes transforms and reconstruction of probability measures}
\label{subsectionStieltjesTransforms}
We recall an important tool and motivate the introduction of the
auxiliary upper-half-plane-valued random variable $\zbold$.

\subsubsection{Stieltjes transforms}
In general, given a probability measure $\mu$ on the real line,
recall that the \textit{Stieltjes transform} is defined by the formula
\[
S_\mu(z)=\int\frac{\mu(dt)}{t-z} \qquad\mbox{for $z\in\CC\setminus
\supp\mu$.}
\]
We use here the same sign convention as (say) in~\cite{BaiSil}
so that $\Im z>0\Rightarrow\Im S_\mu(z)>0$.
Recall also that
%
%
\begin{equation}\label{equationUniformStieltjesBound}
S_\mu(z^*)\equiv S_\mu(z)^* \quad\mbox{and}\quad |S_\mu(z)\Im z|\leq1.
\end{equation}
In particular, $S_\mu$ is real-valued on $\RR\setminus\supp\mu$.

\subsubsection{The auxiliary random variable $\zbold$}\label
{subsubsectionzbold}
Let $m$ be an even positive integer. Let~$\zbold$ be an ${\mathfrak
{h}}$-valued random variable the law of which is specified by the
integration formula
\[
\Ebold\varphi(\zbold)=\int_0^\infty
\int_{-\infty}^\infty\varphi(x+\ii y)
\frac{e^{-(x^2+y^2)/2}y^m}{(m-1)!!\pi} \,dx \,dy.
\]
Note that $\zbold={\mathbf{x}}+\ii\ybold$ where ${\mathbf{x}}$ and
$\ybold$ are independent,
${\mathbf{x}}$ is standard Gaussian, and~$\ybold$ has density that
vanishes to order $m$ at $0$.
We call $m$
the \textit{strength} of the \textit{repulsion} of~$\zbold$
from the real axis. For simplicity, we assume that $\Im\zbold>0$
holds without exception.
In general, we allow $m$ to vary
from one appearance of $\zbold$ to the next. Results below involving
$\zbold$ are
often stated with hypotheses to the effect that $m$ be sufficiently large.
As we will see, the exact distribution of $\zbold$ is not too important.
But it is quite important that ${\Vert1/\Im\zbold \Vert}_p<\infty$ for
$p\in[1,m+1)$.
Thus, by choosing the strength of the repulsion of $\zbold$\vadjust{\goodbreak} from the
real axis large enough,
the random variable $1/\Im\zbold$ can be made to possess as many
finite moments as we like.

The method we will use for reconstructing probability measures from
their Stieltjes transforms is codified by the following lemma
in which the auxiliary random variable $\zbold$ enters as a convenience
for bookkeeping.
\begin{Lemma}\label{LemmaReconstruction}
Let $\varphi\dvtx \RR\rightarrow\RR$ be infinitely differentiable
and compactly supported. Then there exists a function
$\Upsilon\dvtx \CC\rightarrow\CC$ depending on the strength of repulsion
of $\zbold$
from the real axis with the following properties:
\begin{longlist}[(III)]
\item[(I)] $\Upsilon$ is infinitely differentiable and compactly supported.
Furthermore, $\Upsilon$~satisfies $\supp\Upsilon\cap\RR=\supp
\varphi$
and has the symmetry $\Upsilon(z^*)\equiv\Upsilon(z)^*$.
\item[(II)] For any open set $D\subset\CC$ such that $D^*=D\supset
\supp\Upsilon$
and analytic function $b\dvtx  D\rightarrow\CC$ such that $b(z^*)\equiv
b(z)^*$, we have
$\Re\Ebold\Upsilon(\zbold)b(\zbold)=0$.
\item[(III)] For probability measures $\mu$ on $\RR$, we have
$\Re\Ebold\Upsilon(\zbold)S_\mu(\zbold)=\int\varphi \,d\mu$.
\end{longlist}
\end{Lemma}

The lemma mildly refines a procedure buried in the proof of~\cite{AGZ}, Lem\-ma~5.5.5.
\begin{pf*}{Proof of Lemma~\ref{LemmaReconstruction}}
We identify $\CC$ with $\RR^2$ in the customary way.
We switch back and forth between writing $x+\ii y$ and $(x,y)$ as it
suits us.
To begin the construction, let $\theta\dvtx \RR\rightarrow[0,1]$ be an even
infinitely differentiable function supported
in the interval $[-1,1]$ and identically equal to $1$ on the subinterval
$[-\frac{1}{2},\frac{1}{2}]$. Let $m$ denote the strength of the
repulsion of $\zbold$
from the real axis.
Put
\[
\Gamma(x, y)=\frac{1}{2\pi}\theta(y)\sum_{j=0}^{m}\frac{(\ii
y)^j}{j!}\varphi^{(j)}(x),
\]
noting that\vspace*{-1pt} $\Gamma$ is supported in $\supp\varphi\times[-1,1]$. Put
$\Gamma'(x,y)=(\frac{\partial}{\partial x}+\ii\frac{\partial
}{\partial y})\Gamma(x,y)$,
noting that $\Gamma'(z^*)\equiv\Gamma'(z)^*$.
The significance of the differential
operator $\frac{\partial}{\partial x}+\ii\frac{\partial}{\partial
y}$ is that it kills all analytic functions,
that is, it encodes the Cauchy--Riemann equations.
The sum defining $\Gamma(x,y)$ is contrived so that
\[
\Gamma'(x,y)=\frac{1}{2\pi}\frac{(\ii y)^m}{m!}\varphi^{(m+1)}(x)
\qquad
\mbox{for }(x,y)\in\RR\times\biggl(-\frac{1}{2},\frac{1}{2}\biggr).
\]
Let\vspace*{1.5pt} $\rho(x,y)=\frac{y^me^{-(x^2+y^2)/2}}{(m-1)!!\pi}$.
Then we have
$2\Gamma'(x,y)=
\Upsilon(x,y)\rho(x,y)$ for some function $\Upsilon$ satisfying (I).
For any Borel measurable function $h\dvtx \CC\rightarrow\CC$ satisfying
$h(z)^*\equiv h(z^*)$ almost everywhere with respect to Lebesgue measure
we have
%
%
\begin{equation}\label{equationUglyZboldFormula}
\Re\Ebold
\Upsilon(\zbold)h(\zbold)
=\int_{-\infty}^\infty\int_{-\infty}^\infty\Gamma'(x,y)h(x,y)
\,dx \,dy
\end{equation}
provided that the integral on the right is absolutely convergent,
as follows directly from the definition of $\Upsilon$.
Furthermore, for any compact set\vadjust{\goodbreak} $T\subset\RR^2$ with a polygonal boundary
and analytic function $h$ defined in a neighborhood of $T$ we have
%
%
\begin{equation}\label{equationGreensTheorem}
\int_T
\Gamma' h \,dx \,dy
=-\ii\int_{\partial T} \Gamma h (dx+\ii \,dy)
\end{equation}
by Green's theorem and the fact that $h$ is killed by $\frac{\partial
}{\partial x}+\ii\frac{\partial}{\partial y}$.
To prove (II), take~$T$ such that $\supp\Gamma\subset T\setminus
\partial T\subset T\subset D$
and take $h=b$. Then formulas \eqref{equationUglyZboldFormula} and
\eqref{equationGreensTheorem} yield the result.
To prove (III),
assume at first that $\mu=\delta_t$ for some real $t$
and hence $S_\mu(z)=\frac{1}{t-z}$.
Take $T$ to be an annulus centered at $t$ and take $h=\frac{1}{t-z}$.
In the limit as the inner radius tends to $0$ and the outer radius
tends to $\infty$,
formulas \eqref{equationUglyZboldFormula} and \eqref{equationGreensTheorem}
yield the result.
Finally, to get (III) in general, use Fubini's theorem---the hypotheses
of the latter
hold by \eqref{equationUniformStieltjesBound} and the fact that
$m\geq1$.
\end{pf*}

\subsection{The main technical result}\label{subsectionModel}
Now we introduce a new model, the one we actually study through most of
the paper,
and we formulate a general statement about it, namely Theorem \ref
{TheoremMainResultBis} below.
All the hypotheses for this model are ones
we deserve to make after performing the truncation step in the proofs of
Theorems~\ref{TheoremPreMainResult} and~\ref{TheoremMainResult}.
Furthermore, in certain respects, hypotheses are actually weakened in
comparison to those for
Theorems~\ref{TheoremPreMainResult} and~\ref{TheoremMainResult}.

\subsubsection{Data}\label{subsubsectionModelData}
For integers $\ell,N\geq1$, fix a random element
$\Xi^N_\ell$ of $\Mat_N(\CC)_\sa$.
Fix also an independent family $\{\FFF(i,j)\}_{1\leq i\leq j<\infty}$
of $\sigma$-fields.
Let $\FFF$ denote the $\sigma$-field generated by all the $\FFF(i,j)$.

\subsubsection{Assumptions}
We assume for each $p\in[1,\infty)$ and index $\ell$ the following:
%
%
\begin{eqnarray}
\label{equationPepper1}
\sup_N\bigvee_{i,j=1}^N{\Vert\Xi^N_\ell(i,j) \Vert}_p&<&\infty,\\
\label{equationPepper2}
\sup_N{\biggl\Vert\biggl[\!\biggl[\frac{\Xi^N_\ell}{\sqrt{N}}\biggr]\!\biggr]
\biggr\Vert}_p&<&\infty.
\end{eqnarray}
Furthermore, we assume for all indices $\ell$ and $N$ the following:
%
%
\begin{eqnarray}
\label{equationPepper4}
&\mbox{$\Xi^N_\ell$ is the upper left $N$-by-$N$ block of $\Xi
^{N+1}_\ell$,}&\\
\label{equationPepper5}
&\mbox{$(\Xi_\ell^N)^\T=(-1)^\ell\Xi_\ell^N$,}&\\
\label{equationPepper6}
&\mbox{$\Xi_\ell^N(i,j)$ is $\FFF(i\wedge j,i\vee j)$-measurable and
$\Ebold\Xi^N_\ell(i,j)=0$}&
\nonumber
\\[-8pt]
\\[-8pt]
\eqntext{\mbox{for }i,j=1,\ldots,N.}\\
\label{equationPepper7}
& {\Vert\Xi^N_\ell(i,j) \Vert}_2=1\qquad \mbox{for }1\leq i<j\leq N.&
\end{eqnarray}
%
Finally, we assume that
%
%
\begin{equation}
\label{equationPepper8}
\Ebold\Xi^N_\ell(i,j)\Xi^N_{m}(i,j)=0\qquad
\mbox{for $1\leq i<j\leq N$
and $1\leq\ell<m<\infty$}
\end{equation}
for all positive integers $i$, $j$, $N$, $\ell$ and $m$ subject to the
indicated constraints.

\begin{Remark}
While moment assumptions here are extremely generous
in comparison to those of Theorem~\ref{TheoremMainResult},
in certain other respects we have significantly weakened assumptions.
Firstly, we do not require the entries $\Xi^N_\ell(i,j)$
to have law depending only on $\ell$ and the sign of $i-j$.
Secondly, we assume somewhat less than strict independence of the
matrices $\Xi^N_\ell$ for fixed $N$ and varying~$\ell$.
\end{Remark}

\subsubsection{Random matrices and empirical distributions of eigenvalues}
For each fixed $N$, we form a sequence
$\Xi^N=\{\Xi^N_\ell\}_{\ell=1}^\infty\in\Mat_N(\CC)_{\sa
}^\infty$
of random hermitian matrices.
Given any self-adjoint $f\in\Mat_n({\CC\langle\Xbold\rangle})_\sa
$, let
$\mu^N_f$ denote\vspace*{-1pt} the empirical distribution of eigenvalues of the random
hermitian matrix $f(\frac{\Xi^N}{\sqrt{N}})\in\Mat
_{nN}(\CC)_{\sa}$.
\subsubsection{The auxiliary random variable $\zbold$}
We adjoin the auxiliary random variable $\zbold$ figuring in Lemma
\ref{LemmaReconstruction}
to our model. We assume that $\zbold$ is independent of~$\FFF$.



\begin{Theorem}\label{TheoremMainResultBis}
Fix $f\in\Mat_n(\CC\langle\Xbold\rangle)_\sa$ arbitrarily
and let $\mu_f$ be defined as in Theorem~\ref{TheoremMainResult}.
Then there exists a sequence
\[
\{{\mathrm{bias}}_f^N\dvtx \CC\setminus\supp\mu_f\rightarrow\CC\}
_{N=1}^\infty
\]
of deterministic analytic functions with the symmetry
%
%
\begin{equation}\label{equationBiasSymmetry}
{\mathrm{bias}}_f^N(z^*)\equiv{\mathrm{bias}}_f^N(z)^*
\end{equation}
such that for every $p\in[1,\infty)$
we have
%
%
\begin{eqnarray}
\label{equationMainResultBis1}
\sup_N N^{1/2}{\Vert S_{\mu^N_f}(\zbold)-S_{\mu_f}(\zbold )
\Vert}_p&<&\infty,\\
\label{equationMainResultBis2}
\sup_N N^{3/2}{\Vert S_{\mu^{N+1}_f}(\zbold)-S_{\mu^N_f}(\zbold)
\Vert}_p
&<&\infty,\\
\label{equationMainResultBis3}
\sup_N {\Vert{\mathrm{bias}}_f^N(\zbold) \Vert}_p&<&\infty
\end{eqnarray}
and
\begin{equation}
\label{equationMainResultBis4}
\sup_N N^2{\biggl\Vert\Ebold(S_{\mu^N_f}(\zbold) \vert \zbold
)-S_{\mu_f}(\zbold)-\frac{{\mathrm{bias}}_f^N(\zbold )}{N} \biggr\Vert
}_p<\infty,
\end{equation}
provided the strength of the repulsion of $\zbold$ from the real axis
is sufficiently great, depending on $p$.
\end{Theorem}

Once we have deduced Theorems~\ref{TheoremPreMainResult}
and~\ref{TheoremMainResult}
from Theorem~\ref{TheoremMainResultBis} in Section~\ref
{subsectionFromMainToBis} below, the proof of Theorem \ref
{TheoremMainResultBis} will take up the rest of the paper.\vadjust{\goodbreak}

\begin{Remark}\label{RemarkWillChooseSpecialXi}
The theorem if true for some sequence $\Xi$ of free semicircular
noncommutative random variables
in a faithful $C^*$-probability space is true for all; only the joint
law of $\Xi$ is important. Taking advantage of this freedom,
we will make a special choice of $\Xi$ below
in Section~\ref{subsectionBoltzmannFock} which is adapted to the
symmetry present
in our model by virtue of assumption \eqref{equationPepper5}.
\end{Remark}

\begin{Remark}\label{equationSmearToPoint}
Fix a point $z_0\in{\mathfrak{h}}$ arbitrarily.
For any analytic function $g\dvtx {\mathfrak{h}}\rightarrow\CC$,
we can recover the value $g(z_0)$ as the average
of $g(z)$ over the disc\vspace*{1pt} $|z-z_0|\leq\frac{1}{2}\Im z_0$.
Thus statement \eqref{equationMainResultBis1} for, say, $p=4$ implies
that\break
$S_{\mu_f^N}(z_0)\rightarrow_{N\rightarrow\infty} S_{\mu_f}(z_0)$,
almost surely,
by Jensen's inequality in conditional form and the Borel--Cantelli lemma.
In short, the averaged result stated in Theorem \ref
{TheoremMainResultBis} easily yields
pointwise results.
\end{Remark}

\begin{Remark}
An explicit if rather complicated description of ${\mathrm
{bias}}_f^N(z)$ in operator-theoretic terms
will be developed below. See Remark~\ref{RemarkCobblingTer} below and the
discussion of the Schwinger--Dyson equation which immediately precedes
that remark.
The role played by ${\mathrm{bias}}_f^N(z)$ below,\vspace*{-1pt}
justified by Lemma~\ref{LemmaReconstruction} and relation~\eqref
{equationMainResultBis4},
is to make possible an estimate of $\int\varphi \,d(\mu_f^N-\mu_f)$
accurate to an order in $1/N$ sufficiently high so that we can achieve
sensitivity to the movement of
individual eigenvalues.
\end{Remark}

\begin{Remark}\label{RemarkFurKomDeployment}
Suppose that in the setup for Theorem~\ref{TheoremMainResultBis},
we discard assumptions~\eqref{equationPepper1} and \eqref{equationPepper2}
and make in their place for every index $\ell$ the assumption
%
%
\begin{equation}\label{equationPepper9}
\sup_{N=1}^\infty\bigvee_{i,j=1}^N{\Vert\Xi^N_\ell(i,j) \Vert
}_\infty
<\infty.
\end{equation}
Then of course assumption \eqref{equationPepper1} holds trivially,
but furthermore and importantly, assumption \eqref{equationPepper2}
and the bound
%
%
\begin{equation}
\label{equationPepper3}
{\biggl\Vert\limsup_{N\rightarrow\infty}\biggl[\!\biggl[\frac{\Xi^N_\ell}{\sqrt
{N}}\biggr]\!\biggr]\biggr \Vert}_\infty<\infty
\end{equation}
hold by the classical result of~\cite{FurKom} recalled immediately
below in a convenient form.
\end{Remark}

\begin{Proposition}\label{PropositionFurKom}
For each $N\geq1$, let $Y^N$ be a random $N$-by-$N$ hermitian matrix
whose entries on or above the diagonal are independent.
Assume furthermore that the entries of the
matrices $Y^N$ are essentially bounded uniformly in $N$ and have mean zero.
Fix any sequence $\{k_N\}_{N=1}^\infty$ of positive integers such that
$\frac{k_N}{\log N}\rightarrow\infty$ but
$\frac{k_N}{N^{1/6}}\rightarrow0$.
Then
$\sum_N \Ebold[\![\frac{Y^N}{c\sqrt{N}}]\!]^{2k_N}<\infty$
for some (finite) constant $c>0$.
\end{Proposition}

Here and elsewhere throughout the paper constants in
estimates are denoted mostly by $c$, $C$ or $K$.
The numerical values of these constants may of course vary from context
to context and even from line to line.
\begin{pf*}{Proof of Proposition~\ref{PropositionFurKom}}
We will use the result of F\"{u}redi--Koml\'{o}s as cast in the form of
the combinatorial estimate
\cite{AGZ}, Lemma~2.1.23. By the cited lemma, for any constants
\[
\frac{c}{2}>K\geq\sup_N \bigvee_{i,j=1}^N {\Vert Y^N(i,j) \Vert
}_\infty,
\]
we deduce via ``opening of the brackets'' and counting of nonzero terms that
\begin{eqnarray*}
\Ebold\biggl[\!\biggl[\frac{Y^N}{c\sqrt{N}}\biggr]\!\biggr]^{2k_N} &\leq& \Ebold
\trace\biggl(\frac{Y^N}{c\sqrt{N}}\biggr)^{2k_N}\\
&\leq&\frac{1}{c^{2k_N}}\sum_{t=1}^{k_N+1} 2^{2k_N}
(2k_N)^{3(2k_N-2t+2)} N^t\frac{K^{2k_N}}{N^{k_N}}\\
&=& N\biggl(\frac{2K}{c}\biggr)^{2k_N}\sum_{t=1}^{k_N+1}\biggl(\frac
{2k_N}{N^{1/6}}\biggr)^{6(k_N-t+1)},
\end{eqnarray*}
whence the result, since the last expression summed on $N$ is finite.
\end{pf*}
\begin{Remark}
The argument presented immediately after~\cite{AGZ}, Lem-\break ma~2.1.23~gives
the analogous result for Wigner-like random matrices whose
entries have $L^p$-norms uniformly under a bound polynomial in $p$. We
do not need the stronger result here for any of our proofs, but we
mention it because it easily produces many natural examples of data
satisfying the assumptions
of Theorem~\ref{TheoremMainResultBis}
and furthermore satisfying \eqref{equationPepper3}.
\end{Remark}

\subsection{\texorpdfstring{Deduction of Theorems \protect\ref{TheoremPreMainResult} and \protect\ref{TheoremMainResult} from Theorem \protect\ref{TheoremMainResultBis}}
{Deduction of Theorems 1 and 2 from Theorem 4}}
\label{subsectionFromMainToBis}

\subsubsection{Common setup for the proofs}
In order to deduce Theorems~\ref{TheoremPreMainResult} and~\ref
{TheoremMainResult}
from Theorem~\ref{TheoremMainResultBis},
we may and we do make the additional assumptions stated in Proposition
\ref{PropositionTruncationStep}.
In turn, in order to apply Theorem~\ref{TheoremMainResultBis},
we now take
\[
\Xi^N_\ell=X^N_\ell
\]
for all $N$ and $\ell$ and we put
\[
\FFF(i,j)=\sigma(\{x_\ell(i,j)\}_{\ell=1}^\infty).
\]
Then---hypothesis \eqref{equationPepper2} excepted---the data
\[
\{\Xi^N_\ell\}\cup\{\FFF(i,j)\}
\]
trivially satisfy
all hypotheses of Theorem~\ref{TheoremMainResultBis}, for example,
the algebraic assumptions \eqref{equationBaiPoly6} and \eqref
{equationBaiPoly7}
imply the (anti)symmetry \eqref{equationPepper5}.\vadjust{\goodbreak}
And hypothesis \eqref{equationPepper2} is fulfilled---not so
trivially---by Remark~\ref{RemarkFurKomDeployment}. So Theorem \ref
{TheoremMainResultBis}
is indeed applicable in the present case. Note that \eqref
{equationPepper3} holds as well,
either by yet another application of Theorem~\ref{TheoremBaiYin} or
by Remark~\ref{RemarkFurKomDeployment}.

\subsubsection{\texorpdfstring{Proof of Theorem \protect\ref{TheoremPreMainResult} with Theorem \protect\ref{TheoremMainResultBis} granted}
{Proof of Theorem 1 with Theorem 4 granted}}
\label{subsubsectionByproduct}
By \eqref{equationPepper3}, there exists a constant $A>0$
such that $\supp\mu_f^N\subset[-A,A]$ for $N\gg0$, almost surely,
so we have tightness.
By Remark~\ref{equationSmearToPoint}, we have
$S_{\mu_f^N}(\ii+1/k)\rightarrow_{N\rightarrow\infty} S_{\mu
_f}(\ii+1/k)$,
almost surely, for every integer $k>0$. The latter statement
by standard subsequencing arguments (which we omit) implies
that $\mu_f^N$ indeed converges weakly to~$\mu_f$, almost surely.

To derive Theorem
\ref{TheoremMainResult} from Theorem~\ref{TheoremMainResultBis},
we need a final lemma variants of which have long been in use.
\begin{Lemma}\label{LemmaGiddyUp}
Let $\{Y_N\}_{N=1}^\infty$ be a sequence of nonnegative random variables.
Assume that
\[
\sup_N N\Ebold Y_N<\infty
\quad\mbox{and}\quad \sup_N N^{1/2}{\Vert Y_{N+1}-Y_N \Vert}_4<\infty.
\]
Then $Y_N\rightarrow_{N\rightarrow\infty} 0$, almost surely.
\end{Lemma}

\begin{pf}
We have $Y_{\lfloor k^{5/4}\rfloor}\rightarrow_{k\rightarrow\infty} 0$,
almost surely, by the Chebyshev inequality and the Borel--Cantelli lemma.
Here, $\lfloor x\rfloor$ denotes the greatest integer not exceeding $x$.
Put $[N]=\bigvee_{k=1}^\infty\lfloor k^{5/4}\rfloor\one_{k^{5/4}<N}$.
Clearly,\vspace*{1pt} we have\break $Y_{[N]}\rightarrow_{N\rightarrow\infty}0$, almost surely.
Since $N-[N]=O(N^{1/5})$, we have
$\Vert Y_N-\break Y_{[N]} \Vert_4=O(N^{-3/10})$
by the Minkowski inequality.
Thus $Y_N-Y_{[N]}\rightarrow_{N\rightarrow\infty}0$, almost surely,
by the Chebyshev inequality and the Borel--Cantelli lemma.
The result follows.
\end{pf}

\subsubsection{\texorpdfstring{Proof of Theorem \protect\ref{TheoremMainResult} with Theorem \protect\ref{TheoremMainResultBis} granted}
{Proof of Theorem 2 with Theorem 4 granted}}
\label{subsubsectionSidestep}
Take $\zbold$ to have a strength of repulsion from the real axis large
enough so that all statements of
Theorem~\ref{TheoremMainResultBis} hold for the given matrix $f\in
\Mat_n(\CC\langle\Xbold\rangle)_\sa$ in the case $p=4$.
As in the proof of Theorem~\ref{TheoremPreMainResult}, fix $A>0$ such
that $\supp\mu_f^N\subset[-A,A]$ for $N\gg0$, almost surely. Fix
$\varepsilon>0$ arbitrarily.
Fix an infinitely differentiable function $\varphi\dvtx \RR\rightarrow
[0,1]$ with the following support properties:
\begin{itemize}
\item$\varphi$ is identically equal to $1$ on $[-A,A]$ minus the
$\varepsilon$-neighborhood of $\supp\mu_f$.
\item$\varphi$ is supported in some compact set disjoint from $\supp
\mu_f$.
\end{itemize}
For $N>0$ consider the nonnegative random variable
$Y_N=nN\int\varphi \,d\mu_f^{N}$
the value of which for $N\gg0$ bounds the number of eigenvalues of
the random hermitian matrix $f(\frac{\Xi^{N}}{\sqrt{N}})$
straying outside the $\varepsilon$-neighborhood of $\supp\mu_f$, almost surely.
It will be enough to show that $Y_N\rightarrow_{N\rightarrow\infty
}0$, almost surely.
Now by Lemma~\ref{LemmaReconstruction} and Fubini's theorem, for some
compactly supported
infinitely differentiable function
$\Upsilon\dvtx \CC\rightarrow\CC$ with support disjoint from
$\supp\mu_f$,
we have for each $N>0$ the representation
$Y_N=nN\Re\Ebold( \Upsilon(\zbold)S_{\mu_f^N}(\zbold)\vert\FFF
)$, almost surely.
Furthermore, by similar reasoning,
for any analytic function $b\dvtx \CC\setminus\supp\mu_f\rightarrow\CC
$ satisfying $b(z^*)\equiv b(z)^*$,
we have $\Re\Ebold(\Upsilon(\zbold)b(\zbold)\vert\FFF)=0$,
almost surely.
From statements \eqref{equationMainResultBis1}
and \eqref{equationMainResultBis2} with $p=4$ we deduce that
$\sup_N N^{{1}/{2}}{\Vert Y_{N+1}-Y_N \Vert}_4<\infty$
via Jensen's inequality in conditional form.
From statements \eqref{equationMainResultBis3} and
\eqref{equationMainResultBis4},
we deduce that $\sup_N N\Ebold Y_N<\infty$.
Thus $Y_N\rightarrow_{N\rightarrow\infty}0$, almost surely, by Lemma
\ref{LemmaGiddyUp},
which finishes the proof. 

\subsection{The self-adjoint linearization trick}\label{subsectionTheLinearizationTrick}
We now present a simple self-adjoint-ness-preserving variant of
the celebrated linearization trick of~\cite{HT} and~\cite{HST},
and we explain how the trick gives access to the
Stieltjes transforms $S_{\mu^N_f}(z)$ and $S_{\mu_f}(z)$ figuring in
Theorem~\ref{TheoremMainResultBis}.
This will motivate our focus on block-decomposed matrices and operators
in the sequel.

\subsubsection{Schur complements}
Recall the familiar formula
%
%
\begin{equation}\label{equationFamiliarInversion}
\left[
\matrix{
a&b\vspace*{2pt}\cr
c&d}
\right]^{-1}=\left[
\matrix{
0&0\vspace*{2pt}\cr
0&d^{-1}}
\right]
+\left[
\matrix{
1\vspace*{2pt}\cr
-d^{-1}c}\right]
(a-bd^{-1}c)^{-1}
\left[
\matrix{
1&-bd^{-1}}
\right]
\end{equation}
for inverting a block-decomposed matrix.
Formula \eqref{equationFamiliarInversion} holds whenever $d$ is invertible
and at least one of $\bigl[{
a\enskip b\atop
c\enskip d}
\bigr]$ and the so-called \textit{Schur complement}
$a-bd^{-1}c$ are invertible, in which case both of the latter two
matrices are invertible.

\begin{Proposition}\label{PropositionNaiveSALT}
Fix $f\in\Mat_n({\CC\langle\Xbold\rangle})_\sa$ arbitrarily.
Then for some integer $s>n$ there exists
$\tilde{f}\in\Mat_s({\CC\langle\Xbold\rangle})_\sa$ all entries
of which are of
degree $\leq1$ in the variables
$\Xbold_\ell$ and which furthermore admits a block decomposition
$\tilde{f}=\bigl[
{0\enskip b\atop
b^*\enskip d}\bigr]$ such that $d\in\Mat_{s-n}({\CC\langle\Xbold\rangle})_\sa$ is
invertible
and $f=-bd^{-1}b^*$.
\end{Proposition}

We call $\tilde{f}$ a \textit{self-adjoint linearization} of $f$.
Later we will upgrade this definition to a slightly more sophisticated
form. See Definition~\ref{DefinitionSALTupgrade} below.
\begin{pf*}{Proof of Proposition~\ref{PropositionNaiveSALT}} Let $\deg f$ be the maximum of the degrees of the
noncommutative polynomials appearing in $f$.
We proceed by induction on $\deg f$.
There is no difficulty to find matrices $f_1,f_2\in\Mat_{n\times
N}({\CC\langle\Xbold\rangle})$
for some $N$ such that
$f=f_1f_2^*+f_2f_1^*$ and $\deg f_1\vee\deg f_2\leq1\vee\frac{2\deg f}{3}$.
Put $b=\bigl[
{
f_1\enskip f_2}
\bigr]$ and $d=-\bigl[
{
0\atop\Ibold_N}\enskip{
\Ibold_N\atop 0}
\bigr]$. Note that $d=d^{-1}$ and $f=-bdb^*$.
If $\deg f\leq1$, then $g=\bigl[
{
0\atop b^*}\enskip
{b\atop d}
\bigr]$ is already a self-adjoint linearization of $f$
and we are done. Otherwise, by induction on $\deg f$, the matrix
$g$ has a self-adjoint linearization $\tilde{g}=\left[\fontsize{8.36}{10.36}\selectfont{\matrix{0&0&x\cr
0&0&b_1\cr
x^{*}&  b_1^{*}&  d_1}}
\right]$\vspace*{1pt}
where the zero block in the upper left is $n$-by-$n$,
the central zero block is $2N$-by-$2N$
and the other blocks are of appropriate sizes. By \eqref
{equationFamiliarInversion} and the definitions, we then have
\[
\left[
\matrix{
0&b\vspace*{2pt}\cr
b^*&d}
\right] = -\left[
\matrix{
x\vspace*{2pt}\cr
b_1}
\right]d_1^{-1}\left[
\matrix{
x^*&b_1^*}
\right],\vspace*{-6pt}
\]
\begin{eqnarray*}
&&-\left[
\matrix{
0&x}
\right]
\left[
\matrix{
0&b_1\vspace*{2pt}\cr
b^*_1&d_1}
\right]^{-1}\left[
\matrix{
0\vspace*{2pt}\cr
x^*}
\right]\\
&&\qquad=-xd_1^{-1}x^*+xd_1^{-1}b_1^*(b_1d_1^{-1}b_1^*)^{-1}b_1d^{-1}x^*=-bd^{-1}b^*=f.
\end{eqnarray*}
Thus $\tilde{g}$ is a self-adjoint linearization of $f$ as well as of
$g$.
\end{pf*}

\begin{Proposition}\label{PropositionNaiveSALTbis}
Let $f\in\Mat_n({\CC\langle\Xbold\rangle})$ and $\tilde{f}=\bigl[
{0\enskip b\atop
b^*\enskip d}
\bigr]\in\Mat_s({\CC\langle\Xbold\rangle})$ be as in Proposition \ref
{PropositionNaiveSALT}.
Let $\AAA$ be any algebra, let $a=\{a_\ell\}_{\ell=1}^\infty$ be
any sequence,
fix $\Lambda\in\Mat_n(\AAA)$ and let $\tilde{\Lambda}=\bigl[
{
\Lambda\atop 0}\enskip
{0\atop 0}
\bigr]\in\Mat_s(\AAA)$.
Then $f(a)-\Lambda$ is invertible if and only if $\tilde{f}(a)-\tilde
{\Lambda}$ is invertible
and under these equivalent conditions
the upper left $n$-by-$n$ block of $(\tilde{f}(a)-\tilde{\Lambda
})^{-1}$ equals $(f(a)-\Lambda)^{-1}$.
\end{Proposition}

\begin{pf} This follows immediately from \eqref
{equationFamiliarInversion} and the definitions.
\end{pf}
\begin{Remark}\label{RemarkPrePreCobbling}
(In this remark, we take for granted the notation and constructions of
Section~\ref{subsectionBoltzmannFock} and Section \ref
{subsectionTensorProducts} below. There is no danger of circularity
since we will not pick up the discussion thread
initiated here again until Section~\ref{sectionUpgrade}.)
Proposition~\ref{PropositionNaiveSALTbis} says in the context of
Theorem~\ref{TheoremMainResultBis}
that
%
%
\begin{equation}\label{equationEmpiricalAccess}
 S_{\mu_f^N}(z)=
\frac{1}{Nn}
\sum_{i=1}^{Nn}
\biggl(\tilde{f}\biggl(\frac{\Xi^N}{\sqrt{N}}\biggr)-\left[
\matrix{
z\Ibold_{Nn}&0\vspace*{2pt}\cr
0&0}
\right]\biggr)^{-1}(i,i)
\end{equation}
for $z\in{\mathfrak{h}}$ and combined with Proposition \ref
{PropositionPreMainResult} below says
furthermore that
%
%
\begin{eqnarray}
\label{equationPreStieltjesAccess}
\qquad\supp\mu_f& = &\Spec(f(\Xi))
= \left\{z\in\CC\Big\vert\tilde{f}(\Xi)-\left[
\pmatrix{
z\Ibold_{n}&0\vspace*{2pt}\cr
0&0}
\right] \mbox{ not invertible}\right\},\\
\label{equationStieltjesAccess}
S_{\mu_f}(z)& = &
\varphi^\BF\Biggl(\frac{1}{n}\sum_{i=1}^n\left(\tilde{f}(\Xi
)-\left[
\pmatrix{
z\Ibold_{n}&0\vspace*{2pt}\cr
0&0}
\right]\right)^{-1}(i,i)\Biggr)
\end{eqnarray}
for $z\in\CC\setminus\supp\mu_f$.
Now for $m\gg0$ and suitable
$a_0,\ldots,a_m\in\Mat_s(\CC)_\sa$ we have
%
%
\begin{equation}\label{equationTypicalBlock}
\tilde{f}=a_0\otimes1_{{\CC\langle\Xbold\rangle}}+\sum_{\ell
=1}^m a_\ell\otimes
\Xbold_\ell.
\end{equation}
In order to gain control of the Stieltjes transforms $S_{\mu_f^N}(z)$
and $S_{\mu_f}(z)$,
all our efforts hereafter will be directed toward understanding
the special properties of block-decomposed random matrices of the form
%
%
\begin{equation}\label{equationTypicalBlockRandom}
\tilde{f}\biggl(\frac{\Xi^N}{\sqrt{N}}\biggr)=a_0\otimes\Ibold
_N+\sum_{\ell=1}^m a_\ell\otimes\frac{\Xi^N_\ell}{\sqrt{N}}\in
\Mat_{sN}(\CC)_\sa
\end{equation}
and of block-decomposed operators of the form
%
%
\begin{equation}\label{equationTypicalBlockFock}
\tilde{f}(\Xi)=a_0\otimes1_{B(\HHH)}+\sum_{\ell=1}^m a_\ell
\otimes\Xi_\ell\in\Mat_s(B(\HHH))_\sa.
\end{equation}
%
Many tools from analysis and algebra will come into play.
\end{Remark}

\section{Tools from operator theory}
\label{sectionOperatorTheoryTools}
We review some topics in elementary $C^*$-algebra theory.
We emphasize the viewpoint of noncommutative probability.
We recall noncommutative laws, relations between laws and spectra,
the construction of Boltzmann--Fock space,
and Schur complements in the $C^*$-algebra context.
Finally, we solve an abstract algebraic version of the Schwinger--Dyson equation.

\subsection{Warmup exercises} We record without proof
some very elementary facts frequently used below. Recall that we only use
algebras $\AAA$ possessing a unit $1_\AAA$.
\begin{Lemma}
\label{LemmaNeumannExpansion}
Let $x$ and $y$ be elements of a Banach algebra with $x$ invertible
and $2[\![x^{-1}]\!][\![y]\!]\leq1$.
Then $x-y$ is invertible and $[\![(x-y)^{-1}]\!]\leq2[\![x^{-1}]\!]$.
\end{Lemma}

Here and below we invariably
use $[\![\cdot]\!]$ to denote the norm on a Banach algebra.
\subsubsection{The resolvent identity}
We note the \textit{resolvent identity}
%
%
\begin{equation}\label{equationResolventIdentity}
x^{-1}-y^{-1}=y^{-1}(y-x)x^{-1}=x^{-1}(y-x)y^{-1} \qquad(x,y\in\AAA
^\times)
\end{equation}
holding in any algebra $\AAA$ and its infinitesimal variant
$\frac{d}{dt}x^{-1}=-x^{-1}\frac{dx}{dt}x^{-1}$. We also need
the iterated version
%
%
\begin{eqnarray}\label{equationResolventIdentityIterated}
x^{-1}-y^{-1}=y^{-1}(y-x)y^{-1}+y^{-1}(y-x)y^{-1}(y-x)x^{-1}
\nonumber
\\[-8pt]
\\[-8pt]
\eqntext{(x,y\in\AAA^\times).}
\end{eqnarray}
%

\subsection{Positivity}\label{subsectionPositivity}
$\!\!$We recall basic facts about positive elements of $C^*$-algebras.

\subsubsection{Positive elements and their square roots}
If an element $x$ of a $C^*$-algebra $\AAA$ is self-adjoint with
nonnegative spectrum,
we write $x\geq0$; and if furthermore $x$ is invertible, then we write $x>0$.
Elements satisfying $x\geq0$ are called \textit{positive}.
Elements of the form $xx^*$ are automatically positive.
For $x\in\AAA$ such that $x\geq0$, there exists unique $y\in\AAA$
such that $y\geq0$ and $y^2=x$ (see~\cite{Murphy}, Theorem~2.2.1), in
which case we write $x^{1/2}=y$.

\subsubsection{$C^*$-subalgebras and GNS}
Let $\AAA$ be a $C^*$-algebra.
We say that a closed subspace $\AAA_0\subset\AAA$
is a \textit{$C^*$-subalgebra} if $\AAA_0$ is stable under $*$,
closed under multiplication and furthermore $1_\AAA\in\AAA_0$,
in which case $\AAA_0$ is a \mbox{$C^*$-algebra} in its own right for which
$1_{\AAA_0}=1_\AAA$.\vadjust{\goodbreak}
Each $C^*$-algebra is isomorphic to a \mbox{$C^*$-subalgebra} of $B(H)$ for
some Hilbert space $H$ via the GNS construction (see~\cite{Murphy}, Section 3.4).

\begin{Proposition}\label{PropositionSubStable}
For any $C^*$-algebra $\AAA$ and $C^*$-subalgebra $\AAA_0\subset\AAA$
we have $\AAA_0\cap\AAA^\times=\AAA_0^\times$.
\end{Proposition}

(See~\cite{Murphy}, Theorem~2.1.11.) Thus, the spectrum
of $x\in\AAA_0$ is the same whether viewed in $\AAA_0$ or $\AAA$.
In particular, $x$ is positive in $\AAA_0$ if and only if positive in
$\AAA$.

\begin{Proposition}\label{PropositionCstarNormUnique}
For every element $x$ of a $C^*$-algebra $\AAA$, if $x$ is normal,
and in particular, if $x$ is self-adjoint,
then $[\![x]\!]$ equals the spectral radius of~$x$. Consequently,
$[\![x]\!]^2$ equals the spectral radius of $xx^*$ and $x^*x$.
\end{Proposition}

(See~\cite{Murphy}, Theorem~2.1.1 and Corollary~2.1.2.)
It follows that a $*$-algebra can be normed as $C^*$-algebra in at most
one way.
We always use that norm when it exists.

\subsubsection{Real and imaginary parts}
Given any $*$-algebra\vspace*{1pt} and $Z\in\AAA$ we write $\Re Z=\frac{Z+Z^*}{2}$
and $\Im Z=\frac{Z-Z^*}{2\ii}$. (This generalizes the notation we
already have
for real and imaginary parts of a complex number.)

The next elementary result plays an vitally important role in the paper.

\begin{Lemma}
\label{Lemmalambda}
Let $\AAA$ be a $C^*$-algebra. Let $A\in\AAA$ satisfy $\Im A\geq0$
and let $z\in{\mathfrak{h}}$.
Then $A+ z1_\AAA\in\AAA^\times$ and $[\![(A+z1_\AAA )^{-1}]\!]\leq
1/\Im z$.
\end{Lemma}

\begin{pf} To abbreviate, we write $1=1_\AAA$, $z=z1_\AAA$, and so on.
After replacing $A$ by $(A+\Re z)/\Im z$, we may assume without loss of
generality that $z=\ii$.
Write $A=X+\ii Y$ with $X=\Re A$ and $Y=\Im A$. Since $Y\geq0$, we
have $1+Y>0$,
and hence we can write
$A+\ii=(1+Y)^{1/2}(W+\ii)(1+Y)^{1/2}$
where $W=(1+Y)^{-1/2}X(1+Y)^{-1/2}\in\AAA_\sa$.
Since both $(1+Y)^{1/2}$ and $W+\ii$ are normal
and have spectra disjoint from the open unit disc centered at the origin,
both are invertible with inverse of norm $\leq1$ by Proposition \ref
{PropositionCstarNormUnique}.
Thus, $A+\ii$ is invertible with inverse of norm $\leq1$.
\end{pf}

\subsection{States and spectral theory}\label
{subsectionPositivityAndStates} We recall some
basic definitions and results pertaining to $C^*$-probability spaces.
Much of this background is covered in~\cite{Murphy}. The rest of it
is more or less implicit in~\cite{Murphy} and~\cite{VDN} but hard to
extract. Some of this material is also covered in~\cite{AGZ},
Chapter~1,
but unfortunately Lemma~\ref{LemmaDykemaHelp} below is not. For the
reader's convenience, we supply short proofs of some key statements
which are part of standard ``$C^*$-know-how'' but hard to pin down in
the literature.

\subsubsection{States}\label{subsubsectionStates}
Let $\AAA$ be a $C^*$-algebra. Let $\phi\dvtx \AAA\rightarrow\CC$ be
any linear functional (perhaps not bounded).
One calls $\phi$ \textit{positive} if for every $A\in\AAA$, if $A\geq
0$, then $\phi(A)\geq0$,
in which case $\phi$ is automatically bounded and satisfies $\phi
(x^*)=\phi(x)^*$.
One calls $\phi$ a \textit{state} if $\phi$ is positive and $\phi
(1_\AAA)=1$,
in which case $[\![\phi]\!]=1$.
One calls a state $\phi$ \textit{faithful} if
for every $A\in\AAA$, if $A\geq0$ and $A\neq0$, then $\phi(A)>0$.
Note that by Proposition~\ref{PropositionSubStable}, for any
$C^*$-subalgebra $\AAA_0\subset\AAA$
and state $\phi$ on $\AAA$ the restriction of $\phi$ to $\AAA_0$ is
again a state.
(All of this is covered in~\cite{Murphy}, Chapter~3.)

\begin{Definition}
A pair $(\AAA,\phi)$ consisting of a $C^*$-algebra $\AAA$ and a
state $\phi$ is called a \textit{$C^*$-probability space}. We call $(\AAA
,\phi)$ \textit{faithful} 
if $\phi$ is so.
\end{Definition}

\subsubsection{Laws of noncommutative random variables}
Given a $C^*$-probability space $(\AAA,\phi)$
and self-adjoint $A\in\AAA$, there exists a unique
Borel probability measure $\mu_A$ on the spectrum of $A$,
called the \textit{law} of $A$, such that
$\phi(f(A))=\int f\,d\mu_A$ for every continuous $\CC$-valued function
$f$ on the
spectrum of $A$, where $f(A)$ is defined by means of the functional
calculus at~$A$,
that is, the inverse Gelfand transform,
and $\mu_A$ is provided by the Riesz representation theorem.
For convenience, we always extend the law $\mu_A$ to a Borel
probability measure on the real line
supported on $\Spec(A)$.
(See~\cite{AGZ}, Chapter~5, for background on laws.)
We note the important formula
%
%
\begin{equation}\label{equationKeyStieltjesRep}
S_{\mu_A}(z)=\phi\bigl((A-z1_\AAA)^{-1}\bigr)
\end{equation}
for the Stieltjes transform of the law $\mu_A$ which holds for every
$z\in\CC$ belonging neither to the support of $\mu_A$ nor to the
spectrum of $A$.
A simple and useful criterion for equality of the latter two sets is
provided by the next result.

\begin{Lemma}\label{LemmaFaithfulSupport}
Let $(\AAA,\phi)$ be a faithful $C^*$-probability space.
Then, for every $A\in\AAA_\sa$, $\supp\mu_A=\Spec( A)$.
\end{Lemma}

\begin{pf} Let $K=\Spec( A)\subset\RR$, noting that $K$ is compact.
Let $\AAA_0\subset\AAA$ be the $C^*$-subalgebra generated by $A$ and
put $\phi_0=\phi\vert_{\AAA_0}$,
which is a faithful state on $\AAA_0$.
By the theory of the Gelfand transform, $\AAA_0$ can be identified
with the $C^*$-algebra
of continuous complex-valued functions defined on~$K$. Under this
identification,
the operator $A$ becomes the
identity function $\Spec( A)\rightarrow\RR$ and $\phi_0$ becomes
the linear functional
represented by $\mu_A$. By Urysohn's lemma, $\phi_0$ cannot be
faithful unless $\supp\mu_A=K$.
\end{pf}

\begin{Lemma}\label{LemmaMatrixConvention}
If $(\AAA,\phi)$ is a faithful $C^*$-probability space,
then so is $(\Mat_n(\AAA),\break\phi_n)$, where $\phi_n(A)=\frac
{1}{n}\sum_{i=1}^n
\phi(A(i,i))$.
\end{Lemma}

We always follow the procedure of this lemma to equip $\Mat
_n(\AAA)$ with a state
when given a state on $\AAA$.\vadjust{\goodbreak}
\begin{pf*}{Proof of Lemma \protect\ref{LemmaMatrixConvention}} There is exactly one way to norm the $*$-algebra $\Mat
_n(\AAA)$ as a $C^*$-algebra.
(See~\cite{Murphy}, Theorem~3.4.2 and also Section \ref
{subsubsectionCstarTensorProduct} below.) Following our convention to
norm every $*$-algebra as a $C^*$-algebra when possible,
we thus regard $\Mat_n(\AAA)$ as a $C^*$-algebra.
For $0\neq A\in\Mat_n(\AAA)$ such that $A\geq0$,
\[
\phi_n(A)=\phi_n(A^{1/2}A^{1/2})=\frac{1}{n}\sum_{i,j=1}^n\phi
(A^{1/2}(i,j)A^{1/2}(i,j)^*).
\]
This formula first of all make it clear that $\phi_n$ is a state and hence
that $(\Mat_n(\AAA),\phi_n)$ is a $C^*$-probability space.
But furthermore, at least one term on the right is $>0$ since $\phi$
is faithful and $A^{1/2}\neq0$.
Thus, $\phi_n$ is faithful.
\end{pf*}

\begin{Lemma}\label{LemmaDykemaHelp}
Let $H$ be a Hilbert space, let $v\in H$ be a unit vector,
and consider the vectorial state $\phi=(A\mapsto
(v,Av))\dvtx B(H)\rightarrow\CC$ associated with $v$.
Let $\AAA,\hat{\AAA}\subset B(H)$ be $C^*$-subalgebras
such that $A\hat{A}=\hat{A}A$ for all $A\in\AAA$ and $\hat{A}\in
\hat{\AAA}$.
Assume furthermore that the vector $v$ is \textit{cyclic} for $\hat{\AAA
}$, that is, that the set $\{\hat{A}v\vert\hat{A}\in\hat{\AAA}\}$
is dense in $H$.
Then $\phi\vert_\AAA$ is faithful.
\end{Lemma}

\begin{pf}
Fix $A\in\AAA$ such that $A\geq0$ and $A\neq0$.
Clearly, there exists $h\in H$ such that $(h,A^{1/2}h)>0$.
Thus by hypothesis, there exists $\hat{A}\in\hat{\AAA}$ such that
\[
0< (\hat{A}v,A^{1/2}\hat{A}v)=(v,\hat{A}^*A^{1/2}\hat{A}v)=\phi
(\hat{A}^*A^{1/2}\hat{A})=\phi(\hat{A}^*\hat{A}A^{1/2}).
\]
Making further use of the hypothesis that operators in $\AAA$ commute
with operators in $\hat{\AAA}$, we have
\[
0\leq\phi\bigl(\bigl(\sqrt{t}\hat{A}^*\hat{A}-A^{1/2}/\sqrt{t}\bigr)^2\bigr)=t\phi
((\hat{A}^*\hat{A})^2)+\phi(A)/t-2\phi(\hat{A}^*\hat{A}A^{1/2})
\]
for $t>0$. The last inequality forces $\phi(A)>0$. Thus, $\phi\vert
_\AAA$ is indeed faithful.~%
\end{pf}

\subsection{Boltzmann--Fock space}
\label{subsectionBoltzmannFock}
We now recall the standard construction of free semicircular variables
laying stress on some properties of the construction which are not
often exploited
in practice, but which will be important here.
We are essentially just summarizing (using somewhat different notation)
enough material
from~\cite{VDN} to take advantage of~\cite{VDN}, Remark~2.6.7.

\subsubsection{\texorpdfstring{Definition of $\HHH$ and the $C^*$-probability space $(B(\HHH),\phi)$}
{Definition of H and the C*-probability space (B(H),phi)}}

Let $\HHH$ be a Hilbert space canonically equipped with an orthonormal basis
$\{v(i_1\cdots i_k)\}$
indexed by all finite sequences of positive integers, including the
empty sequence.
We write $1_\HHH=v(\varnothing)\in\HHH$.
We equip $B(\HHH)$ with the vectorial state $\phi^{\BF}$ defined by
\[
\phi^{\BF}(A)=(1_\HHH,A1_\HHH),
\]
thus making it into a noncommutative probability space. (Note that we
take Hilbert space
inner products to be linear on the right, antilinear on the left.)
Context permitting, we drop the superscript and write $\phi=\phi^{\BF}$.

\subsubsection{Raising and lowering operators}
Let $\Sigma_i\in B(\HHH)$
act by the rule
\[
\Sigma_iv(i_1\cdots i_k)=v(ii_1\cdots i_k).
\]
Let $p_\HHH\in B(\HHH)$ denote orthogonal projection to the linear
span of $1_\HHH$.
It is easy to verify the following relations, where $i$ and $j$ are any
positive integers:
%
%
\begin{eqnarray}\label{equationSigmaRelations}
 p_\HHH\Sigma_i&=&0 = \Sigma_i^*p_\HHH,\qquad
\Sigma_i^*\Sigma_j = \delta_{ij}1_{B(\HHH)},\qquad [\![\Sigma _i]\!] = [\![\Sigma_i^*]\!] = 1,\\
\label{equationSigmaMoments}
\phi(\Sigma_i)&=&\phi(\Sigma_i^*) = 0,\qquad
\phi(\Sigma_i\Sigma_j)=\phi(\Sigma_i\Sigma_j^*) = \phi(\Sigma
_i^*\Sigma_j^*) = 0,
\nonumber
\\[-8pt]
\\[-8pt]
\nonumber
\qquad\phi(\Sigma_i^*\Sigma_j) &=& \delta_{ij}.
\end{eqnarray}

\subsubsection{\texorpdfstring{The semicircular variables $\Xi_\ell$}
{The semicircular variables Xi l}}\label{subsubsectionXiConstruction}

We now make the special choice of sequence $\Xi$ mentioned in Remark
\ref{RemarkWillChooseSpecialXi} above.
Put
\[
\Xi_\ell=\ii^\ell\Sigma_\ell+\ii^{-\ell}\Sigma_\ell^*
\]
for all $\ell$ and put $\Xi=\{\Xi_\ell\}_{\ell=1}^\infty\in
B(\HHH)^\infty$.
It is well known that the joint law of the sequence $\{\Sigma_\ell
+\Sigma_\ell^*\}_{\ell=1}^\infty$
is free semicircular.
See~\cite{VDN} or~\cite{AGZ}, Chapter~5.
It is easy to see that the sequences $\Xi$ and $\{\Sigma_\ell+\Sigma
_\ell^*\}_{\ell=1}^\infty$ have the same joint law.
Indeed, the former is conjugate to the latter by a unitary
transformation preserving $\phi^\BF$.
We work exclusively with the choice of $\Xi$ made in this paragraph
for the rest of the paper.

\subsubsection{Right raising and lowering operators}\label
{subsubsectionRightAnalogues}
For each integer $i>0$, let $\hat{\Sigma}_i\in B(\HHH)$
act on $\HHH$ by the rule
\[
\hat{\Sigma}_iv(i_1\cdots i_k)=v(i_1\cdots i_ki).
\]
In direct analogy to \eqref{equationSigmaRelations},
we have
%
%
\begin{equation}\label{equationRhoRelations}
p_\HHH\hat{\Sigma}_i=0 = \hat{\Sigma}_i^*p_\HHH,\qquad
\hat{\Sigma}_i^*\hat{\Sigma}_j = \delta_{ij}1_{B(\HHH)},\qquad
[\![\hat{\Sigma}_i]\!] = [\![\hat{\Sigma}_i^*]\!] = 1.
\end{equation}
We also have right analogues of the relations
\eqref{equationSigmaMoments} but we will not need them.
It is easy to verify the following relations, where $i$ and $j$ are any
positive integers:
%
%
\begin{equation}\label{equationRaisingLoweringRelations}
\hat{\Sigma}_i\Sigma_j=\Sigma_j\hat{\Sigma}_i,\qquad \hat
{\Sigma}_j^*\Sigma_i=\Sigma_i\hat{\Sigma}_j^*+\delta_{ij}p_\HHH
,\qquad
\Sigma_i p_\HHH=\hat{\Sigma}_ip_\HHH.
\end{equation}
It is trivial but important to note that every relation above implies
another by taking adjoints on both sides.

\begin{Proposition}\label{PropositionPreMainResult}
For all $f\in\Mat_n({\CC\langle\Xbold\rangle})_{\sa}$, $\supp
\mu_f=\Spec(f(\Xi))$.
\end{Proposition}

\begin{pf} Let $\AAA\subset B(\HHH)$ be the $C^*$-subalgebra generated
by the sequence $\Xi$.
By Lemmas~\ref{LemmaFaithfulSupport} and~\ref{LemmaMatrixConvention},
it is enough to show that $\phi\vert_\AAA$ is faithful.
Let
\[
\hat{\Xi}_\ell=\ii^{\ell}\hat{\Sigma}_\ell+\ii^{-\ell}\hat
{\Sigma}_\ell^*\in B(\HHH)_\sa
\]
for positive integers $\ell$
and let $\hat{\AAA}\subset B(\HHH)$ be the $C^*$-subalgebra of
$B(\HHH)$
generated by the sequence $\{\hat{\Xi}_\ell\}$.
Using \eqref{equationRaisingLoweringRelations}, one verifies that
\[
\Xi_\ell\hat{\Xi}_m
=\hat{\Xi}_m\Xi_\ell\vadjust{\goodbreak}
\]
for all $\ell$ and $m$. Here, we are using the powerful insight of
\cite{VDN}, Remark~2.6.7.
Thus, every element of $\AAA$ commutes with every element of $\hat
{\AAA}$.
It is also easy to see that $1_\HHH$ is cyclic for $\hat{\AAA}$. Therefore,
$\phi\vert_\AAA$ is faithful by Lemma~\ref{LemmaDykemaHelp}.
\end{pf}

We conclude our discussion of Boltzmann--Fock space by recording the
following ``hypothesis-checking'' result for use in
Section~\ref{sectionOpTheoSchwingerDyson}.
\begin{Lemma}\label{LemmaQuasiCasimir}
Fix a positive integer $m$. Let
\[
x\in\bigl\{1_{B(\HHH)}\bigr\}\cup\{\Sigma_j,\Sigma_j^*\vert j=1,\ldots,m\}.
\]
The following hold:
\begin{itemize}
\item
$\hat{\Sigma}_i^*x\hat{\Sigma}_j=\delta_{ij}x$ for all $i$ and $j$.
\item$p_\HHH x\hat{\Sigma}_i=p_\HHH x\hat{\Sigma}_ip_\HHH$
and $p_\HHH\hat{\Sigma}_i^*xp_\HHH=\hat{\Sigma}_i^*xp_\HHH$
for all $i$.
\item$x$ commutes with $\hat{\Sigma}_i$ and $\hat{\Sigma}_i^*$
for all $i>m$.
\item$x$ commutes with $p_\HHH+\sum_{i=1}^m \hat{\Sigma}_i\hat
{\Sigma}_i^*$.
\end{itemize}
\end{Lemma}

\begin{pf} The first three statements follow straightforwardly
from \eqref{equationSigmaRelations}, \eqref{equationRhoRelations}
and \eqref{equationRaisingLoweringRelations}, so we just supply a
proof for the last statement.
We write $[A,B]=AB-BA$. Note that $[A,BC]=[A,B]C+B[A,C]$.
Fix $j\in\{1,\ldots,m\}$. We then have
\begin{eqnarray*}
 \Biggl[\Sigma_j,p_\HHH+\sum_{i=1}^m \hat{\Sigma}_i\hat
{\Sigma}_i^*\Biggr] &=& [\Sigma_j,p_\HHH]+\sum_{i=1}^m ([\Sigma_j,\hat
{\Sigma}_i]\hat{\Sigma}_i^*+\hat{\Sigma}_i[\Sigma_j,\hat{\Sigma
}_i^*])\\
&=& \Sigma_j p_\HHH+\hat{\Sigma}_j[\Sigma_j,\hat{\Sigma
}_j^*] = \Sigma_jp_\HHH-\hat{\Sigma}_jp_\HHH=0.
\end{eqnarray*}
The analogous relation with $\Sigma_j^*$ in place of $\Sigma_j$
follows by taking adjoints.~%
\end{pf}

\subsection{Projections, inverses and Schur complements}\label
{subsectionProjectionsAndInverses}
We make an ad hoc extension of the Schur complement concept to the
context of $C^*$-algebras.

\subsubsection{\texorpdfstring{Projections and $\pi$-inverses}
{Projections and pi-inverses}}\label{subsubsectionProjections}

Let $\AAA$ be a $C^*$-algebra.
A \textit{projection} $\pi\in\AAA$ by definition satisfies $\pi=\pi
^*=\pi^2$.
A family $\{\pi_i\}$ of projections is called \textit{orthonormal} if
$\pi_i\neq0$
and $\pi_i\pi_j=\delta_{ij}\pi_i$ for all $i$ and $j$.
Given $x\in\AAA$ and a projection $0\neq\pi\in\AAA$, we denote by
$x_\pi^{-1}$ the inverse of $\pi x\pi$ in the
$C^*$-algebra $\pi\AAA\pi$, if it exists, in which case it is
uniquely defined.
We call $x_\pi^{-1}$ the \textit{$\pi$-inverse} of $x$.
Note that $x_\pi^{-1}=(\pi x\pi)_\pi^{-1}$.

\begin{Proposition}\label{PropositionVeryDry}
Let $\AAA$ be a $C^*$-algebra.
Let $\{\pi,\pi^\perp\}$ be an orthonormal system of
projections in $\AAA$ and put $\sigma=\pi+\pi^\perp$.
Let $x\in\AAA$ satisfy $\pi^\perp x\pi^\perp\in(\pi^\perp\AAA
\pi^\perp)^{\times}$.
Then we have
%
%
\begin{equation}\label{equationPrePreVeryDry}
\sigma x\sigma\in(\sigma\AAA\sigma)^\times\Leftrightarrow
\pi(x-xx_{\pi^\perp}^{-1}x)\pi\in(\pi\AAA\pi)^\times
\end{equation}
and under these equivalent
conditions we have
%
%
\begin{eqnarray}
\label{equationPreVeryDry}
\pi x_\sigma^{-1}\pi&=&(x-x x_{\pi^\perp}^{-1} x)_\pi^{-1},\\
\label{equationVeryDry}
x^{-1}_\sigma-x_{\pi^\perp}^{-1}&=&
(\pi-x_{\pi^\perp}^{-1}x\pi) x_\sigma^{-1}(\pi-\pi xx_{\pi^\perp}^{-1}).
\end{eqnarray}
\end{Proposition}

Thus, the expression $\pi(x-xx_{\pi^\perp}^{-1}x)\pi$ is
a sort of generalized Schur complement.
\begin{pf*}{Proof of Proposition~\ref{PropositionVeryDry}} By exploiting the GNS construction,
one can interpret the proposition as an instance of \eqref
{equationFamiliarInversion}.
We omit the details.~%
\end{pf*}
%

\subsection{Cuntz frames and quasi-circularity}\label{subsectionSDAxiomatics}
We elaborate
upon a suggestion made in the last exercise of~\cite{AGZ}, Chapter~5.
We fix a $C^*$-algebra $\AAA$.
\begin{Definition}
Suppose we are given a collection
$\{\pi\}\cup\{\rho_i\}_{i=1}^\infty$
of elements of $\AAA$ satisfying the following conditions:
%
%
\begin{equation}
\label{equationSD05}\qquad
\mbox{$\pi$ is a nonzero projection, $\pi\rho_i=0$ and $\rho
_{i}^{*}\rho_j=\delta_{ij}1_{\AAA}$} \mbox{ for all $i$ and $j$.}
\end{equation}
We call $\{\pi\}\cup\{\rho_i\}_{i=1}^\infty$ a \textit{Cuntz frame} in
$\AAA$.
Note that $\{\pi\}\cup\{\rho_i\rho_i^*\}_{i=1}^\infty$ is
automatically an orthonormal system of projections.
\end{Definition}

\begin{Remark}
The relations $\rho_i^*\rho_j=\delta_{ij}1_\AAA$ are those defining
the \textit{Cuntz algebra}~\cite{Cuntz}, hence our choice of terminology.
\end{Remark}

\subsubsection{Quasi-circular operators}
Suppose we are given a Cuntz frame
$\{\pi\}\cup\{\rho_i\}_{i=1}^\infty$ in $\AAA$.
We say that an operator $A\in\AAA$ is \textit{quasi-circular} (with
respect to the given Cuntz frame)
if the following statements hold:
%
%
\begin{eqnarray}
\label{equationSD2}
&{\rho}_i^*A{\rho}_j=\delta_{ij} A\qquad \mbox{for all $i$ and $j$,}&\\
\label{equationSD3}
&{\pi} A{\rho}_i{\pi}={\pi} A{\rho}_i\quad \mbox{and}\quad {\pi}
{\rho}_i^*A{\pi}={\rho}_i^*A{\pi}\qquad \mbox{for all $i$,}&
\end{eqnarray}\vspace*{-16pt}
\begin{equation}\label{equationSD1}
\begin{tabular}{p{260pt}@{}}
\mbox{There exists an integer $k_A\geq0$ such
that $A$ commutes}
\mbox{with $\pi+\sum_{i=1}^{k_A} {\rho}_i{\rho}_i^*$
and also with ${\rho}_i$ and ${\rho}_i^*$ for all $i>k_A$.}
\end{tabular}
\end{equation}

\begin{Proposition}\label{PropositionAbstractSchwingerDyson}
Let $\{\pi\}\cup\{\rho_i\}_{i=1}^\infty$ be a Cuntz frame in $\AAA$.
Let $A\in\AAA^\times$ be quasi-circular with respect to the given frame.
Choose any integer $k\geq k_A$. Then
%
%
\begin{equation}
\label{equationSDMachine3}
\pi A^{-1}\pi= \Biggl(\pi A\pi-
\sum_{i=1}^k\pi A\rho_i\pi A^{-1}\pi\rho_i^*A\pi\Biggr)^{-1}_\pi.
\end{equation}
In particular, one automatically has $\pi A^{-1}\pi\in(\pi\AAA\pi
)^\times$.
\end{Proposition}

Identity \eqref{equationSDMachine3} is an abstract
algebraic version of the Schwinger--Dyson equation.
See the proof of Proposition~\ref{PropositionConcreteSchwingerDyson}
below for the application.\vadjust{\goodbreak}

\begin{pf*}{Proof of Proposition~\ref{PropositionAbstractSchwingerDyson}} Consider the projections
$\sigma=\pi+\sum_{i=1}^k\rho_i\rho_i^*$
and $\pi^\perp=\sigma-\pi$.
We claim that
%
%
\begin{equation}
\label{equationSDMachine2}
A_{\pi^\perp}^{-1}=\sum_{i=1}^k\rho_i A^{-1}\rho_i^*.
\end{equation}
In any case, we have
$\pi^\perp A\pi^\perp=\sum_{i=1}^k \rho_i A\rho_i^*$
by \eqref{equationSD2}. Furthermore, we have
\[
\Biggl(\sum_{i=1}^k \rho_i A\rho_i^*\Biggr)
\Biggl(\sum_{j=1}^k \rho_j A^{-1}\rho_j^*\Biggr)=
\pi^\perp=\Biggl(\sum_{i=1}^k \rho_i A^{-1}\rho_i^*\Biggr)
\Biggl(\sum_{j=1}^k \rho_j A\rho_j^*\Biggr)
\]
by \eqref{equationSD05}. Thus, claim \eqref{equationSDMachine2} holds.
To prove \eqref{equationSDMachine3}, we calculate
as follows:
\begin{eqnarray*}
\pi A^{-1}\pi&=&\pi\sigma A^{-1}\sigma\pi= \pi A_\sigma^{-1}\pi
= (A-AA_{\pi^\perp}^{-1}A)^{-1}_\pi\\
&=&
(\pi A\pi-\pi AA_{\pi^\perp}^{-1}A\pi)^{-1}_\pi
=\Biggl(\pi A\pi-\sum_{i=1}^k \pi A\rho_iA^{-1}\rho_i^*A\pi\Biggr)_\pi^{-1}\\
&=& \Biggl(\pi A\pi-\sum_{i=1}^k \pi A\rho_i\pi A^{-1}\pi\rho
_i^*A\pi\Biggr)_\pi^{-1}.
\end{eqnarray*}
The first step is simply an exploitation of orthonormality of $\{\pi
,\pi^\perp\}$.
Since $A$ commutes with $\sigma$ by \eqref{equationSD1}, we have
$\sigma A^{-1}\sigma A=\sigma=1_{\sigma A\sigma}=A\sigma A^{-1}\sigma$,
which justifies the second step.
The third step is
an application of \eqref{equationPreVeryDry} and
the fourth step is a trivial consequence of the definition of $\pi$-inverse.
The fifth step is an application of \eqref{equationSDMachine2}
and the last step is an application of \eqref{equationSD3}.
The proof of \eqref{equationSDMachine3} is complete.
\end{pf*}

\begin{Remark}
The preceding calculation will obviate consideration of combinatorics of
free semicircular variables in the sequel. We present this approach as
counterpoint to the nowadays standard combinatorial approach
discussed briefly in~\cite{AGZ}, Chapter~5 and developed at length in
\cite{NicaSpeicher}.
\end{Remark}

\section{Tools for concentration}
\label{sectionConcentration}
In this section, we introduce an ensemble of tools we use to
(partially) replace
the Poincar\'{e}-type inequalities
used in~\cite{HT,HST,Male,CapDon} and~\cite{Schultz}.
We speak of an ensemble because no one tool seems to contribute more
than incrementally.

\subsection{Tensor products and norming rules} \label
{subsectionTensorProducts}
We rehearse the most basic rules of calculation and estimation
used in the paper.

\subsubsection{Tensor products of vector spaces and algebras}
Given vector spaces $\AAA$ and $\BBB$ over $\CC$,
let $\AAA\otimes\BBB$ denote the tensor product of $\AAA$ and $\BBB
$ formed over~$\CC$. If $\AAA$ and $\BBB$ are both algebras, we
invariably endow $\AAA\otimes\BBB$ with algebra structure by the\vadjust{\goodbreak}
rule $(a_1\otimes b_1)(a_2\otimes b_2)=
a_1a_2\otimes b_1b_2$.
If $\AAA$ and $\BBB$ are both $*$-algebras, we invariably endow $\AAA
\otimes\BBB$
with $*$-algebra structure by the rule $(a\otimes b)^*=a^*\otimes b^*$.

\subsubsection{Tensor notation for building matrices}
Let $\AAA$ be an algebra.
We identify the algebra $\Mat_n(\CC)\otimes\AAA$ with $\Mat_n(\AAA)$
by the rule $(X\otimes a)(i,j)=x(i,j)a$ and more generally use the same
rule to identify the space\break $\Mat_{k\times\ell}(\CC)\otimes\AAA$
with the space of rectangular matrices $\Mat_{k\times\ell}(\AAA)$.
Furthermore, in the case $\AAA=\Mat_s(\CC)$,
we identify $X\otimes a$ with an element of $\Mat_{ks\times\ell
s}(\CC)$
by viewing $X\otimes a$ as a $k$-by-$\ell$
arrangement of $s$-by-$s$ blocks $X(i,j)a$.
In other words, we identify $X\otimes a$ with the usual Kronecker
product of $X$ and $a$.

\subsubsection{Banach spaces}
Banach spaces always have $\CC$ as scalar field,
and bounded (multi)linear maps between Banach spaces
are always $\CC$-(multi)li\-near, unless explicitly noted otherwise.
To avoid collision with the notation ${\Vert\cdot \Vert}_p$,
we let $[\![\cdot]\!]_\VVV$ denote the norm of a Banach space $\VVV$
and context permitting (nearly always), we drop the subscript.

\subsubsection{(Multi)linear maps between Banach spaces}
Given Banach spaces $\VVV$ and $\WWW$, let $B(\VVV,\WWW)$ denote
the space of bounded linear maps $\VVV\rightarrow\WWW$.
Let $B(\VVV)=B(\VVV,\VVV)$ and let $\VVV^\star$ denote the linear
dual of $\VVV$.
Given $T\in B(\VVV,\WWW)$,
let $[\![T]\!]=[\![T]\!]_{B(\VVV,\WWW)}$ be the best constant
such that
$[\![Tv]\!]\leq[\![T]\!][\![v]\!]$. We always use the norm on
$B(\VVV,\WWW)$ so defined.
More generally, let $B(\VVV_1,\ldots,\VVV_r;\WWW)$ denote the space
of bounded
$r$-linear maps $\VVV_1\times\cdots\times\VVV_r\rightarrow\WWW$
and given $T\in B(\VVV_1,\ldots,\VVV_r;\WWW)$,
let
$[\![T]\!]=[\![T]\!]_{B(\VVV_1,\ldots,\VVV_r;\WWW)}$ be the best
constant such that
$[\![T(v_1,\ldots,v_r)]\!]\leq[\![T]\!][\![v_1]\!]\cdots[\![v_r]\!]$.
We always use the norm on $B(\VVV_1,\ldots,\VVV_r; \WWW)$ so defined.

\subsubsection{Matrix spaces over $C^*$-algebras}\label
{subsubsectionNormingRules}
Let $\AAA$ be any $C^*$-algebra. We have already noted
in the proof of Lemma~\ref{LemmaMatrixConvention} that there is a
unique way to norm
the $*$-algebra $\Mat_n(\AAA)$ as a $C^*$-algebra. In turn, we always
norm the space of rectangular matrices $\Mat_{k\times\ell}(\AAA)$
by the formula $[\![A]\!]=[\![AA^*]\!]^{1/2}$. Note that
%
%
\begin{equation}\label{equationBlockEstimate}
\bigvee_{i=1}^k\bigvee_{j=1}^\ell[\![A(i,j)]\!]\leq[\![A]\!]\leq
\sum_{m=-\infty}^\infty\bigvee_{i=1}^k\bigvee_{j=1}^\ell[\![A(i,j)]\!]\one_{i-j=m}.
\end{equation}
Moreover, given $B\in\Mat_{\ell\times m}(\AAA)$, we have $[\![AB]\!]\leq[\![A]\!][\![B]\!]$.
In particular, for every square or rectangular matrix $A$
with complex number entries,
$[\![A]\!]$ is the largest singular value of $A$.

\subsubsection{Tensor products of $C^*$-algebras}\label
{subsubsectionCstarTensorProduct}
Given $C^*$-algebras $\AAA$ and $\BBB$ with at least one of them
finite-dimensional,
the $*$-algebra $\AAA\otimes\BBB$ has exactly one $C^*$-algebra norm.
To see this, only existence requires comment since uniqueness we have\vadjust{\goodbreak}
already noted
after Proposition~\ref{PropositionCstarNormUnique}.
We proceed as follows. Firstly, we observe that since $\AAA\otimes
\BBB$ and $\BBB\otimes\AAA$ are isomorphic $*$-algebras,
we may assume that $\AAA$ is finite-dimensional.
Then, after reducing to the case $\AAA=\Mat_n(\CC)$ and $\BBB=B(H)$
by using the GNS construction, we can make identifications $\AAA
\otimes\BBB=\Mat_n(B(H))=B(H^n)$ yielding the desired norm. Thus,
existence is settled. The preceding argument shows that
for all $a\in\AAA$ and $b\in\BBB$ we have $[\![a\otimes b]\!]=[\![a]\!][\![b]\!]$.
In a similar vein, we have the following useful general observation.

\begin{Lemma}\label{LemmaBlockization}
Let $\SSS$ be a finite-dimensional $C^*$-algebra.
Let $\{e_i\}_{i=1}^n$ be any linearly independent family of elements of
$\SSS$.
Then for all $C^*$-algebras $\AAA$ and families $\{a_i\}_{i=1}^n$ of
elements of $\AAA$ we have
\[
\frac{1}{C}\bigvee_{i=1}^n [\![a_i]\!]\leq\Biggl[\!\Biggl[\sum_{i=1}^n e_i\otimes
a_i\Biggr]\!\Biggr]
=\Biggl[\!\Biggl[\sum_{i=1}^n a_i\otimes e_i\Biggr]\!\Biggr]
\leq C\sum_{i=1}^n [\![a_i]\!]
\]
for a constant $C\geq1$ depending only on $\SSS$
and $\{e_i\}$.
\end{Lemma}
\begin{pf} We may assume that $\SSS=\Mat_s(\CC)$.
Furthermore, there is no loss of generality to assume that $n=s^2$ and thus
that \vspace*{-1pt}$\{e_i\}_{i=1}^{s^2}$ is a basis for $\Mat_s(\CC)$.
Finally, there is no loss of generality
to assume that $\{e_i\}_{i=1}^{s^2}$ consists of elementary matrices,
in which case the lemma at hand reduces to \eqref{equationBlockEstimate}.
\end{pf}

\subsubsection{Block algebras}
We define a \textit{block algebra} to be a $C^*$-algebra
isomorphic to $\Mat_s(\CC)$ for some positive integer $s$.
The point of the definition is merely to compress notation and to put
some psychological distance between
us and the too-familiar algebra $\Mat_s(\CC)$.
(Later we will refine the notion of block algebra by keeping track not
only of the transpose conjugate operation
but also of the transpose.
See Section~\ref{subsectionBlockAlgebraUpgrade} below.)

\subsection{Quadratic forms in independent random vectors}
Variants of the next result
are in common use in RMT. (See, e.g.,~\cite{BaiSil}, Lemma 2.7.)

\begin{Proposition}\label{PropositionWhittle}
Let $Y_1,\ldots,Y_n$ and $Z_1,\ldots,Z_n$ be $\CC$-valued random variables
which for some $p\in[2,\infty)$ all belong to $L^{2p}$ and have mean zero.
Let $A\in\Mat_n(\CC)$ be a (deterministic) matrix.
Assume furthermore that the family of $\sigma$-fields $\{
\sigma(Y_i,Z_i)\}_{i=1}^n$ is independent.
Then we have
\[
{\Biggl\Vert\sum_{i,j=1}^n A(i,j)(Y_iZ_j-\Ebold Y_iZ_j) \Biggr\Vert}_p
\leq c \Biggl(\sum_{i,j=1}^n|A(i,j)|^2{\Vert Y_i \Vert}_{2p}^2{\Vert Z_j
\Vert}_{2p}^2\Biggr)^{1/2}
\]
for a constant $c$ depending only on $p$.
\end{Proposition}

Results of this type are well-known.
The earliest reference we know is~\cite{Whittle}.
In that reference, the result above is proved
in the special case in which\vadjust{\goodbreak} $Y_i=Z_i=Y_i^*=Z_i^*$ and $A$ has real entries.
From that special case, the general case
of the proposition above can be deduced easily by algebraic manipulation.

We now generalize in an innocuous if superficially complicated way.

\begin{Proposition}\label{PropositionKeyEstimate}
Fix constants $p\in[2,\infty)$ and $K\in(0,\infty)$.
Let $\VVV$ be a finite-dimensional Banach space,
let $\SSS$ be a block algebra and let $\GGG$ be a $\sigma$-field.
Let $Y\in\Mat_{1\times n}(\SSS)$ and $Z\in\Mat_{n\times1}(\SSS)$
be random such that
\[
\Bigl(\bigvee{\Vert[\![Y(1,j)]\!] \Vert}_{2p}\Bigr)\vee
\Bigl(\bigvee{\Vert[\![Z(i,1)]\!] \Vert}_{2p}\Bigr)\leq K\quad \mbox{and}\quad \Ebold
Y=0=\Ebold Z.
\]
Assume also that the family
$\GGG\cup\{\sigma(Y(1,i),Z(i,1))\}_{i=1}^n$
of $\sigma$-fields is independent.
Then for any $\GGG$-measurable random bilinear map
\[
R\in B(\Mat_{1\times n}(\SSS),\Mat_{n\times1}(\SSS);\VVV)
\]
such that ${\Vert[\![R]\!] \Vert}_p<\infty$
we have
\[
{\Vert[\![R(Y,Z)-\Ebold(R(Y,Z)\vert\GGG)]\!] \Vert}_p\leq CK^2{\Vert
[\![R]\!] \Vert}_p \sqrt{n},
\]
where the constant $C$ depends only on $p$, $\SSS$ and $\VVV$.
\end{Proposition}

We need two lemmas, the first of which actually proves more
than we immediately need but has several further uses in the paper.

\begin{Lemma}\label{LemmaHilbertSchmidt}
Let $\SSS$ be a block algebra of dimension $s^2$.
\textup{(i)} For $X\in\Mat_{k\times\ell}(\SSS)$, we have
\[
\frac{1}{s}[\![X]\!]^2\leq\sum_{i=1}^k\sum_{j=1}^\ell
[\![X(i,j)]\!]^2\leq s(k\wedge\ell)[\![X]\!]^2.\vspace*{-6pt}
\]
\begin{longlist}[(iii)]
\item[(ii)] For $X\in\Mat_{k\times\ell}(\SSS)$ and $Y\in\Mat_{\ell
\times k}(\SSS)$, we have
\[[\![\trace_\SSS( XY)]\!]\leq s\ell[\![X]\!][\![Y]\!].
\]
\item[(iii)] For $X\in\Mat_n(\SSS)$, we have
\[\Biggl[\!\Biggl[\sum_{i,j=1}^n X(i,j)^{\otimes2}\Biggr]\!\Biggr]\leq sn[\![X]\!]^2.
\]
\end{longlist}
\end{Lemma}

\begin{pf} Statement (i) is an assertion concerning the Hilbert--Schmidt
norm which is easy
to verify. Statements (ii) and (iii) follow from statement (i) via the
Cauchy--Schwarz inequality.
\end{pf}

\begin{Lemma}\label{LemmaHorribleIndices}
Let $\SSS$ be a block algebra
and let $\{e_i\}_{i=1}^\ell$ be a basis of the underlying vector space.
Let $\VVV$ be a finite-dimensional Banach space
and let $\{v_k\}_{k=1}^m$ be a basis of the underlying vector space.
Fix matrices $R_{ij}^k\in\Mat_n(\CC)$\vadjust{\goodbreak}
for $i,j=1,\ldots,\ell$ and $k=1,\ldots,m$.
Define $R\in B(\Mat_{1\times n}(\SSS),\break\Mat_{n\times1}(\SSS);\VVV
)$ by
requiring that
\[
R(x\otimes e_i,
y\otimes e_j)=\sum_k
(xR_{ij}^ky)v_k
\]
for $i,j=1,\ldots,\ell$,
$x\in\Mat_{1\times n}(\CC)$ and $y\in\Mat_{n\times1}(\CC)$.
Then
\[
\frac{1}{C}\bigvee_{i,j,k}[\![R_{ij}^k]\!]\leq
[\![R]\!]\leq C\sum_{i,j,k}[\![R_{ij}^k]\!]
\]
for a constant $C\geq1$ which depends only on the data
$(\SSS,\{e_i\},\VVV,\{v_k\})$
and in particular is independent of $n$.
\end{Lemma}

\begin{pf} By Lemma~\ref{LemmaBlockization} and the fact that the map
\[
\bigl(A\mapsto\bigl((x,y)\mapsto xAy\bigr)\bigr)\dvtx \Mat_n(\CC)\rightarrow B(\Mat_{1\times
n}(\CC),
\Mat_{n\times1}(\CC);\CC)
\]
is an isometric isomorphism,
the proof of the lemma at hand reduces to a straightforward calculation
the remaining details of which we can safely omit.
\end{pf}

\subsubsection{\texorpdfstring{Proof of Proposition \protect\ref{PropositionKeyEstimate}}
{Proof of Proposition 11}}
After using standard properties of conditional expectation,
we may assume that $R$ is deterministic.
We may also assume that $\SSS$ is isomorphic to $\Mat_s(\CC)$ for
some $s$
and in turn Lemma~\ref{LemmaHorribleIndices} permits us to assume
that $\SSS=\CC$.
Finally, by Lemma~\ref{LemmaHilbertSchmidt},
the proposition at hand reduces
to Proposition~\ref{PropositionWhittle}.

\begin{Remark}\label{RemarkTwoTypesBilinear}
In applications of Proposition~\ref{PropositionKeyEstimate},
we will only use two special types of bilinear map $R$.
We describe these types and estimate $[\![R]\!]$ for each.
(They conform to the patterns set\vspace*{1pt} by the objects $Q^N_{I,J,j_1,j_2}$
and $P^N_{I,J,j_1,j_2}$ defined in Section \ref
{sectionUrMatrixIdentities} below, resp.)
(i) In the ``$Q$-type''
first case of interest, we have $\VVV=\SSS$ and
for some $A\in\Mat_n(\SSS)$ we have $R(y,z)=yAz$, in which case
$[\![R]\!]\leq[\![A]\!]$.
(ii) In the ``$P$-type''
second case of interest,
we have $\VVV=B(\SSS)$, and for some $A\in\Mat_n(\SSS)$
we have $R(y,z)=(B\mapsto\trace_\SSS(AzByA))$,
in which case $[\![R]\!]\leq s[\![A]\!]^2$ for $s$ equal to the
square root
of the dimension of $\SSS$ over the complex numbers
by Lemma~\ref{LemmaHilbertSchmidt}.
\end{Remark}

\subsection{A conditional variance bound}
We present a result which harmlessly generalizes
the well-known subadditivity of variance
to a situation involving vector-valued random variables
and some mild dependence.
\subsubsection{Setup for the result}
Let $\VVV$ be a finite-dimensional Banach space (either real or
complex scalars).
Let
$\{\EEE\}\cup\{\GGG(i,j)\}_{1\leq i\leq j\leq N}$ be a family of independent
$\sigma$-fields and
let $\GGG$ be the $\sigma$-field generated by this family.
Let $Z\in\VVV$ be a $\GGG$-measurable random vector such that
${\Vert[\![Z]\!] \Vert}_p<\infty$ for $p\in[1,\infty)$.
For $k=1,\ldots,N$, let $\widehat{\GGG}_k$
be the $\sigma$-field generated\vadjust{\goodbreak} by the subfamily
$\{\EEE\}\cup\{\GGG(i,j)\vert k\notin\{i,j\}\}$
and let $Z_k\in\VVV$ be a $\widehat{\GGG}_k$-measurable random vector
such that ${\Vert[\![Z_k]\!] \Vert}_p<\infty$ for
$p\in[1,\infty)$.
\begin{Proposition}\label{PropositionVarianceBound}
Notation and assumptions are as above. For every constant $p\in
[1,\infty)$, we have
%
%
\begin{equation}\label{equationVarianceBound}
{\bigl\Vert\Ebold\bigl([\![Z-\Ebold(Z\vert\EEE)]\!]^2\vert\EEE\bigr)\bigr \Vert
}_p\leq
c\sum_{k=1}^N {\Vert[\![Z-Z_k]\!] \Vert}_{2p}^2
\end{equation}
for a constant $c$ depending only on $\VVV$ and in particular
independent of $p$.
\end{Proposition}

\begin{pf} We may assume that $\VVV$ is a (finite-dimensional) real
Hilbert space, and in this case we will prove the claim
with a constant $c=1$. After a routine application of Minkowski and Jensen
inequalities, it is enough to prove
%
%
\begin{equation}\label{equationNufVB1}
\Ebold\bigl([\![Z-\Ebold(Z\vert\EEE)]\!]^2\vert\EEE\bigr)\leq
\sum_{k=1}^N\Ebold([\![Z-Z_k]\!]^2\vert\EEE),
\end{equation}
almost surely. There is also no harm in assuming that $\VVV=\RR$.
For $k=0,\ldots,N$, let $\GGG_k$ be the $\sigma$-field generated
by the subfamily
$\{\EEE\}\cup\{\GGG(i,j)\vert1\leq i\leq j\leq k\}$.
In any case, by orthogonality of martingale increments, we have
\[
\Ebold\bigl([\![Z-\Ebold(Z\vert\EEE)]\!]^2\vert\EEE\bigr)=
\sum_{i=1}^N \Ebold\bigl([\![\Ebold(Z\vert\GGG_k)-\Ebold(Z\vert \GGG
_{k-1})]\!]^2\vert\EEE\bigr),
\]
almost surely. Furthermore,
we have
\[
\Ebold(\Ebold(Z\vert\widehat{\GGG}_k)\vert\GGG_k)=
\Ebold(Z\vert\GGG_{k-1}),
\]
almost surely.
Finally, we have
\begin{eqnarray*}
\Ebold\bigl([\![\Ebold(Z\vert\GGG_k)-\Ebold(Z\vert\GGG _{k-1})]\!]^2\vert\EEE\bigr)
&=&
\Ebold\bigl(\bigl[\!\bigl[\Ebold\bigl(Z-\Ebold(Z\vert\widehat{\GGG }_k)\vert\GGG
_{k}\bigr)\bigr]\!\bigr]^2\big\vert\EEE\bigr)\\
&\leq& \Ebold\bigl([\![Z-\Ebold(Z\vert\widehat{\GGG }_k)]\!]^2\vert\EEE\bigr)
\leq\Ebold([\![Z-Z_k]\!]^2\vert\EEE),
\end{eqnarray*}
almost surely, whence \eqref{equationNufVB1}.
\end{pf}
\begin{Definition}
The random variable $\Ebold([\![Z-\Ebold(Z\vert\EEE)]\!]^2\vert
\EEE)$ appearing on the
left-hand side of \eqref{equationVarianceBound} will be denoted by $\Var
_\VVV(Z\vert\EEE)$ in the sequel.
\end{Definition}

\subsection{Estimates for tensor-cubic forms}
We work out a specialized estimate involving three-fold tensor products
and partitions of a set of cardinality six.
The combinatorial apparatus introduced here will have further uses.

\subsubsection{Set partitions and related apparatus}
A \textit{set partition} of $k$ is a disjoint family $\Pi$ of nonempty
subsets of the set $\{1,\ldots,k\}$ whose union is $\{1,\ldots,k\}$.
Each member of a set partition is called a \textit{part}.\vadjust{\goodbreak}
Let $\Part(k)$ be the family of set partitions of $k$.
Let $\Part^*(2k)$ be the subset of $\Part(2k)$ consisting
of set partitions having no singleton as a part,
nor having any of the sets $\{2i-1,2i\}$ for $i=1,\ldots,k$ as a part.
Let $\Part^*_2(2k)\subset\Part^*(2k)$ be the subfamily
consisting of partitions all of whose parts have cardinality $2$.
For each positive integer $k$ let $S_k$ be the group of permutations of
$\{1,\ldots,k\}$.
Let $\Gamma_k\subset S_{2k}$ be the subgroup
centralizing the involutive permutation $(12)\cdots(2k-1,2k)$.
Then $\Gamma_k$
acts on the set $\Part^*(2k)$. For $\Pi_1,\Pi_2\in\Part^*(2k)$
belonging to the same $\Gamma_k$-orbit we write $\Pi_1\sim\Pi_2$.

\subsubsection{Explicit descriptions of $\Part^*(4)$ and $\Part^*(6)$}
To describe $\Part^*(4)$, we can easily enumerate it, thus:
%
%
\begin{equation}\label{equationPartitionEnumeration}
\{\{1,2,3,4\}\},\qquad\{\{1,3\},\{2,4\}\},\qquad\{\{1,4\},\{2,3\}\}.
\end{equation}
It can be shown (we omit the tedious details) that for every $\Pi\in
\Part^*(6)$
there exists exactly one set partition on the list
%
%
\begin{eqnarray}\label{equationPartitionList}
&&\{\{1,2,3,4,5,6\}\},\qquad \{\{1,6\},\{2,3,4,5\}\}, \qquad\{\{1,3,5\},\{2,4,6\}
\},
\nonumber
\\[-8pt]
\\[-8pt]
\nonumber
&&\qquad\{\{1,6\},\{2,3\},\{4,5\}\},\qquad \{\{1,2,3\},\{4,5,6\}\}
\end{eqnarray}
belonging to the $\Gamma_3$-orbit of $\Pi$.
\subsubsection{Sequences and associated partitions}
For any finite set $I$, we write
\[
\Seq(k,I)=\bigl\{{\mathbf{i}}\dvtx \{1,\ldots,k\}\rightarrow I\bigr\}.
\]
Given ${\mathbf{i}}\in\Seq(k,I)$,
let $\Pi({\mathbf{i}})\in\Part(k)$ be the set partition generated
by ${\mathbf{i}}$,
that is, the coarsest set partition on the parts of which ${\mathbf
{i}}$ is constant.
If $I=\{1,\ldots,n\}$, we write $\Seq(k,I)=\Seq(k,n)$ by abuse of notation.
Sometimes we represent elements of $\Seq(k,I)$ as ``words'' $i_1\cdots i_k$
spelled with ``letters'' $i_1,\ldots, i_k\in I$.

\subsubsection{Setup for the main result}
Let $\SSS$ be a block algebra.
Let a set partition $\Pi\in\Part^*(6)$ and matrices $M_1,M_2,M_3\in
\Mat_n(\SSS)$ be given.
Put
\[
\Mfrak_\Pi=
\cases{
\displaystyle\Biggl[\!\Biggl[\mathop{\mathop{\sum_{{\mathbf{i}}=i_1\cdots i_6}}_{\in\Seq(6,n)}}_{\mathrm{s.t.\ }
\Pi({\mathbf{i}})=\Pi}M_1(i_1,i_2)\otimes M_2(i_3,i_4)\otimes M_3(i_5,i_6)\Biggr]\!\Biggr],&\quad
$\mbox{if $\Pi\in\Part^*_2(6)$},$\vspace*{2pt}\cr
\displaystyle\mathop{\mathop{\sum_{{\mathbf{i}}=i_1\cdots i_6}}_{\in\Seq(6,n)}}_{\mathrm{s.t.\ }
\Pi({\mathbf{i}})=\Pi}[\![M_1(i_1,i_2)]\!][\![M_2(i_3,i_4)]\!][\![M_3(i_5,i_6)]\!],&\quad
$\mbox{if $\Pi\notin\Part^*_2(6)$.}$}
\]
\begin{Proposition}\label{PropositionCubic} Notation and assumptions
are as above.
For $\Pi\in\Part^*(6)$,
unless
$\Pi\sim\{\{1,2,3\},\{4,5,6\}\}$,
we have
$\Mfrak_\Pi\leq c n[\![M_1]\!][\![M_2]\!][\![M_3]\!]$
for a constant $c$ depending only on $\SSS$.\vadjust{\goodbreak}
\end{Proposition}
\begin{pf} We may assume that $\SSS$ is isomorphic to $\Mat_s(\CC)$
and thus by Lemma~\ref{LemmaBlockization} that $\SSS=\CC$.
After replacing $(M_1,M_2,M_3)$
by $(M_{\sigma(1)}^{\T^{\nu_1}},M^{\T^{\nu_2}}_{\sigma(2)}
,M_{\sigma(3)}^{\T^{\nu_3}})$
for suitably chosen $\sigma\in S_3$ and $\nu_1,\nu_2,\nu_3\in\{
0,1\}$,
we may assume that $\Pi$ appears on the list \eqref{equationPartitionList}.
We may also assume that each matrix $M_\alpha$
is either diagonal or else vanishes identically on the diagonal.
Finally, we may assume that $\Mfrak_\Pi>0$.
Let $d$ be the number of matrices $M_\alpha$ which are diagonal.
Consider the following mutually exclusive and exhaustive collection of cases:
\begin{longlist}[(iii)]
\item[(i)] $\Pi= \{\{1,6\},\{2,3\},\{4,5\}\}$ and hence
$d=0$.
\item[(ii)] $\Pi= \{\{1,3,5\},\{2,4,6\}\}$
and hence $d=0$.
\item[(iii)] $\Pi= \{\{1,6\},\{2,3,4,5\}\}$ and hence $d=1$.
\item[(iv)] $\Pi=\{\{1,2,3,4,5,6\}\}$ and hence $d=3$.
\end{longlist}
In case (i) we have
$\Mfrak_\Pi=|\trace M_1M_2M_3|\leq n[\![M_1M_2M_3]\!]
\leq n[\![M_1]\!][\![M_2]\!][\![M_3]\!]$.
In case (ii) we have
\begin{eqnarray*}
\Mfrak_\Pi&\leq& [\![M_{1}]\!]\sum_{i,j=1}^n
|M_{2}(i,j)M_{3}(i,j)|\\
&\leq& [\![M_{1}]\!]\prod_{\alpha\in\{2,3\}}\Biggl(\sum
_{i,j=1}^n |M_{\alpha}(i,j)|^2\Biggr)^{1/2} \leq n[\![M_1]\!][\![M_2]\!][\![M_3]\!].
\end{eqnarray*}
In case (iii), similarly, we have
\[
\Mfrak_\Pi\leq[\![M_2]\!]\sum_{i,j=1}^n
|M_1(i,j)M_3(j,i)|\leq n[\![M_1]\!][\![M_2]\!][\![M_3]\!].
\]
Finally, in case (iv) we have $\Mfrak_\Pi\leq n[\![M_1]\!][\![M_2]\!][\![M_3]\!]$
simply by counting.
\end{pf}

\section{Transpositions, SALT block designs and the secondary trick}
\label{sectionUpgrade}
We introduce algebraic tools
for exploiting the symmetry which the random matrices \eqref
{equationTypicalBlockRandom}
possess
as a consequence of assumption \eqref{equationPepper5}
and which the operators \eqref{equationTypicalBlockFock} also possess
by virtue
of the special choice of $\Xi$ made in Section~\ref{subsectionBoltzmannFock}.
We introduce SALT block designs, we show how the self-adjoint
linearization trick
generates examples of such, and we use SALT block designs to put the
crucial equations~\eqref{equationEmpiricalAccess} and \eqref{equationPreStieltjesAccess}
into streamlined form. See Remark~\ref{RemarkPreCobbling} below. We
introduce the secondary trick
which produces new SALT block designs from old and in particular
produces some SALT block designs
which do not come from the self-adjoint linearization trick.

\subsection{Transpositions}

\begin{Definition}\label{DefinitionTransposition}
Let $\AAA$ be a $*$-algebra. A \textit{transposition} $a\mapsto a^\T$ of
$\AAA$ is
a $\CC$-linear map such that $(a^\T)^\T=a$, $(a^*)^\T
=(a^\T)^*$ and $(ab)^\T=b^\T a^\T$ for all $a,b\in\AAA$.
Necessarily $1_\AAA^\T=1_\AAA$. A $*$-algebra (resp., $C^*$-algebra)
equipped with a transposition $\T$ will be called
a \textit{$(*,\T)$-algebra} (resp.,
\textit{$C^{*,\T}$-algebra}).\vadjust{\goodbreak}
\end{Definition}

\begin{Remark}
Of course $\Mat_n(\CC)$ is a $C^{*,\T}$-algebra.
More generally, for any Hilbert space $H$ equipped with an orthonormal
basis $\{h_i\}$,
there exists a unique structure of $C^{*,\T}$-algebra for $B(H)$
such that $(h_i,A h_j)=(h_j,A^\T h_i)$ for all
operators $A\in B(H)$ and indices $i$ and $j$. The concept of $C^{*,\T
}$-algebra is
essentially equivalent to that of a real $C^*$-algebra.
\end{Remark}

\subsubsection{Transpositions, tensor products and matrices}\label
{subsubsectionSesquilinear}
Given $(*,\T)$-alge\-bras $\AAA$ and $\BBB$,
we invariably equip $\AAA\otimes\BBB$
with a transposition by the rule $(a\otimes b)^\T=a^\T\otimes b^\T$,
thus equipping $\AAA\otimes\BBB$ with the structure of $(*,\T
)$-algebra. Note that if $\AAA$ and $\BBB$ are $C^{*,\T}$-algebras
at least one of which is finite-dimensional, then $\AAA\otimes\BBB$
is again a $C^{*,\T}$-algebra.
For any $(*,\T)$-algebra $\AAA$ and matrix $A\in\Mat_{k\times\ell
}(\AAA)$,
we define $A^\T\in\Mat_{\ell\times k}(\AAA)$ by $A^\T
(i,j)=A(j,i)^\T$.
Thus, in particular, $\Mat_n(\AAA)$ is automatically a $(*,\T
)$-algebra (resp., $C^{*,\T}$-algebra)
whenever $\AAA$ is.

\subsubsection{\texorpdfstring{Transpositions on $\CC\langle\Xbold\rangle$ and $B(\HHH)$}
{Transpositions on C<X> and B(H)}}
\label{subsubsectionNaturalTranspositions}
We equip the noncommutative polynomial algebra ${\CC\langle\Xbold
\rangle}$
with a transposition by the rule $\Xbold_\ell^\T=(-1)^\ell\Xbold
_\ell$ for every $\ell$.
Note that the $C^*$-algebra $B(\HHH)$ is
canonically equipped with a transposition because
Boltzmann--Fock space $\HHH$ is
canonically equipped with an orthonormal basis.

\begin{Remark}\label{RemarkCStarTeeMain}
We claim that the evaluation maps
\begin{eqnarray*}
\biggl(f\mapsto f\biggl(\frac{\Xi^N}{\sqrt{N}}\biggr)\biggr)
\dvtx
\Mat_n({\CC\langle\Xbold\rangle})&\rightarrow&\Mat_{nN}(\CC),\\
\bigl(f\mapsto f(\Xi)\bigr) \dvtx  \Mat_n({\CC\langle\Xbold\rangle}
)&\rightarrow&\Mat_n(B(\HHH))
\end{eqnarray*}
figuring in Theorems~\ref{TheoremPreMainResult},~\ref{TheoremMainResult}
and~\ref{TheoremMainResultBis} are $(*,\T)$-algebra homomorphisms.
In any case, it is clear that each homomorphism is a $*$-algebra homomorphism,
so we have only to verify that each map commutes with the transposition.
The former commutes with $\T$
by assumption \eqref{equationPepper5} which (recall) says that
$(\Xi^N_\ell)^\T=(-1)^\ell\Xi^N_\ell$ for all $\ell$.
The latter commutes with $\T$ because
(recall from Section~\ref{subsectionBoltzmannFock}) by definition
$\Xi_\ell=\ii^{\ell}\Sigma_\ell+\ii^{-\ell} \Sigma_\ell^*$,
clearly $\Sigma_\ell^\T=\Sigma_\ell^*$,
and hence $\Xi_\ell^\T=(-1)^\ell\Xi_\ell$
for all $\ell$. The claim is proved.
\end{Remark}

\begin{Lemma}\label{LemmaCStarTee}
If $x$ is an element of a $C^{*,\T}$-algebra $\AAA$,
then $(x^{-1})^\T=(x^\T)^{-1}$, $x\in\AAA_\sa\Rightarrow x^\T\in
\AAA_{\sa}$,
$\Spec( x)=\Spec( x^\T)$, $x\geq0\Rightarrow x^\T\geq0$
and\break \mbox{$[\![x^\T]\!]=[\![x]\!]$}.
\end{Lemma}
\begin{pf} The first two claims are obvious. The third claim follows
from the first.
The second and third claims imply the fourth.
The fifth holds for self-adjoint~$x$ by Proposition \ref
{PropositionCstarNormUnique}
along with the second and third claims.
The fifth claim holds in general because $[\![x^\T]\!]^2=[\![(x^\T
)^*x^\T]\!]
=[\![(xx^*)^\T]\!]=\break[\![xx^*]\!]=[\![x]\!]^2$.
\end{pf}

\begin{Definition}
Given a $C^{*,\T}$-algebra $\AAA$ and a state $\phi\in\AAA^\star$,
we say that $\phi$ is \textit{$\T$-stable}
if $\phi(A^\T)=\phi(A)$ for all $A\in\AAA$.
A pair $(\AAA,\phi)$ consisting of a $C^{*,\T}$-algebra and a $\T
$-stable state
$\phi$ will be called a \textit{$C^{*,\T}$-probability space}.
\end{Definition}

\begin{Remark}
It is easy to see that both $(\Mat_N(\CC),\frac{1}{N}\trace)$ and
$(B(\HHH),\phi^\BF)$ are in fact
$C^{*,\T}$-probability spaces.
\end{Remark}

\subsection{Block algebras (enhanced version)}\label
{subsectionBlockAlgebraUpgrade}
We re-introduce the notion of block algebra, this time with structure
enriched by a transposition.
We also introduce the notion of $\SSS$-linear form in terms of which
we will handle much bookkeeping below.
\begin{Definition}\label{DefinitionBlockAlgebra}
A \textit{block algebra} is a $C^{*,\T}$-algebra
isomorphic to the matrix algebra $\Mat_s(\CC)$ for some integer $s>0$.
A basis $\{e_{ij}\}_{i,j=1}^{s}$ for $\SSS$ such that
$e_{ij}e_{i'j'}=\delta_{ji'}e_{ij'}$ and
$e_{ij}^*=e_{ji}=e_{ij}^\T$ will be called \textit{standard}.
\end{Definition}

\begin{Remark}
A choice of standard basis of a block algebra is the same thing as a
choice of a $C^{*,\T}$-algebra isomorphism with $\Mat_s(\CC)$.
\end{Remark}

\begin{Remark}
The tensor product of block algebras is again a block algebra.
Furthermore, for every block algebra $\SSS$,
the tensor product algebra ${\CC\langle\Xbold\rangle}\otimes\SSS$
[resp.,
$B(\HHH)\otimes\SSS$] is a $(*,\T)$-algebra (resp., $C^{*,\T}$-algebra).
\end{Remark}

\begin{Definition}
Let $\SSS$ be any block algebra.
An \textit{$\SSS$-linear form} $L$ is an element of
the tensor product algebra
$\CC\langle\Xbold\rangle\otimes\SSS$
of the form $L=\sum_{\ell=1}^\infty\Xbold_\ell\otimes a_\ell$
for some elements $a_\ell\in\SSS$ vanishing for $\ell\gg0$.
We refer to the sum $\sum_\ell\Xbold_\ell\otimes a_\ell$ as the
\textit{Hamel expansion}
of $L$ and to the elements $a_\ell\in\SSS$ as the \textit{Hamel
coefficients} of $L$.
Given a sequence $\xi=\{\xi_\ell\}_{\ell=1}^\infty\in\AAA^\infty$
in an algebra $\AAA$,
we define
$L(\xi)=\sum_\ell\xi_\ell\otimes a_\ell\in\AAA\otimes\SSS$,
calling this the \textit{evaluation of $L$ at $\xi$}.
It is especially important to notice that if $\AAA=\Mat_N(\CC)$,
then $L(\xi)\in\Mat_N(\CC)\otimes\SSS=\Mat_N(\SSS)$. This is
the reason
for putting the tensor factors in $\CC\langle\Xbold\rangle\otimes
\SSS$ in the ``wrong'' order.
\end{Definition}

\begin{Definition}\label{DefinitionPhiPsi}
Let $\SSS$ be a block algebra and let $L$ be an $\SSS$-linear form with
Hamel expansion $L=\sum\Xbold_\ell\otimes a_\ell$. We define $\Phi
_L\in B(\SSS)$
by the formula $\Phi_L(\zeta)=\sum a_\ell\zeta a_\ell$ for $\zeta
\in\SSS$
and we define $\Psi_L=\sum(-1)^\ell a_\ell^{\otimes2}\in\SSS
^{\otimes2}$.
We call $\Phi_L$ the \textit{covariance map} attached to $L$.
We call $\Psi_L$ the \textit{covariance tensor} attached
to~$L$.
\end{Definition}

\begin{Definition}\label{DefinitionSSStracialState}
Each block algebra $\SSS$ is equipped with a unique state $\tau_\SSS$
satisfying $\tau_\SSS(e_{ij})=(\dim\SSS)^{-1/2}\delta_{ij}$
for any standard basis $\{e_{ij}\}$. Necessarily $\tau_\SSS$ is
$\T$-stable.
More generally, for each projection $e\in\SSS$, there exists a unique
state $\tau_{\SSS,e}\in\SSS^\star$
such that $\tau_{\SSS,e}\vert_{e\SSS e}=\tau_{e\SSS e}$.
\end{Definition}

\begin{Remark}
If $\SSS\,{=}\,\Mat_s(\CC)$ and $e\,{=}\,\bigl[
{
\Ibold_n\atop 0}\enskip
{0\atop 0}
\bigr]$, then $\tau_{S,e}(A)\,{=}\,\frac{1}{n}\sum_{i=1}^n A(i,i)$.
\end{Remark}

\subsubsection{The bullet map}
Given a block algebra $\SSS$, we define a linear isomorphism
$(A\mapsto A^\bullet)\dvtx \SSS^{\otimes2}\rightarrow B(\SSS)$ by the formula
$(x\otimes y)^\bullet=(z\mapsto xzy)$.
That the bullet map is indeed a
linear isomorphism one can check by calculating with a standard basis.
The bullet map
in general neither preserves norms nor algebra structure.

\subsubsection{The half-transpose map}
Given a block algebra $\SSS$, we define a linear isomorphism
$(A\mapsto A^{1\otimes\T})\in B(\SSS^{\otimes2})$
by the formula $(x\otimes y)^{1\otimes\T}=x\otimes y^\T$. The
half-transpose map
in general neither preserves norms nor algebra structure.
\begin{Remark}\label{RemarkPlotThickens}
Strangely enough, the composite map
\[
\bigl((x\otimes y)\mapsto\bigl((x\otimes y)^{1\otimes\T}\bigr)^\bullet\bigr)\dvtx \SSS
^{\otimes2}\rightarrow B(\SSS)
\]
is an isomorphism of algebras, as one verifies by calculating with a
standard basis.
(But this map still does not in general preserve norms.)
\end{Remark}

\subsection{$\SSS$-(bi)linear constructions}

\subsubsection{$\SSS$-linear extension of states}\label
{subsubsectionSlinearExtStates}
Given any $C^*$-probability space $(\AAA,\phi)$ and block algebra
$\SSS$,
we define the \textit{$\SSS$-linear extension}
$\phi_\SSS\dvtx \AAA\otimes\SSS\rightarrow\SSS$
of $\phi$ by the formula
\[
\phi_\SSS(x\otimes y)=\phi(x)y.
\]
Note that since $\phi$ commutes with the involution,
the same is true for $\phi_\SSS$, that is,
%
%
\begin{equation}\label{equationExtensionStarSymmetry}
\phi_\SSS(A^*)=\phi_\SSS(A)^*
\end{equation}
for $A\in\AAA\otimes\SSS$. Suppose now that
$(\AAA,\phi)$ is a $C^{*,\T}$-probability space.
Note that since $\phi$ is $\T$-stable,
$\phi_\SSS$ commutes with $\T$, that is,
%
%
\begin{equation}\label{equationExtensionTransSymmetry}
\phi_\SSS(A^\T)=\phi_\SSS(A)^\T
\end{equation}
for $A\in\AAA\otimes\SSS$.

\begin{Remark}
Objects like $\phi_\SSS$ are the stock-in-trade of operator-valued
free probability theory.
See~\cite{HRS07} for a useful introduction to this point of view
and see~\cite{NicaSpeicher} for in-depth treatment.
\end{Remark}

\begin{Remark}\label{RemarkUtterly}
Consider the case $(\AAA,\phi)=(\Mat_N(\CC),\frac{1}{N}\trace)$.
We have
\[
\phi_\SSS=\Biggl(A\mapsto\frac{1}{N}\trace_\SSS A=\sum_{i=1}^N
A(i,i)\Biggr)\dvtx \Mat_N(\SSS)\rightarrow\SSS.
\]
Thus the ad hoc construction $\trace_\SSS$ fits into a more
general conceptual framework.
\end{Remark}

\begin{Remark}\label{RemarkCoffeePress}
The $\SSS$-linear extension $\phi^{\BF}_\SSS$
of the state $\phi^{\BF}$ with which $B(\HHH)$
is canonically equipped satisfies
%
%
\begin{equation}
\label{equationStateBound}
(p_\HHH\otimes1_\SSS)A (p_\HHH\otimes1_\SSS)=p_\HHH\otimes\phi
^{\BF}_\SSS(A) \qquad\mbox{hence } [\![\phi^{\BF}_\SSS(A)]\!]\leq
[\![A]\!]
\end{equation}
for all $A\in B(\HHH)\otimes\SSS$ and hence $[\![\phi^{\BF }_\SSS]\!]=1$.
In fact, in full generality, we have $[\![\phi_\SSS]\!]=1$
by a similar argument using the GNS construction, which we omit.
\end{Remark}

\begin{Definition}
For any $(*,\T)$-algebra $\AAA$,
let $\AAA_{\alt}^\infty$ denote the space of
sequences $\xi=\{\xi_\ell\}_{\ell=1}^\infty$ in $\AAA$
such that $\xi_\ell^\T=(-1)^\ell\xi_\ell$ for all $\ell$.
Also put $\AAA^\infty_{\salt}=\AAA_\sa^\infty\cap\AAA_\alt
^\infty$.
\end{Definition}
\begin{Remark}
Let $\AAA$ be a $(*,\T)$-algebra, $\xi=\{\xi_\ell\}_{\ell
=1}^\infty\in\AAA^\infty_\salt$ a sequence and $L$ an $\SSS
$-linear form.
Then we have $L^\T(\xi)=L(\xi)^\T$ and $L^*(\xi)=L(\xi)^*$.
In particular, this observation
applies to the sequences
$\Xi^N\in\Mat_N(\CC)_\salt^\infty$
and $\Xi\in B(\HHH)_\salt^\infty$
figuring in Theorem~\ref{TheoremMainResultBis}.
\end{Remark}

\subsubsection{$\SSS$-bilinear extension of states}
Let $\SSS$ be a block algebra and let $(\AAA,\phi)$ be a
$C^*$-probability space.
We define the \textit{$\SSS$-bilinear extension}
\[
\phi_{\SSS,\SSS}\dvtx \AAA\otimes\SSS\times\AAA\otimes\SSS
\rightarrow
\SSS^{\otimes2}
\]
of $\phi$ by the formula
\[
\phi_{\SSS,\SSS}(x_1\otimes y_1,x_2\otimes y_2)=\phi
(x_1x_2)y_1\otimes y_2.
\]

\begin{Remark}
Let $\SSS$ be a block algebra and consider the $C^*$-probability space
$(\Mat_N(\CC),\frac{1}{N}\trace)$.
For $R_1,R_2\in\Mat_N(\SSS)$, we have
\[
\phi_{\SSS,\SSS}(R_1,R_2)=\frac{1}{N}\sum_{i,j=1}^N
R_1(i,j)\otimes R_2(j,i)\in\SSS^{\otimes2}.
\]
\end{Remark}
\begin{Remark}\label{RemarkCoffeePressBis}
Consider the $C^*$-algebra embeddings
\[
\left.
\begin{array}{l}
\iota^{(1)}=(x\otimes y\mapsto x\otimes y\otimes1_\SSS)\\
\iota^{(2)}=(x\otimes y\mapsto x\otimes1_\SSS\otimes y)
\end{array}
\right\}\dvtx \AAA\otimes\SSS\rightarrow\AAA\otimes\SSS^{\otimes2}.
\]
One has
%
%
\begin{equation}\label{equationPreBulletCuriosity}
\phi_{\SSS,\SSS}(A,B)=\phi_{\SSS^{\otimes2}}\bigl(\iota^{(1)}(A)\iota
^{(2)}(B)\bigr)
\end{equation}
and thus $[\![\phi_{\SSS,\SSS}]\!]=1$ since $[\![\phi_{\SSS ^{\otimes
2}}]\!]=1$.
In a similar vein, we have the formula
%
%
\begin{equation}\label{equationBulletCuriosity}
\phi_{\SSS,\SSS}(A,B)^\bullet(\zeta)=\phi_\SSS\bigl(A(1_\AAA\otimes
\zeta)B\bigr)
\end{equation}
which we will use below in Section~\ref{sectionOpTheoSchwingerDyson}
to study the secondary Schwinger--Dyson equation.
\end{Remark}

\begin{Remark}\label{RemarkUtterlyTer}
Consider the $C^{*,\T}$-probability space $(\AAA,\phi)=(\Mat_N(\CC
),\break \frac{1}{N}\trace)$.
Let $\SSS$ be any block algebra. Let $R\in\Mat_N(\SSS)$ be any matrix.
We have\looseness=-1
%
%
\begin{eqnarray}
\label{equationBulletTransposeRationale}
\phi_{\SSS,\SSS}(R,R)^\bullet&=&\Biggl(\zeta\mapsto\frac
{1}{N}\sum_{i,j=1}^N
R(i,j)\zeta R(j,i)\Biggr)\in B(\SSS),\\[-3pt]
\label{equationHalfTransposeRationale}
\phi_{\SSS,\SSS}(R,R^\T)^{1\otimes\T}&=&\frac{1}{N}\sum_{i,j=1}^N
R(i,j)^{\otimes2}\in\SSS^{\otimes2}.
\end{eqnarray}\looseness=0
We will study the expressions of the form on the left in Section \ref
{sectionOpTheoSchwingerDyson} below.
We will study expressions of the form on the right in
Section~\ref{sectionUrMatrixIdentities} and Section \ref
{sectionMatrixIdentities} below.\vspace*{-2pt}
\end{Remark}

\subsection{An upgrade of the self-adjoint linearization trick}\vspace*{-2pt}

\begin{Definition}\label{DefinitionSALT}
A \textit{SALT block design} is a quadruple $(\SSS,L,\Theta,e)$
consisting of
\begin{itemize}
\item a block algebra $\SSS$,
\item a self-adjoint $\SSS$-linear form $L$,
\item an element $\Theta\in\SSS$ (perhaps not self-adjoint) and
\item a projection $e\in\SSS$
\end{itemize}
such that for every $C^{*,\T}$-algebra $\AAA$, sequence $\xi\in\AAA
_\salt^\infty$,
point $z\in{\mathfrak{h}}$ and parameter $t\geq0$ we have
%
%
\begin{equation}\label{equationSALTDef2}
L(\xi)-1_\AAA\otimes(\Theta+ze+\ii t1_\SSS)\in(\AAA\otimes
\SSS)^\times
\end{equation}
and
\begin{equation}
\label{equationSALTDef3}
\qquad\bigl[\!\bigl[\bigl(L(\xi)-1_\AAA\otimes(\Theta+ze+\ii t1_\SSS )\bigr)^{-1}\bigr]\!\bigr]
\leq
\frac{c_0}{2}\bigl(1+[\![L(\xi)]\!]\bigr)^{c_1}(1+1/\Im z)^{c_2}
\end{equation}
for some constants $c_0,c_1,c_2\geq1$ depending only on $(\SSS
,L,\Theta,e)$
and thus independent of $\AAA$, $\xi$, $z$ and $t$.
We declare any finite constant
$\Tfrak\geq[\![\Im\Theta]\!]+4(1+[\![\Phi_L]\!])$
to be a \textit{cutoff} for the design, where $\Phi_L\in B(\SSS)$ is as
in Definition~\ref{DefinitionPhiPsi}.\vspace*{-2pt}
\end{Definition}

\begin{Proposition}\label{PropositionSALTnouveau}
Fix $f\in\Mat_n({\CC\langle\Xbold\rangle})_\sa$ arbitrarily.
Either let $\tilde{f}\in\break
\Mat_s({\CC\langle\Xbold\rangle})_\sa$
be a self-adjoint linearization of $f$ as provided by Proposition \ref
{PropositionNaiveSALT}
or else let $\tilde{f}=f$ if all entries of $f$
are already of degree $\leq1$ in the variables $\Xbold_\ell$. Write
\[
\tilde{f}=a_0\otimes1_{{\CC\langle\Xbold\rangle}}+\sum_{\ell
=1}^m a_\ell\otimes
\Xbold_\ell
\]
for $m\gg0$ and $a_0,\ldots,a_m\in\Mat_s(\CC)_\sa$.
Let
\[
\SSS=\Mat_s(\CC),\qquad
L=\sum_{\ell=1}^m \Xbold_\ell\otimes a_\ell,\qquad \Theta=-a_0,\qquad
e=\left[
\matrix{
\Ibold_n&0\vspace*{2pt}\cr
0&0}
\right].
\]
Then $(\SSS,L,\Theta,e)$ is a SALT block design for which $\Theta
=\Theta^*$ and $c_2=1$.\vadjust{\goodbreak}
\end{Proposition}
\begin{pf}
The case $s=n$ is easy to check using Lemma~\ref{Lemmalambda}. We
leave the details to the reader.
We assume $s>n$ for the rest of the proof.
Then property~\eqref{equationSALTDef2} holds by Proposition \ref
{PropositionNaiveSALTbis} and Lemma~\ref{Lemmalambda}.
The latter lemma, the inversion formula \eqref{equationFamiliarInversion},
Lemma~\ref{LemmaBlockization} and the definitions yield an estimate
\[\bigl[\!\bigl[\bigl(L(\xi)-1_\AAA\otimes(\Theta+ze)\bigr)^{-1}\bigr]\!\bigr]
\leq
\frac{c_0}{4}\bigl(1+[\![L(\xi)]\!]\bigr)^{c_1}(1+1/\Im z)
\]
for $\AAA$, $\xi$ and $z$ as in \eqref{equationSALTDef3}. Estimate
\eqref{equationSALTDef3} for general $t\geq0$ (but still with
$c_2=1$) then follows via Lemma~\ref{LemmaNeumannExpansion} if $\frac
{c_0}{2}(1+[\![L(\xi)]\!])^{c_1}(1+1/\Im z)t\leq1$
and otherwise follows via Lemma~\ref{Lemmalambda}.
\end{pf}

\begin{Definition}\label{DefinitionSALTupgrade}
For each $f\in\Mat_n({\CC\langle\Xbold\rangle})_\sa$,
any SALT block design $(\SSS,L,\Theta,e)$ arising via Proposition
\ref{PropositionSALTnouveau}
will be called a \textit{self-adjoint linearization} of $f$.
(Self-adjoint linearizations in the relatively naive sense defined
immediately after Proposition~\ref{PropositionNaiveSALT}
will hereafter no longer be used.)
\end{Definition}

\begin{Remark}\label{RemarkPreCobbling}
(This is a continuation of Remark~\ref{RemarkPrePreCobbling}.)
In the sequel, we will prefer to write formulas \eqref
{equationEmpiricalAccess} and \eqref{equationPreStieltjesAccess}
above in the following form.
Fix $f\in\Mat_n({\CC\langle\Xbold\rangle})_\sa$
and let $(\SSS,L,\Theta,e)$ be any self-adjoint linearization of $f$.
Then, with $\tau_{\SSS,e}$ as in
Definition~\ref{DefinitionSSStracialState}, we have
%
%
\begin{eqnarray}
\label{equationPreCobblingRandom}
S_{\mu_f^N}(z)&=&\tau_{\SSS,e}
\biggl(\frac{1}{N} \trace_\SSS\biggl(\biggl(L\biggl(\frac{\Xi
^N}{\sqrt{N}}\biggr)-\Ibold_N\otimes(\Theta+ze)\biggr)^{-1}
\biggr)\biggr),\\
\label{equationPreCobblingOperator}
S_{\mu_f}(z)&=&\tau_{\SSS,e}
\circ\phi^{\BF}_\SSS\bigl(\bigl(L(\Xi)-1_{B(\HHH)}\otimes
(\Theta+ze)\bigr)^{-1}\bigr)
\end{eqnarray}
for $z\in{\mathfrak{h}}$ (resp., for $z\in\CC\setminus\supp\mu_f$).
Furthermore, for distinct indices $i,j=1,\ldots,N$
and arbitrary $A\in\SSS$, with the covariance map $\Phi_L\in B(\SSS
)$ as in Definition~\ref{DefinitionPhiPsi}, we have
%
%
\begin{eqnarray}\label{equationPhiEllRationale}
\Phi_L(A)&=&\Ebold((L(\Xi^N)(i,j)
)A(L(\Xi^N)(j,i)))
\nonumber
\\[-8pt]
\\[-8pt]
\nonumber
\nonumber&=&\phi_\SSS
\bigl(L(\Xi)\bigl(1_{B(\HHH)}\otimes A\bigr)L(\Xi)\bigr).
\end{eqnarray}
[Formula \eqref{equationPhiEllRationale} is actually valid for any
self-adjoint $\SSS$-linear form $L$ whether or not it is part of a
self-adjoint linearization.] One verifies the first equality of \eqref
{equationPhiEllRationale} by applying assumptions \eqref{equationPepper5},
\eqref{equationPepper6}, \eqref{equationPepper7} and \eqref
{equationPepper8}
which fix the covariance structure of the sequence $\Xi^N$.
One verifies the second equality of \eqref{equationPhiEllRationale}
by applying the relations
\eqref{equationSigmaMoments} which analogously fix
the (noncommutative) covariance structure
of the sequence $\Xi$.
We also have for distinct indices $i,j=1,\ldots,N$ the formula
%
%
\begin{equation}\label{equationPsiEllRationale}
\Psi_L=\Ebold((L(\Xi^N)(i,j))^{\otimes2})=\phi_{\SSS
,\SSS}(L(\Xi),L(\Xi))
\end{equation}
which is proved more or less the same way as formula \eqref
{equationPhiEllRationale}.
[Formula \eqref{equationPsiEllRationale} is actually valid for any
$\SSS$-linear form $L$.]
\end{Remark}

\begin{Remark}\label{RemarkCutoff}
We return to the setting of Definition~\ref{DefinitionSALT}.
We provide some amplification of relations \eqref{equationSALTDef2}
and \eqref{equationSALTDef3}.
Firstly, we have for $\zeta\in\SSS$ that
%
%
\begin{eqnarray}\label{equationSALTDef5}
&&2\bigl[\!\bigl[\bigl(L(\xi)-1_\AAA\otimes(\Theta+ze+\ii t1_\SSS )\bigr)^{-1}\bigr]\!\bigr][\![\zeta]\!]\leq1\nonumber\\
&&\quad\Rightarrow\quad
L(\xi)-1_\AAA\otimes(\Theta+ze+\ii t1_\SSS+\zeta)\in(\AAA
\otimes\SSS)^\times \quad\mbox{and}
\nonumber
\\[-8pt]
\\[-8pt]
\nonumber
&&\hspace*{41pt}\bigl[\!\bigl[\bigl(L(\xi)-1_\AAA\otimes(\Theta+ze+\ii t1_\SSS+\zeta
)\bigr)^{-1}\bigr]\!\bigr]\\
&&\hspace*{41pt}\qquad
\leq 2\bigl[\!\bigl[\bigl(L(\xi)-1_\AAA\otimes(\Theta+ze+\ii t1_\SSS
)\bigr)^{-1}\bigr]\!\bigr]\nonumber
\end{eqnarray}
via Lemma~\ref{LemmaNeumannExpansion}.
Secondly, we have
%
%
\begin{equation}\label{equationSALTDef4}
\bigl[\!\bigl[\bigl(L(\xi)-1_\AAA\otimes(\Theta+ze+\ii t1_\SSS)\bigr)^{-1}\bigr]\!\bigr]
\leq
\frac{1}{t-[\![\Im\Theta]\!]}\qquad \mbox{for $t> [\![\Im \Theta]\!]$}
\end{equation}
via Lemma~\ref{Lemmalambda}. The bound latter helps to explain the
meaning of the cutoff~$\Tfrak$.
\end{Remark}

\begin{Remark}\label{RemarkYCapsuleComplement}
Let $(\SSS,L,\Theta,e)$ be a SALT block design and let
$c_0$, $c_1$, $c_2$ and $\Tfrak$ be the constants
from Definition~\ref{DefinitionSALT}. Let $Y\in\Mat_N(\SSS)_\sa$
be of the form $Y=L(\eta)$ for some $\eta\in\Mat_N(\CC)_\salt
^\infty$.
Then for every $z\in{\mathfrak{h}}$ and $t\in[0,\infty)$ we have
\begin{eqnarray}
&&Y-\Ibold_N\otimes(\Theta+ze+\ii t1_\SSS)\in\GL_N(\SSS)
\quad\mbox{and}\nonumber\\
&&\bigl[\!\bigl[\bigl(Y-\Ibold_N\otimes(\Theta+ze+\ii t1_\SSS )\bigr)^{-1}\bigr]\!\bigr]\\
&&\qquad\leq \cases{
c_0(1+[\![Y]\!])^{c_1}(1+1/\Im z)^{c_2},&\quad$\mbox{in general
and}$\vspace*{2pt}\cr
\displaystyle\frac{1}{2},
&\quad $\mbox{for $t\geq\Tfrak$}$}\nonumber
\end{eqnarray}
by definition of a SALT block design
along with Remark~\ref{RemarkCutoff}.
By the resolvent identity \eqref{equationResolventIdentity}, the
following important (if trivial) observation also holds:
\begin{eqnarray}\label{equationCrucialLipschitz}
\begin{tabular}{p{260pt}@{}}
\mbox{For each fixed\vspace*{2pt} $z\in{\mathfrak{h}}$, the Lipschitz constant
of the map}
$(t\mapsto(Y-\Ibold_N\otimes
(\Theta+ze+\ii t1_\SSS))^{-1}
)\dvtx [0,\infty)\rightarrow\Mat_N(\SSS)$\vspace*{2pt}
$\mbox{does not exceed } c_0^2(1+[\![Y]\!])^{2c_1}(1+1/\Im z)^{2c_2}.$
\end{tabular}
\end{eqnarray}
\end{Remark}
%

\subsection{The secondary trick}\label{subsectionCoffeePress}
We first present the underline construction
and then the secondary trick itself. As a byproduct, we will find
natural examples of SALT block designs beyond those produced by
Proposition~\ref{PropositionSALTnouveau}.
See Remark~\ref{RemarkForceOnUs} below.

\subsubsection{The underline construction for a block algebra}\label
{subsubsectionUnderlineBlock}
Let $\MMM_3$ denote a block algebra equipped with a standard basis $\{
e_{ij}\}_{i,j=1}^3$.
Let $\SSS$ be any block algebra.
We define
\[
\underline{\SSS}=\SSS^{\otimes2}\otimes\MMM_3,\qquad
\Diamond_\SSS=1_\SSS^{\otimes2}\otimes(e_{12}+e_{13})\in
\underline{\SSS}.
\]
Furthermore, given $\Lambda\in\SSS$ we define
\begin{eqnarray*}
\underline{\Lambda}_1&=&\Lambda\otimes1_\SSS\in\SSS^{\otimes2},
\qquad
\underline{\Lambda}_2=1_\SSS\otimes\Lambda\in\SSS^{\otimes2},
\\
\underline{\Lambda}&=&\underline{\Lambda}_1\otimes e_{11}+\underline
{\Lambda}_2\otimes e_{22}+
\underline{\Lambda}_2^\T\otimes e_{33}\in\underline{\SSS}.
\end{eqnarray*}
We also define linear maps
\[
\partial_1\in B(\underline{\SSS},B(\SSS))\quad \mbox{and}\quad
\partial_2\in B(\underline{\SSS},\SSS^{\otimes2})
\]
by the formulas
\[
\partial_1(A\otimes e_{ij})=A^\bullet\delta_{1i}\delta_{2j},\qquad
\partial_2(A\otimes e_{ij})=A^{1\otimes\T}\delta_{1i}\delta_{3j}
\]
for $A\in\SSS^{\otimes2}$ and $i,j=1,\ldots,3$.

\subsubsection{The underline construction for $\SSS$-linear forms}
\label{subsubsectionUnderlineLinearForm}
Let $\SSS$ be a block algebra and let $\underline{\SSS}$ be the corresponding
``underlined'' block algebra as defined in the preceding paragraph.
Given an $\SSS$-linear form $L$ with Hamel expansion
$L=\sum\Xbold_\ell\otimes a_\ell$, we define $\SSS^{\otimes
2}$-linear forms
$\underline{L}_1$ and $\underline{L}_2$ by the formulas
\[
\underline{L}_1=\sum_\ell\Xbold_\ell\otimes a_\ell\otimes1_\SSS
,\qquad
\underline{L}_2=\sum_\ell\Xbold_\ell\otimes1_\SSS\otimes a_\ell
\]
and in turn we define an $\underline{\SSS}$-linear form $\underline
{L}$ by the formula
\[
\underline{L}=\underline{L}_1\otimes e_{11}+\underline{L}_2\otimes e_{22}+
\underline{L}_2^\T\otimes e_{33}.
\]

\subsubsection{The trick itself}\label{subsubsectionSecondaryTrick}
Let $\SSS$ be a block algebra,
let $\Lambda\in\SSS$ be an element and let $L$ be an $\SSS$-linear form.
Let $\underline{\SSS}$, $\underline{\Lambda}_1$, $\underline
{\Lambda}_2$,
$\underline{\Lambda}$, $\underline{L}_1$, $\underline{L}_2$
and $\underline{L}$ be as defined in the preceding two paragraphs.
Let $(\AAA,\phi)$ be a $C^{*,\T}$-probability space
and fix a sequence $\xi\in\AAA^\infty_\salt$.
For the rest of this paragraph, we abuse notation by writing $x=1_\AAA
\otimes x$
for $x\in\SSS$, $x\in\SSS^{\otimes2}$ or $x\in\underline{\SSS}$.
We assume that $L(\xi)-\Lambda\in(\AAA\otimes\SSS)^\times$,
in which case $\underline{L}(\xi)-\underline{\Lambda}\in(\AAA
\otimes\underline{\SSS})^\times$ via Lemma~\ref{LemmaCStarTee}.
In turn,
one can verify that
\[
\underline{L}(\xi)-\underline{\Lambda}-\Diamond_\SSS\in(\AAA
\otimes\underline{\SSS})^\times
\]
and more precisely
\begin{eqnarray*}
&&\bigl(\underline{L}(\xi)-\underline{\Lambda}-\Diamond_\SSS\bigr)^{-1}\\
&&\qquad=\bigl(\underline{L}_1(\xi)-\underline{\Lambda}_1\bigr)^{-1}\otimes e_{11}+
\bigl(\underline{L}_2(\xi)-\underline{\Lambda}_2\bigr)^{-1}\otimes e_{22}\\
&&\qquad\quad{}+
\bigl(\bigl(\underline{L}_2(\xi)-\underline{\Lambda}_2\bigr)^{-1}\bigr)^\T\otimes
e_{33}\\
&&\qquad\quad{}+\bigl(\bigl(\underline{L}_1(\xi)-\underline{\Lambda}_1\bigr)^{-1}\bigl(\underline
{L}_2(\xi)-\underline{\Lambda}_2\bigr)^{-1}\bigr)
\otimes e_{12}\\
&&\qquad\quad{}+\bigl(\bigl(\underline{L}_1(\xi)-\underline{\Lambda}_1\bigr)^{-1}\bigl(\bigl(\underline
{L}_2(\xi)-\underline{\Lambda}_2\bigr)^{-1}\bigr)^\T\bigr)\otimes e_{13},
\end{eqnarray*}
by direct calculation.
It follows that
%
%
\begin{eqnarray}
\label{equationUnderlineConstruction1}
&&[\![L(\xi)]\!]=[\![\underline{L}(\xi)]\!],
\nonumber
\\[-8pt]
\\[-8pt]
\nonumber
&&\bigl[\!\bigl[\bigl(\underline{L}(\xi)-\underline{\Lambda}-\Diamond_\SSS \bigr)^{-1}\bigr]\!\bigr] \leq3\bigl(1\vee\bigl[\!\bigl[\bigl(L(\xi)-\Lambda\bigr)^{-1}\bigr]\!\bigr]\bigr)^2,\\
\label{equationUnderlineConstruction3}
&&\partial_1\circ\phi_{\underline{\SSS}}\bigl(\bigl(\underline{L}(\xi
)-\underline{\Lambda}-\Diamond_\SSS\bigr)^{-1}\bigr)
\nonumber
\\[-8pt]
\\[-8pt]
\nonumber
&&\qquad = \phi_{\SSS,\SSS
}\bigl(\bigl(L(\xi)-\Lambda\bigr)^{-1},\bigl(L(\xi)-\Lambda\bigr)^{-1}\bigr)^\bullet,\\
\label{equationUnderlineConstruction4}
&&\partial_2\circ\phi_{\underline{\SSS}}\bigl(\bigl(\underline{L}(\xi
)-\underline{\Lambda}-\Diamond_\SSS\bigr)^{-1}\bigr)
\nonumber
\\[-8pt]
\\[-8pt]
\nonumber
&&\qquad= \phi_{\SSS,\SSS
}\bigl(\bigl(L(\xi)-\Lambda\bigr)^{-1},\bigl(\bigl(L(\xi)-\Lambda\bigr)^{-1}\bigr)^\T\bigr)^{1\otimes\T},
\end{eqnarray}
where to get the last two identities we use the trivial formula \eqref
{equationPreBulletCuriosity}.

\begin{Lemma}\label{LemmaForcedOnUs}
For any SALT block design $(\SSS,L,\Theta,e)$,
again
$(\underline{\SSS},\underline{L},\underline{\Theta}+\Diamond_\SSS
,\underline{e})$
is a SALT block design.
More precisely, if $c_0$, $c_1$ and $c_2$ are constants
rendering the estimate
\eqref{equationSALTDef3} valid for $(\SSS,L,\Theta,e)$,
then one can take the corresponding constants
$\underline{c}_0$, $\underline{c}_1$ and $\underline{c}_2$
for $(\underline{\SSS},\underline{L},\underline{\Theta}+\Diamond
_\SSS,\underline{e})$
to be $\underline{c}_0=3c_0^2$, $\underline{c}_1=2c_1$
and $\underline{c}_2=2c_2$.
\end{Lemma}
\begin{pf} One can read off
the necessary estimates from \eqref{equationUnderlineConstruction1}.
\end{pf}

\begin{Remark}\label{RemarkForceOnUs}
Note that $\underline{\Theta}+\Diamond_\SSS\in\underline{\SSS}$
is in general not self-adjoint and
(more significantly) $\underline{c}_2=2c_2$.
Thus, the extra generality in
Definition~\ref{DefinitionSALT} not used by Proposition \ref
{PropositionSALTnouveau} is forced on us in order to make Lemma \ref
{LemmaForcedOnUs} hold.
\end{Remark}

\section{Construction of solutions of the Schwinger--Dyson equation}
\label{sectionOpTheoSchwingerDyson}
We construct solutions of the Schwinger--Dyson equation by using the
Boltzmann--Fock apparatus reviewed in Section \ref
{sectionOperatorTheoryTools} above
along with the $\SSS$-linear machinery introduced in Section \ref
{sectionUpgrade} above.
See Proposition~\ref{PropositionConcreteSchwingerDyson} below.
Following~\cite{HST}, we then express the Stieltjes transform $S_{\mu
_f}(z)$ figuring in
Theorem~\ref{TheoremMainResultBis}
in terms of one of the solutions so constructed. See Remark \ref
{RemarkCobblingBis} below.
We also construct solutions of
a secondary version of the Schwinger--Dyson equation by using the $\SSS
$-bilinear machinery of Section~\ref{sectionUpgrade}.
See Proposition~\ref{PropositionSecondarySchwingerDyson} below.
We then define our candidate for the correction ${\mathrm
{bias}}_f^N(z)$ figuring in Theorem~\ref{TheoremMainResultBis}.
See Remark~\ref{RemarkCobblingTer} below.




\subsection{The Schwinger--Dyson equation and its differentiated
form}\label{subsectionSD}

\begin{Definition}\label{DefinitionSD}
Let $\SSS$ be a block algebra.
Let $\DDD\subset\SSS$ be a (nonempty) open subset.
Let $\Phi\in B(\SSS)$ be a linear map.
We say that an analytic function $G\dvtx \DDD\rightarrow\SSS$
satisfies the \textit{Schwinger--Dyson (SD) equation} with \textit{covariance}
\textit{map} $\Phi$
if
\[
1_\SSS+\bigl(\Lambda+\Phi(G(\Lambda))\bigr)G(\Lambda)=0
\]
for all $\Lambda\in\DDD$. Necessarily one has $G(\Lambda)\in\SSS
^\times$ for all $\Lambda\in\DDD$.\vadjust{\goodbreak}
\end{Definition}

See~\cite{VDN,NicaSpeicher} or~\cite{AGZ} for background.

\subsubsection{Notation for derivatives}
Given an analytic function $G\dvtx \DDD\rightarrow\SSS$ defined on an
open subset $\DDD$
of a block algebra $\SSS$ and $\Lambda\in\DDD$, we define
\[
\DD[G](\Lambda)=\biggl(\zeta\mapsto\frac{d}{dt}G(\Lambda+t\zeta
)\bigg\vert_{t=0}\biggr)\in B(\SSS).
\]
For $\zeta\in\SSS$ we write
$\DD[G](\Lambda;\zeta)=\DD[G](\Lambda)(\zeta)$ to compress notation.

\begin{Proposition}\label{PropositionSchwingerDysonCorrection}
Let $\SSS$ be a block algebra and let $\DDD\subset\SSS$ be an open set.
Let $G\dvtx \DDD\rightarrow\SSS$ be a solution of the SD equation
with covariance map $\Phi\in B(\SSS)$. Then for every $\Lambda\in
\DDD$
and $\zeta\in\SSS$ we have
%
%
\begin{eqnarray}\label{equationSchwingerDysonCorrection}
\zeta&=&G(\Lambda)^{-1}\DD[G](\Lambda;\zeta)G(\Lambda)^{-1}-\Phi
(\DD[G](\Lambda;\zeta))
\nonumber
\\[-8pt]
\\[-8pt]
\nonumber
&=&\DD[G]\bigl(\Lambda;G(\Lambda)^{-1}\zeta G(\Lambda)^{-1}-\Phi(\zeta
)\bigr),\nonumber\\
\label{equationSchwingerDysonCorrectionSpecial}
0&=&G(\Lambda)+\DD[G](\Lambda;\Lambda)+2\DD[G](\Lambda;\Phi
(G(\Lambda))).
\end{eqnarray}
\end{Proposition}

Relation \eqref{equationSchwingerDysonCorrectionSpecial}
plays a key role
in proving the bound \eqref{equationMainResultBis2} asserted in
Theorem~\ref{TheoremMainResultBis}.
\begin{pf} To compress notation further, we write $G=G(\Lambda)$ and
$G'=\DD[G](\Lambda)$.
By differentiation of the SD equation we obtain
$(\zeta+\Phi(G'(\zeta)))G+(\Lambda+\Phi(G))G'(\zeta)=0$
and hence $\zeta=G^{-1}G'(\zeta)G^{-1}-\Phi(G'(\zeta))$.
Thus, the first equality in \eqref{equationSchwingerDysonCorrection} holds.
Now for any linear operators $A$ and $B$ on a finite-dimensional vector
space we have
$AB=1\Rightarrow BA=1$. Thus
$\zeta=G'(G^{-1}\zeta G^{-1}-\Phi(\zeta))$,
and hence the second equality in \eqref
{equationSchwingerDysonCorrection} holds.
Finally, \eqref{equationSchwingerDysonCorrectionSpecial}
follows by taking $\zeta=G(\Lambda)$ in \eqref
{equationSchwingerDysonCorrection}.
\end{pf}

\subsection{The solution of the SD equation attached to an $\SSS
$-linear form}

\subsubsection{\texorpdfstring{Definition of $G_L(\Lambda)$}{Definition of G L(Lambda)}}

Fix a block algebra $\SSS$ and an $\SSS$-linear form
$L$ with Hamel expansion $L=\sum\Xbold_\ell\otimes a_\ell$.
We define the set
\[
\DDD_L=\bigl\{\Lambda\in\SSS\vert L(\Xi)-1_{B(\HHH)
}\otimes\Lambda\in\bigl(B(\HHH)\otimes\SSS\bigr)^\times\bigr\}\subset\SSS.
\]
It is clear that $\DDD_L$ is nonempty and Lemma~\ref{LemmaNeumannExpansion}
implies that $\DDD_L$ is open.
For $\Lambda\in\DDD_L$, we put
\[
G_L(\Lambda)=\phi^\BF_\SSS\bigl(\bigl(L(\Xi)- 1_{B(\HHH)}\otimes\Lambda
\bigr)^{-1}\bigr)\in\SSS,
\]
where $\phi^\BF_\SSS$ is the $\SSS$-linear extension of $\phi^\BF
\in B(\HHH)^\star$.
By direct manipulation of series expansions one can verify that
$G_L\dvtx \DDD_L\rightarrow\SSS$ is an analytic function.
Recall that by Definition~\ref{DefinitionPhiPsi},
we have $\Phi_L=(\zeta\mapsto\sum a_\ell\zeta a_\ell)\in B(\SSS)$.

\begin{Proposition}\label{PropositionConcreteSchwingerDyson}
The function
$G_L\dvtx \DDD_L\rightarrow\SSS$
is a solution of the SD equation with covariance map $\Phi_L$.\vadjust{\goodbreak}
\end{Proposition}
\begin{pf}
We specialize Proposition~\ref{PropositionAbstractSchwingerDyson} by taking
\begin{eqnarray*}
\AAA&=&B(\HHH)\otimes\SSS,\qquad \pi=p_\HHH\otimes1_\SSS,\qquad \rho
_\ell=\hat{\Sigma}_\ell
\otimes1_\SSS \quad\mbox{and}
\\
A&=&L(\Xi)-1_{B(\HHH)}\otimes\Lambda=-1_{B(\HHH)}\otimes
\Lambda+\sum_\ell
( \ii^\ell\Sigma_\ell\otimes a_\ell+\ii^{-\ell}\Sigma
_\ell^*\otimes a_\ell).
\end{eqnarray*}
To verify that the family, $\{\pi\}\cup\{\rho_i\}_{\ell=1}^\infty$
is a Cuntz frame in $\AAA$
we use \eqref{equationRhoRelations}. To verify that $A$
is quasi-circular, we use Lemma~\ref{LemmaQuasiCasimir}.
Now in view of \eqref{equationStateBound},
the left side of \eqref{equationSDMachine3} specializes to $p_\HHH
\otimes G_L(\Lambda)$
and moreover necessarily $G_L(\Lambda)\in\SSS^\times$.
But we also have
\[
\pi A\pi=-p_\HHH\otimes\Lambda,\qquad
\pi A\rho_\ell\pi=\ii^\ell p_\HHH\otimes a_\ell,\qquad
\pi\rho_\ell^*A\pi=\ii^{-\ell}p_\HHH\otimes a_\ell
\]
as one verifies by using \eqref{equationSigmaRelations} and
\eqref{equationRaisingLoweringRelations}.
Thus, the inverse in the algebra $\pi\AAA\pi$ of the right side of
\eqref{equationSDMachine3}
specializes to $-p_\HHH\otimes(\Lambda+\Phi_L(G_L(\Lambda)))$.
\end{pf}

\begin{Remark}
Proposition~\ref{PropositionConcreteSchwingerDyson} is essentially well-known
apart from one small detail.
For comparison with a typical proof, see~\cite{AGZ}, Chapter~5, Sections 4 and 5 (main text, not the exercises), and in particular
\cite{AGZ}, Chapter~5, Lemma~5.5.10.
That proof falls a bit short of proving Proposition \ref
{PropositionConcreteSchwingerDyson} as stated
because it relies on an analytic continuation argument to extend a generating
function identity proved by combinatorics throughout a connected open set.
But we do not know a priori that $\DDD_L$ is connected. (It
would be a surprise if it were not but we leave the question aside.)
Thus, we have presented the operator-theoretic proof
of Proposition~\ref{PropositionConcreteSchwingerDyson}
suggested by the last exercise in~\cite{AGZ}
(which does not otherwise seem to be present in the literature in
detail) because
it makes connectedness of $\DDD_L$ a nonissue.
\end{Remark}

\begin{Remark}\label{RemarkCobblingBis}
(This is a continuation of the thread of remarks including Remarks \ref
{RemarkPrePreCobbling} and~\ref{RemarkPreCobbling}.)
If $(\SSS,L,\Theta,e)$ is a self-adjoint linearization of $f\in\Mat
_n({\CC\langle\Xbold\rangle})_\sa$,
then we have the simple formula
%
%
\begin{equation}\label{equationDeeplyRelevantBis}
S_{\mu_f}(z)=\tau_{\SSS,e}\bigl(G_L(\Theta+ze)\bigr)
\end{equation}
for $z\in\CC\setminus\supp\mu_f$, which is just a rewrite of
\eqref{equationPreCobblingOperator}.
This is one way---but not the only way---in which solutions of the SD
equation enter the proof of Theorem~\ref{TheoremMainResultBis}.
\end{Remark}

\subsection{Derivatives, symmetries and estimates}
We record some immediate consequences of the construction of
$G_L(\Lambda)$ for later use.
\subsubsection{Operator-theoretic representation of the derivative}
By means of the resolvent identity \eqref{equationResolventIdentity}
in infinitesimal form, one verifies that
%
%
\begin{eqnarray}
&&\DD[G_L](\Lambda;\zeta)
\nonumber
\\[-8pt]
\\[-8pt]
\nonumber
&&\qquad=\phi^\BF_\SSS\bigl(\bigl(L(\Xi)-1_{B(\HHH
)}\otimes\Lambda\bigr)^{-1}
(1\otimes\zeta)\bigl(L(\Xi)-1_{B(\HHH)}\otimes\Lambda\bigr)^{-1}\bigr),
\end{eqnarray}
for all $\zeta\in\SSS$.
\subsubsection{Symmetries}
Note that
%
%
\begin{eqnarray}\label{equationGConjSymm}
\Lambda&\in&\DDD_L\quad\Leftrightarrow\quad\Lambda^*\in\DDD
_{L^*}\quad\Rightarrow\quad
G_L(\Lambda)^*=G_{L^*}(\Lambda^*),\\
\label{equationGTransSymm}
\Lambda&\in&\DDD_L\quad\Leftrightarrow\quad\Lambda^\T\in\DDD_{L^\T
}\quad\Rightarrow\quad
G_L(\Lambda)^\T=G_{L^\T}(\Lambda^\T).
\end{eqnarray}
Relation \eqref{equationGConjSymm}
holds by the symmetry \eqref{equationExtensionStarSymmetry}
along with the observation that $*$ commutes with inversion.
Relation \eqref{equationGTransSymm}
can be verified by a straightforward calculation exploiting Lemma \ref
{LemmaCStarTee}
and relation \eqref{equationExtensionTransSymmetry}.

\subsubsection{Estimates}
For any $\SSS$-linear form $L$
and points $\Lambda,\Lambda_1,\Lambda_2\in\DDD_L$, we have estimates
%
%
\begin{eqnarray}\label{equationGLBound}
&&[\![G_L(\Lambda)]\!]\leq\bigl[\!\bigl[\bigl(L(\Xi)- 1_{B(\HHH)}\otimes \Lambda
\bigr)^{-1}\bigr]\!\bigr],\\
\label{equationGLBoundBis}
&&[\![G_L(\Lambda_1)-G_L(\Lambda_2)]\!]
\nonumber
\\
&&\qquad\leq[\![\Lambda_1-\Lambda_2]\!]
\bigl[\!\bigl[\bigl(L(\Xi)-
1_{B(\HHH)}\otimes\Lambda_1\bigr)^{-1}\bigr]\!\bigr]\\
&&\qquad\quad{}\times
\bigl[\!\bigl[\bigl(L(\Xi)- 1_{B(\HHH)}\otimes\Lambda_2\bigr)^{-1}\bigr]\!\bigr],\nonumber\\
\label{equationGLBoundTer}
&&[\![\DD[G_L](\Lambda)]\!]\leq[\![(L(\Xi)- 1_{B(\HHH )}\otimes
\Lambda)^{-1}]\!]^2,\\
\label{equationGLBoundQuad}
&&[\![G_L(\Lambda_1)-G_L(\Lambda_2)-\DD[G_L](\Lambda_2;\Lambda
_1-\Lambda_2)]\!]\nonumber\\
&&\qquad\leq[\![\Lambda_1-\Lambda_2]\!]^2
\bigl[\!\bigl[\bigl(L(\Xi)-
1_{B(\HHH)}\otimes\Lambda_1\bigr)^{-1}\bigr]\!\bigr]\\
&&\qquad\quad{}\times\bigl[\!\bigl[\bigl(L(\Xi)- 1_{B(\HHH)}\otimes\Lambda_2\bigr)^{-1}\bigr]\!\bigr]^2\nonumber
\end{eqnarray}
which follow directly from the resolvent identity \eqref
{equationResolventIdentity},
the iterated resolvent identity \eqref{equationResolventIdentityIterated},
the estimate \eqref{equationStateBound}
and the definitions.

\subsection{The secondary SD equation}\label
{subsectionSecondarySchwingerDyson}

We construct solutions of a secondary form of the
Schwinger--Dyson equation by using a variant of the secondary trick.
See Proposition~\ref{PropositionSecondarySchwingerDyson} below and
its proof.

\subsubsection{\texorpdfstring{The special function $G_{L_1,L_2}(\Lambda_1,\Lambda_2)$}
{The special function G L1, L2(Lambda1, Lambda2)}}
\label{subsubsectionGDoubleL}

Let $\SSS$ be a block algebra.
For $j=1,2$, let $L_j$ be an $\SSS$-linear form
and let $\Lambda_j\in\DDD_{L_j}$ be a point. We define
\[
G_{L_1,L_2}(\Lambda_1,\Lambda_2)
=\phi^{\BF}_{\SSS,\SSS}\bigl(\bigl(L_1(\Xi) - 1_{B(\HHH)}\otimes\Lambda
_1 \bigr)^{-1},
\bigl(L_2(\Xi) - 1_{B(\HHH)}\otimes\Lambda_2\bigr)^{-1}\bigr),
\]
where $\phi^{\BF}_{\SSS,\SSS}$
is the $\SSS$-bilinear extension of $\phi^{\BF}$. It is easy to see
that $G_{L_1,L_2}(\Lambda_1,\break\Lambda_2)$ depends
analytically on
$(\Lambda_1,\Lambda_2)$.
By Remark~\ref{RemarkCoffeePressBis}, we have
%
%
\begin{eqnarray}\label{equationGLLBound}
&&[\![G_{L_1,L_2}(\Lambda_1,\Lambda_2)]\!]
\nonumber
\\[-8pt]
\\[-8pt]
\nonumber
&&\qquad\leq\bigl[\!\bigl[\bigl(L_1(\Xi) - 1_{B(\HHH)}\otimes\Lambda_1 \bigr)^{-1}\bigr]\!\bigr]\bigl[\!\bigl[{\bigl(L_2(\Xi) - 1_{B(\HHH)}\otimes\Lambda_2\bigr)^{-1}}\bigr]\!\bigr],
\end{eqnarray}
which is an estimate analogous to \eqref{equationGLBound}.

\begin{Proposition}\label{PropositionSecondarySchwingerDyson}
Let $\SSS$ be a block algebra.
For $j=1,2$, let $L_j=\sum\Xbold_\ell\otimes a_{\ell j}$ be an $\SSS
$-linear form\vadjust{\goodbreak}
and let $\Lambda_j\in\DDD_{L_j}$ be a point. Then the \textit{secondary
SD equation}
%
%
\begin{eqnarray}\label{equationSecondaryMcGuffinBis}
&&G_{L_1,L_2}(\Lambda_1,\Lambda_2)
\nonumber
\\[-8pt]
\\[-8pt]
\nonumber
&&\qquad=
\Bigl(\Bigl(\Bigl(G_{L_1}(\Lambda_1)^{-1}\otimes G_{L_2}(\Lambda
_2)^{-1}-\sum a_{\ell1}
\otimes a_{\ell2}\Bigr)^{1\otimes\T}\Bigr)^{-1}
\Bigr)^{1\otimes\T}
\end{eqnarray}
holds. In particular, the expression on the right side is well-defined.
\end{Proposition}

It is worth noting as a consistency check that the expression
on the right side
remains invariant if we replace the transposition $\T$ by any other
transposition of $\SSS$.
\begin{pf} By Remark~\ref{RemarkPlotThickens}
it suffices to prove that
%
%
\begin{eqnarray}\label{equationSecondaryMcGuffin}
\zeta&=&G_{L_1}(\Lambda_1)^{-1}G_{L_1,L_2}(\Lambda
_1,\Lambda_2)^\bullet(\zeta)G_{L_2}(\Lambda_2)^{-1}
\nonumber
\\[-8pt]
\\[-8pt]
\nonumber
&&{}-\sum a_{\ell1} G_{L_1,L_2}(\Lambda_1,\Lambda_2)^\bullet(\zeta
)a_{\ell2}
\end{eqnarray}
holds for all $\zeta\in\SSS$.
Let $\MMM_2$ be a block algebra equipped with a standard basis
$\{e_{ij}\}_{i,j=1}^2$. Fix $\zeta\in\SSS$ arbitrarily and put
\[
\Lambda=
\Lambda_1\otimes e_{11}+\Lambda_2\otimes e_{22}+\zeta\otimes
e_{12}\in\SSS\otimes\MMM_2.
\]
Consider also the $\SSS\otimes\MMM_2$-linear form
\[
L=L_1\otimes e_{11}+L_2\otimes e_{22}.
\]
To compress notation put
\[
A_j=L_j(\Xi)-1_{B(\HHH)}\otimes\Lambda_j\in(\AAA\otimes\SSS
)^\times
\]
for $j=1,2$ and put
\[
A=L(\Xi)-1_{B(\HHH)}\otimes\Lambda\in\AAA\otimes\SSS\otimes
\MMM_2.
\]
In fact $A\in(\AAA\otimes\SSS\otimes\MMM_2)^\times$, and
more precisely
\[
A^{-1}= A_1^{-1}\otimes e_{11}+A_2^{-1}\otimes e_{22}
+\bigl(A_1^{-1}(1_\AAA\otimes\zeta)A_2^{-1}\bigr)\otimes e_{12}
\]
as one immediately verifies.
Thus by the trivial identity \eqref{equationBulletCuriosity},
we have
\[
G_L(\Lambda)=
G_{L_1}(\Lambda_1)\otimes e_{11}+G_{L_2}(\Lambda_2)\otimes e_{22}+
G_{L_1,L_2}(\Lambda_1,\Lambda_2)^\bullet(\zeta)\otimes e_{12}.
\]
By Proposition~\ref{PropositionConcreteSchwingerDyson}, the SD equation
\[
0=1_\SSS\otimes1_{\MMM_2}+\bigl(\Lambda+\Phi_L(G_L(\Lambda
))\bigr)G_L(\Lambda)
\]
is satisfied.
By expanding the right side in the form
$\ldots+b\otimes e_{12}+\ldots$ we find that
\begin{eqnarray*}
0&=&\bigl(\Lambda_1+\Phi_{L_1}(G_{L_1}(\Lambda_1))\bigr)G_{L_1,L_2}(\Lambda
_1,\Lambda_2)^\bullet(\zeta)\\
&&{} +\Bigl(\zeta+\sum a_{\ell1}G_{L_1,L_2}(\Lambda_1,\Lambda
_2)^\bullet(\zeta) a_{\ell2}\Bigr)G_{L_2}(\Lambda_2),
\end{eqnarray*}
which yields \eqref{equationSecondaryMcGuffin} after some
further manipulation which we omit.
\end{pf}

\begin{Remark} \label{RemarkBiasWellDefined}
Fix an $\SSS$-linear form $L$ and a point $\Lambda\in\DDD_L$.
Then we have
%
%
\begin{equation}\label{equationBiasWellDefined}
\DD[G_L](\Lambda)=G_{L,L}(\Lambda,\Lambda)^\bullet
\end{equation}
as one verifies by exploiting the infinitesimal form of the resolvent
identity \eqref{equationResolventIdentity}.
Note that the equation
\eqref{equationSecondaryMcGuffin}
in the case $(L_1,L_2,\Lambda_1,\Lambda_2)=(L,L,\Lambda,\Lambda)$
specializes to the equation \eqref{equationSchwingerDysonCorrection}
obtained through differentiation.
\end{Remark}

\begin{Remark} \label{RemarkBiasWellDefinedBis}
Fix an $\SSS$-linear form $L$ and a point $\Lambda\in\DDD_L$.
Let $\Psi_L$ be as in Definition~\ref{DefinitionPhiPsi}.
Recall that if $L=\sum\Xbold_\ell\otimes a_\ell$ is the Hamel
expansion of $L$
then $\Psi_L=\sum(-1)^\ell a_\ell^{\otimes2}$. Then we have
%
%
\begin{equation}\label{equationBiasWellDefinedBis}
\bigl((G_L(\Lambda)^{-1})^{\otimes2}-\Psi_L\bigr)^{-1}=G_{L,L^\T
}(\Lambda,\Lambda^\T)^{1\otimes\T}
\end{equation}
by the secondary SD equation \eqref{equationSecondaryMcGuffinBis}
in the case
$(L_1,L_2,\Lambda_1,\Lambda_2)=(L,L^\T,\Lambda,\Lambda^\T)$
along with the symmetry \eqref{equationGTransSymm}. In turn, we have
%
%
\begin{eqnarray}\label{equationGLBizarreBound}
\bigl[\!\bigl[\bigl((G_L(\Lambda)^{-1})^{\otimes2}-\Psi_L\bigr)^{-1}\bigr]\!\bigr]&\leq&[\![1\otimes
\T]\!]
\bigl[\!\bigl[\bigl(L(\Xi)-1_{B(\HHH)}\otimes\Lambda\bigr)^{-1}\bigr]\!\bigr]^2
\end{eqnarray}
by Remark~\ref{RemarkCoffeePressBis}, \eqref{equationGLLBound}, and
Lemma~\ref{LemmaCStarTee}.
\end{Remark}

\begin{Remark}\label{RemarkRaisonDetre}
Fix a SALT block design $(\SSS,L,\Theta,e)$ and a point $\Lambda\in
\DDD_L$.
We automatically have $\underline{\Lambda}+\Diamond_\SSS\in\DDD
_{\underline{L}}$
and
%
%
\begin{eqnarray}\label{equationRaisonDetre2}
\partial_1G_{\underline{L}}(\underline{\Lambda}+\Diamond_\SSS)
&=&\DD[G_L](\Lambda)\quad \mbox{and}
\nonumber
\\[-8pt]
\\[-8pt]
\nonumber
\partial_2 G_{\underline{L}}(\underline{\Lambda}+\Diamond_\SSS)
&=&\bigl((G_L(\Lambda)^{-1})^{\otimes2}-\Psi_L\bigr)^{-1}
\end{eqnarray}
by \eqref{equationUnderlineConstruction3} and
\eqref{equationUnderlineConstruction4} along with
\eqref{equationBiasWellDefined} and \eqref{equationBiasWellDefinedBis}.
\end{Remark}

\subsection{\texorpdfstring{The universal correction $\Bias_L^N$}{The universal correction Bias L N}}\label
{subsectionUniversalCorrection}
We first present a general construction needed for the proof of Theorem
\ref{TheoremMainResultBis}
which involves solutions of both the ``primary'' and secondary SD equations.
Then in Remark~\ref{RemarkCobblingTer} we specialize the construction
to produce our candidate
for the correction ${\mathrm{bias}}_f^N(z)$ figuring in Theorem~\ref
{TheoremMainResultBis}.
Throughout, $\SSS$ denotes a fixed block algebra.

\subsubsection{A tensor generalization of fourth cumulants}
Let $Y$ be any $\SSS$-valued random variable such that ${\Vert[\![Y]\!]
\Vert}_4<\infty$
and $\Ebold Y=0$. Let $Z$ be an independent copy of $Y$.
We define
\begin{eqnarray*}
\Cbold^{(4)}(Y)&=&\Ebold(Y^*\otimes Y\otimes Y^*\otimes Y)-\Ebold
(Y^*\otimes Y\otimes Z^*\otimes Z)\\
&&{}-\Ebold( Y^*\otimes Z\otimes
Z^*\otimes Y)
-\Ebold( Y^*\otimes Z\otimes Y^*\otimes Z)\in\SSS^{\otimes
4}.\nonumber
\end{eqnarray*}
%
\subsubsection{Shuffle notation}
For positive integers $k$, we define bilinear maps
\begin{eqnarray*}
&&[\cdot,\cdot]_k\dvtx \SSS^{\otimes k}\times\SSS^{\otimes k}
\rightarrow\SSS^{\otimes2k},\\
&&\hspace*{40pt}\qquad [x_1\otimes\cdots\otimes
x_{k},y_1\otimes\cdots\otimes y_k]_k=
x_1\otimes y_1\otimes\cdots\otimes x_k\otimes y_k,\hspace*{-40pt}
\\
&&\langle\cdot,\cdot\rangle_k\dvtx \SSS^{\otimes k}\times\SSS^{\otimes k}
\rightarrow\SSS,\\
&&\hspace*{89pt}\qquad\langle x_1\otimes\cdots\otimes x_{k},y_1\otimes\cdots\otimes
y_k\rangle_k=
x_1y_1\cdots x_k y_k.\hspace*{-89pt}
\end{eqnarray*}

\subsubsection{\texorpdfstring{Definition of $\Bias^N_L$}{Definition of Bias N L}}
Let $L=\sum\Xbold_\ell\otimes a_\ell$ be any self-adjoint \mbox{$\SSS$-linear} form.
Let $\Lambda\in\DDD_L$ be a point.
To abbreviate notation, we write
\begin{eqnarray*}
\Phi&=&\Phi_L\in B(\SSS), \qquad\Psi=\Psi_L\in\SSS^{\otimes2},\\
X^N&=&L(\Xi^N)=\sum\Xi^N_\ell\otimes a_\ell\in\Mat_N(\SSS)_\sa,
\\
G &=& G_L(\Lambda)\in\SSS^{\times},\qquad
G'=\DD[G_L](\Lambda)\in B(\SSS),\\
\check{G} &= &\bigl((G^{-1})^{\otimes2}-\Psi_L\bigr)^{-1}\in
(\SSS^{\otimes2})^\times.
\end{eqnarray*}
By Remark~\ref{RemarkBiasWellDefinedBis},
the object $\check{G}$ above is well-defined. We now define
%
%
\begin{eqnarray}
\label{equationUnwrapped}
\widehat{\Bias}_L^N(\Lambda)&=&
\langle[\Psi,\Psi]_2,[\check{G},G^{\otimes2}]_2\rangle_4-\Phi
(G)G
\nonumber\\
&&{}+\frac{1}{N}\sum_{i=1}^N \langle\Ebold
X^N(i,i)^{\otimes2},G^{\otimes2}\rangle_2
\nonumber
\\[-8pt]
\\[-8pt]
\nonumber
&&{}-
\frac{1}{N^{3/2}}\sum_{i=1}^N\langle\Ebold X^N(i,i)^{\otimes
3},G^{\otimes3}\rangle_3\\
&&{} +\frac{1}{N^2}\mathop{\sum_{
i,j=1}}_{
i\neq j}
^N\langle\Cbold^{(4)}(X^N(i,j)),G^{\otimes4}\rangle_4,\nonumber\\
\label{equationWrapped}
\Bias_L^N(\Lambda) &=& G'(\widehat{\Bias}_L^NG^{-1}).
\end{eqnarray}
The analytic functions
\[
\widehat{\Bias}_L^N,\Bias_L^N\dvtx \DDD_L\rightarrow\SSS
\]
thus defined we call the
\textit{unwrapped universal correction}
and \textit{universal correction}
indexed by $L$ and $N$, respectively.
We only define the former function to expedite certain
calculations---the latter function is the theoretically important one
with good symmetry properties.
It is a straightforward if tedious matter to verify that $\Bias_L^N$
commutes with the $C^*$-algebra involution just as $G_L$ does.
For a constant $c$ independent of $N$, $L$ and $\Lambda$, we have
%
%
\begin{equation}\label{equationVitalBiasBound}
\sup_N[\![\Bias^N_L(\Lambda)]\!]\leq c
\bigl[\!\bigl[\bigl(L(\Xi)-1_{B(\HHH)}\otimes\Lambda\bigr)^{-1}\bigr]\!\bigr]^5
\end{equation}
by estimates \eqref{equationGLBound}, \eqref{equationGLBoundTer}
and \eqref{equationGLBizarreBound} along with assumption \eqref
{equationPepper1}.

\begin{Remark}
Since we are long done with the discussion of Theorems~\ref
{TheoremPreMainResult} and~\ref{TheoremMainResult}
and are focused now on proving Theorem~\ref{TheoremMainResultBis},
we feel free to repurpose the letter $X$ in various ways as, for
example, in the construction
above of $\Bias^L_N$, and later in our discussion of the block Wigner
model. This should not cause confusion.
\end{Remark}

\begin{Remark}\label{RemarkCobblingTer}
(This is a continuation of the thread of remarks including Remarks \ref
{RemarkPrePreCobbling} and~\ref{RemarkCobblingBis}.)
If $(\SSS,L,\Theta,e)$ is a self-adjoint linearization of $f\in\Mat
_n({\CC\langle\Xbold\rangle})_\sa$,
then our candidate for the correction figuring in Theorem \ref
{TheoremMainResultBis}
is defined by the formula
%
%
\begin{equation}
{\mathrm{bias}}_f^N(z)=\tau_{\SSS,e}\bigl(\Bias_L^N(\Theta+ze)\bigr)
\end{equation}
for $z\in\CC\setminus\supp\mu_f$.
This is another distinct way that solutions of the SD equation enter
the proof of Theorem~\ref{TheoremMainResultBis}.
\end{Remark}

\section{Approximation of solutions of the Schwinger--Dyson
equation}\label{sectionMagSetup}
We refine a powerful idea from~\cite{HST}
concerning the approximation of solutions of the Schwinger--Dyson equation.
See Lemma~\ref{LemmaKeyObservation} below
for a short paraphrase of that idea in a simplified geometry.
See Proposition~\ref{PropositionModifiedKahuna} below for the main
result of this section,
which is an estimate tailored to the proof of Theorem~\ref
{TheoremMainResultBis}.

\subsection{SD tunnels}

\begin{Definition}\label{DefinitionTunnels}
Suppose we are given
\begin{itemize}
\item a solution $G\dvtx \DDD\rightarrow\SSS$ of the SD equation with
covariance map $\Phi\in B(\SSS)$,
\item a point $\Lambda_0\in\DDD$ and
\item(finite) constants $\Tfrak>0$ and $\Gfrak\geq1$.
\end{itemize}
Put
\[
\TTT=\{\Lambda_0+\ii t 1_\SSS+\zeta\vert t\in[0,\infty) \mbox
{ and } \zeta\in\SSS\mbox{ s.t. }
[\![\zeta]\!]\leq1/\Gfrak\}.
\]
Suppose that the following conditions hold:
%
%
\begin{eqnarray}
\label{equationMag1}
\TTT&\subset&\DDD,\\
\label{equationMag2}
\sup_{\Lambda\in\TTT}[\![G(\Lambda)]\!]&\leq&\Gfrak,\\
\label{equationMag25}
\sup_{\Lambda\in\TTT}[\![\Dbold[G](\Lambda)]\!]&\leq&\Gfrak
^2,\\
\label{equationMag3}
\mathop{\sup_{
\Lambda,\Lambda'\in\TTT}}_{
\mathrm{s.t.}\ \Lambda\neq\Lambda'}
\frac{[\![G(\Lambda)-G(\Lambda')]\!]}{[\![\Lambda-\Lambda ']\!]}
&\leq&
\Gfrak^2,\\
\label{equationMag4}
\mathop{\sup_{
\Lambda,\Lambda'\in\TTT}}_{
\mathrm{s.t.}\ \Lambda\neq\Lambda'}
\frac{[\![G(\Lambda)-G(\Lambda')-\DD[G](\Lambda';\Lambda -\Lambda
')]\!]}{[\![\Lambda-\Lambda']\!]^2}&\leq&\Gfrak^3,\\
\label{equationMag5}
\sup_{\Lambda\in\TTT}[\![G(\Lambda+\ii\Tfrak1_\SSS)]\!]&\leq&
\frac{1}{2(1+[\![\Phi]\!])}.
\end{eqnarray}
In this situation we say that the collection
$(G\dvtx \DDD\rightarrow\SSS,\Phi,\Lambda_0,\Tfrak,\Gfrak)$
is a \textit{Schwinger--Dyson (SD) tunnel}.
\end{Definition}

\begin{Remark}
If $(G\dvtx \DDD\rightarrow\SSS,\Phi,\Lambda_0,\Tfrak,\Gfrak)$ is an
SD tunnel,
then
for every $t\in[0,\infty)$, so is
$(G\dvtx \DDD\rightarrow\SSS,\Phi,\Lambda_0+\ii t1_\SSS,\Tfrak,\Gfrak)$.
\end{Remark}

\begin{Remark}\label{RemarkCapsules}
All examples of SD tunnels needed for the proof of Theorem~\ref
{TheoremMainResultBis}
arise as follows. Let $(\SSS,L,\Theta,e)$ be a SALT block design.
Let $c_0$, $c_1$, $c_2$ and $\Tfrak$ be be the constants from
Definition~\ref{DefinitionSALT}.
Put
%
%
\begin{equation}\label{equationPreTypicalTunnel}
\Gfrak(z)=c_0\bigl(1+[\![L(\Xi)]\!]\bigr)^{c_1}(1+1/\Im z)^{c_2}
\end{equation}
for $z\in{\mathfrak{h}}$.
We claim that the collection
%
%
\begin{equation}\label{equationTypicalTunnel}
\bigl(G_L\dvtx \DDD_L\rightarrow\SSS,\Phi_L,\Theta+ze,\Tfrak,\Gfrak(z)\bigr)
\end{equation}
is an SD tunnel for each fixed $z\in{\mathfrak{h}}$. To prove the claim,
arbitrarily fix $z\in{\mathfrak{h}}$, $t\in[0,\infty)$ and $\zeta
\in\SSS$ such that $[\![\zeta]\!]\leq1/\Gfrak(z)$,
and put
\[
\Lambda=\Theta+ze+\ii t1_\SSS+\zeta.
\]
We then have
\[
\Lambda\in\DDD_L\quad \mbox{and}\quad
\bigl[\!\bigl[\bigl(L(\Xi)-1_{B(\HHH)}\otimes\Lambda\bigr)^{-1}\bigr]\!\bigr]\leq\cases{
\Gfrak(z),&\quad $\mbox{in general,}$\vspace*{2pt}\cr
\displaystyle\frac{1}{2(1+[\![\Phi]\!])},&\quad$\mbox{for $t\geq\Tfrak$,}$}
\]
by Definition~\ref{DefinitionSALT}
and Remark~\ref{RemarkCutoff}.
In particular, \eqref{equationTypicalTunnel} satisfies \eqref
{equationMag1} for each fixed $z\in{\mathfrak{h}}$.
In turn, it follows by \eqref{equationGLBound} that
\eqref{equationTypicalTunnel} satisfies
\eqref{equationMag2} and \eqref{equationMag5} for each fixed $z\in
{\mathfrak{h}}$.
Given also $\Lambda'\in\DDD_L$ with ``primed'' variables,
we have
\begin{eqnarray*}
[\![G_L(\Lambda)-G_L(\Lambda')]\!]&\leq&[\![\Lambda-\Lambda']\!]\Gfrak(z)\Gfrak(z'),\\{}
[\![{\DD[G_L](\Lambda)}]\!]&\leq&\Gfrak{(z)}^2,\\{}
[\![G_L(\Lambda)-G_L(\Lambda')-\DD[G_L](\Lambda';\Lambda -\Lambda
')]\!]&\leq&{[\![\Lambda-\Lambda']\!]}^2\Gfrak(z)\Gfrak{(z')}^2
\end{eqnarray*}
by \eqref{equationGLBoundBis}, \eqref{equationGLBoundTer}
and \eqref{equationGLBoundQuad}, respectively.
It follows that \eqref{equationTypicalTunnel} also satisfies
\eqref{equationMag2}, \eqref{equationMag25} and \eqref
{equationMag3} for each fixed $z\in{\mathfrak{h}}$.
The claim is proved.
We note also that we have a bound
%
%
\begin{equation}\label{equationVitalBiasBoundBis}
[\![\Bias_L^N(\Lambda)]\!]\leq c \Gfrak(z)^5
\end{equation}
for a constant $c$ independent of $L$, $N$ and $z$ by \eqref
{equationVitalBiasBound}.
This last estimate turns out to be the crucial point for proving the
bound \eqref{equationMainResultBis3}
asserted in Theorem~\ref{TheoremMainResultBis}.
\end{Remark}

\begin{Remark}\label{RemarkRedundancy}
Definition~\ref{DefinitionTunnels} is not particularly delicate or economical.
Indeed, conditions \eqref{equationMag1}, \eqref{equationMag2} and
\eqref{equationMag5}
alone imply that $(G\dvtx \DDD\rightarrow\SSS,\Phi,\Lambda_0, \Tfrak
,c\Gfrak)$ is an SD tunnel,
where $c>1$ is an absolute constant. However, for the present purpose,
no advantage is gained by
reformulating Definition~\ref{DefinitionTunnels} in more economical
fashion since all the properties
\eqref{equationMag1}--\eqref{equationMag5} are needed to prove
Proposition~\ref{PropositionModifiedKahuna} below.
\end{Remark}

\subsection{The key lemma}\label{subsectionLemmasEtc}
Before working out our main estimate, we first prove a simple lemma to
explain the mechanism by which SD tunnels control errors. The lemma\vadjust{\goodbreak} captures
a key idea of~\cite{HST} but works with a simpler geometry. The lemma
uses only the first
and last of the defining conditions of an SD tunnel.

\subsubsection{Setup for the key lemma}
\begin{itemize}
\item Let $(G\dvtx \DDD\rightarrow\SSS,\Phi,\Lambda_0,\Tfrak,\Gfrak)$
satisfy conditions
\eqref{equationMag1} and \eqref{equationMag5}
of the definition of an SD tunnel.
\item Let $F=(t\mapsto F_t)\dvtx [0,\Tfrak]\rightarrow\SSS$ be a
continuous function.
\end{itemize}
For $t\in[0,\Tfrak]$ we put
\[
\Lambda_t=\Lambda_0+\ii t1_\SSS,\qquad G_t = G(\Lambda_t),\qquad
E_t=1_\SSS+\bigl(\Lambda_t+\Phi(F_t)\bigr)F_t.
\]
Note that we have $\Lambda_t\in\DDD$ by definition of an SD tunnel
and hence $G_t$ is well-defined.
In turn we define constants
\[
\Cfrak_0=2(1+[\![\Phi]\!]),\qquad \Ffrak= 1\vee\sup_{t\in
[0,\Tfrak]}[\![F_t]\!],\qquad \Afrak=\sup_{t\in[0,\Tfrak]}[\![E_t]\!].
\]
The quantity $\Afrak$ is a natural measure of the failure of $F$ to
satisfy the SD equation.
We emphasize that we assume nothing of the function $F$ beyond continuity.

\begin{Lemma}\label{LemmaKeyObservation}
If
%
%
\begin{eqnarray}
\label{equationKey1}
\Cfrak_0\Gfrak\Ffrak\Afrak&<&1\quad \mbox{and}\\
\label{equationKey2}
[\![F_\Tfrak]\!]&<&1,
\end{eqnarray}
then for every $t\in[0,\Tfrak]$, the inverse $H_t=-(\Lambda_t+\Phi
(F_t))^{-1}$
exists,
%
%
\begin{eqnarray}
\label{equationKey5}
[\![H_t]\!] &\leq& 2[\![F_t]\!],\\
\label{equationKey3}
[\![\Phi(H_tE_t)]\!]&\leq&1/\Gfrak,\qquad \mbox{(hence) }
\Lambda_t-\Phi(H_tE_t)\in\DDD\quad\mbox{and} \\
\label{equationKey4}
H_t-H_tE_t&=& F_t=G\bigl(\Lambda_t-\Phi(H_tE_t)\bigr)-H_tE_t.
\end{eqnarray}
\end{Lemma}
\begin{pf}
Fix $t\in[0,\Tfrak]$ arbitrarily.
Hypothesis \eqref{equationKey1} implies that $[\![E_t]\!]\leq1/2$.
By Lemma~\ref{LemmaNeumannExpansion}
it follows that $H_t$ is well-defined
and satisfies \eqref{equationKey5}.
Then claim~\eqref{equationKey3} holds
by \eqref{equationMag1}, \eqref{equationKey1}
and \eqref{equationKey5}.
It remains only to prove claim \eqref{equationKey4},
and since the first equality in \eqref{equationKey4} holds by
definition of $H_t$,
we have only to prove the second equality.
By the Weierstrass Approximation theorem,
we may assume that $F$ depends polynomially and \textit{a fortiori}
analytically on $t$. Put
\[
\widehat{H}_t=G\bigl(\Lambda_t-\Phi(H_tE_t)\bigr)\quad \mbox{and} \quad\widehat
{F}_t=\widehat{H}_t-H_tE_t.
\]
Note that $\widehat{F}_t$ depends analytically on $t$.
It is enough to prove $F_t\equiv\widehat{F}_t$.
In any case, since $G$ satisfies the SD equation with covariance map
$\Phi$, we have
\[
1_\SSS+\bigl(\Lambda_t-\Phi(H_tE_t)+\Phi(\widehat{H}_t)\bigr)\widehat
{H}_t=1_\SSS+\bigl(\Lambda_t+\Phi(\widehat{F}_t)\bigr)\widehat{H}_t=0
\]
and hence $\widehat{H}_t=-(\Lambda_t+\Phi(\widehat{F}_t))^{-1}$.
We thus have
\[
F_t-\widehat{F}_t=H_t-\widehat{H}_t
=
H_t
\Phi(F_t-\widehat{F}_t)\widehat{H}_t
=
H_t\Phi(F_t-\widehat{F}_t)G\bigl(\Lambda_t-\Phi(H_tE_t)\bigr),\vadjust{\goodbreak}
\]
where at the second step we use the resolvent identity
\eqref{equationResolventIdentity}.
Finally, by \eqref{equationMag5}, \eqref{equationKey2},
\eqref{equationKey5} and \eqref{equationKey3} we have
\[[\![H_\Tfrak]\!][\![\Phi]\!]\bigl[\!\bigl[G\bigl(\Lambda_\Tfrak-\Phi (H_\Tfrak
E_\Tfrak)\bigr)\bigr]\!\bigr]<1,
\]
hence the difference $F_t-\widehat{F}_t$
vanishes identically for $t$ near $\Tfrak$
and hence $F_t\equiv\widehat{F}_t$ by analytic continuation.
\end{pf}

\begin{Remark}
We work out the simplest concrete example of the phenomenon described
by the lemma.
Let $\sigma\dvtx \CC\setminus[-2,2]\rightarrow\CC$ be the Stieltjes
transform of the semicircle law. As is well-known,
$\sigma(z)$ is the unique bounded analytic solution of the equation
$1+(z+\sigma(z))\sigma(z)=0$
in the domain $\CC\setminus[-2,2]$.
Now fix $z_0\in{\mathfrak{h}}$ arbitrarily.
It is easy to see that
\[
(G\dvtx \DDD\rightarrow\SSS,\Phi,\Lambda_0,\Tfrak,\Gfrak)=
\Biggl(\sigma\dvtx \CC\setminus[-2,2]\rightarrow\CC,1,z_0,4,1\vee\frac
{2}{\Im z_0}\Biggr)
\]
satisfies conditions \eqref{equationMag1} and \eqref{equationMag5}
of the definition of an SD tunnel.
Now fix a continuous function
$(t\mapsto F_t)\dvtx [0,4]\rightarrow\CC$ and put $\Ffrak=1\vee\sup
_{t\in[0,4]}|F_t|$.
In turn put $E_t=1+(z_0+\ii t+F_t)F_t$ for $t\in[0,4]$
and $\Afrak=\sup_{t\in[0,4]}|E_t|$.
Assume that
$|F_{4}|<1$ and $\Ffrak\Afrak<\frac{1}{4}(1\wedge\frac{\Im
z_0}{2})$.
Then by Lemma~\ref{LemmaKeyObservation} we have
$H_0=-(z_0+F_0)^{-1}\neq\infty$, $|H_0|\leq2 |F_0|$, $|H_0E_0|\leq
1\wedge\frac{\Im z_0}{2}$
and finally
$F_0=\sigma(z_0-H_0E_0)-H_0E_0$. This last equation is at first glance
a bit strange but in fact the strategy of writing $\rho(z)=\sigma
(z-\delta)-\delta$ to estimate $\rho(z)-\sigma(z)$ has long been in
use. See, for example,~\cite{BaiRate}, equation~4.11.
\end{Remark}

\begin{Remark}\label{RemarkShoot}
Equation \eqref{equationKey4} is not an obvious target to shoot for!
But once~\eqref{equationKey4} is written down, it is clear that it
offers excellent opportunities
for systematically estimating the difference $[\![F_0-G(\Lambda_0)]\!]$.
This surprising and powerful idea we learned
from~\cite{HST}. The importance and utility of this idea cannot be
overestimated.
\end{Remark}

\subsection{The tunnel estimates}\label{subsectionChocolate}
We now use Lemma~\ref{LemmaKeyObservation} to obtain an estimate
in terms of parameters over which we will be able to gain good control.
In particular, the estimate is designed to take
advantage of Remark~\ref{RemarkYCapsuleComplement} above.
\subsubsection{Setup for the tunnel estimates}
\begin{itemize}
\item Let $(G\dvtx \DDD\rightarrow\SSS,\Phi,\Lambda_0,\Tfrak,\Gfrak)$
be an SD tunnel.
\item Let $\Lfrak\in[1,\infty)$ be a constant.
\item Let $F=(t\mapsto F_t)\dvtx [0,\infty)\rightarrow\SSS$ be a
Lipschitz continuous function
with Lipschitz constant bounded by $\Lfrak$ and satisfying
$\sup_{t\in[\Tfrak,\infty)}[\![F_t]\!]\leq1/2$.
\end{itemize}
For $t\in[0,\infty),$ we put
\[
\Lambda_t=\Lambda_0+\ii t1_\SSS,\qquad
E_t=1_\SSS+\bigl(\Lambda_t+\Phi(F_t)\bigr)F_t,\vadjust{\goodbreak}
\]
and we define constants
\[
\Cfrak=99 e^{2\Tfrak}(1+[\![\Lambda_0]\!]+[\![\Phi]\!]),\qquad
\Efrak=\frac{1}{2}[\![E_0]\!]+\frac{1}{2}\int_0^\infty[\![E_t]\!]e^{-t}\,dt.
\]
The integral converges since $[\![E_t]\!]$ has at worst linear growth
as $t\rightarrow\infty$.
\begin{Proposition}\label{PropositionModifiedKahuna}
Data, notation and assumptions are as above.
We have
%
%
\begin{eqnarray}
\label{equationModifiedKahuna1}
[\![F_0-G(\Lambda_0)]\!]& \leq &(\Cfrak\Gfrak\Lfrak
)^{6}(\Efrak+\Efrak^2),\\
\label{equationModifiedKahuna2}
[\![F_0+\DD[G](\Lambda_0;E_0G(\Lambda _0)^{-1})-G(\Lambda_0)]\!]&
\leq &(\Cfrak\Gfrak\Lfrak
)^{12}(\Efrak^2+\Efrak^4).
\end{eqnarray}
\end{Proposition}

The exponents of $\Cfrak$, $\Gfrak$ and $\Lfrak$ are of no
importance in the sequel.
They could be replaced by any larger absolute constants without
disturbing later arguments.
Only the exponents of $\Efrak$ will be important.
\begin{pf*}{Proof of Proposition~\ref{PropositionModifiedKahuna}}
In anticipation of applying Lemma~\ref{LemmaKeyObservation},
we put
%
%
\begin{eqnarray}\label{equationFfrakCrude}
\Ffrak&=&1\vee\sup_{t\in[0,\Tfrak]}
[\![F_t]\!]=1\vee\sup_{t\in[0,\infty)}
[\![F_t]\!]\leq\sqrt{\Cfrak}\Lfrak,\\
\label{equationAfrakCrude}
\Afrak&=&\sup_{t\in[0,\Tfrak]}[\![E_t]\!]\leq\Cfrak\Ffrak^2.
\end{eqnarray}
We have also noted here some crude bounds needed later. Now write
\begin{eqnarray*}
G_0&=&G(\Lambda_0),\qquad G'_0=\Dbold[G](\Lambda_0),\\
V_0&=& G'_0(E_0G_0^{-1})=G'_0(\Phi(G_0E_0))+G_0E_0
\end{eqnarray*}
in order to abbreviate notation.
The last equality above is an instance of~\eqref
{equationSchwingerDysonCorrection}.

We now claim that
%
%
\begin{eqnarray}
\label{equationMcGuffin1}
[\![F_0-G_0]\!]& \leq &\Cfrak\Gfrak^2\Ffrak
([\![E_0]\!]
+\one_{\Cfrak\Gfrak\Ffrak\Afrak\geq1}),\\
\label{equationMcGuffin2}
[\![F_0+V_0-G_0]\!]& \leq &\Cfrak^2\Gfrak
^5\Ffrak^2
([\![E_0]\!]^2+\one_{\Cfrak\Gfrak\Ffrak\Afrak\geq1}).
\end{eqnarray}
If $\Cfrak\Gfrak\Ffrak\Afrak\geq1$, then
crude estimates based on the definition of an SD tunnel
along with the bound \eqref{equationAfrakCrude} suffice.
We may therefore assume without loss of generality
that $\Cfrak\Gfrak\Ffrak\Afrak<1$,
in which case
the hypotheses of Lemma~\ref{LemmaKeyObservation} are fulfilled.
Thus it follows via \eqref{equationMag3}, \eqref{equationKey5},
\eqref{equationKey3}
and \eqref{equationKey4} that
\[[\![F_0-G_0]\!]\leq[\![H_0-G_0]\!]+[\![H_0E_0]\!]
\leq\Gfrak^2(2[\![\Phi]\!]\Ffrak[\![E_0]\!])+2\Ffrak[\![E_0]\!].
\]
Thus, the claim \eqref{equationMcGuffin1} is proved.
To prove \eqref{equationMcGuffin2}, we begin by noting the identity
\begin{eqnarray*}
F_0+V_0-G_0&=&
G\bigl(\Lambda_0-\Phi(H_0E_0)\bigr)-G_0+G'_0(\Phi(H_0E_0))\\
&&{}+(G_0-H_0)E_0+G'_0\bigl(\Phi\bigl((G_0-H_0)E_0\bigr)\bigr)
\end{eqnarray*}
derived from \eqref{equationKey4}.
Then, reasoning as in the proof of \eqref{equationMcGuffin1},
but now also using~\eqref{equationMag4} and \eqref{equationMag5},
we find that
\begin{eqnarray*}
&&[\![F_0+V_0-G_0]\!]\\
&&\qquad\leq \Gfrak^3(2[\![\Phi]\!]\Ffrak[\![E_0]\!])^2+\Gfrak
^2(2[\![\Phi]\!]\Ffrak[\![E_0]\!])
[\![E_0]\!]\\
&&\qquad\quad{}+\Gfrak^3[\![\Phi]\!](\Gfrak^2(2[\![\Phi]\!]\Ffrak[\![E_0]\!]))[\![E_0]\!].
\end{eqnarray*}
Thus, the claim \eqref{equationMcGuffin2} is proved.

We next claim that
%
%
\begin{equation}\label{equationLipschitzTriangle}
\Afrak
\leq\sqrt{\Cfrak\Ffrak\Lfrak}\bigl(\sqrt{\Efrak}+\Efrak\bigr).
\end{equation}
To prove the claim, consider the function
$b(t)=e^{-t}[\![E_t]\!]$
defined for $t\in[0,\infty)$.
Since $b$ is continuous
and tends to $0$ at infinity, $b$ achieves its maximum
at some point $t_0\in[0,\infty)$. Clearly, we have
\[
b(t_0)e^{\Tfrak}\geq\Afrak, 2\Efrak\geq\int_0^\infty b(t)\,dt.
\]
Now fix $t>t_0$ arbitrarily.
We have
\begin{eqnarray*}
|b(t_0)-b(t)|&\leq&e^{-t_0}[\![E_{t_0}]\!](1-e^{t_0-t})+e^{-t}[\![E_{t_0}-E_t]\!]\\
&\leq&
\bigl(b(t_0)+(1+[\![\Phi]\!]\Lfrak)\Ffrak+
([\![\Lambda_0]\!]+e^{-t}t+[\![\Phi]\!]\Ffrak)\Lfrak
\bigr)|t_0-t|\\
&\leq&\biggl(b(t_0)+\frac{\Cfrak\Ffrak\Lfrak}{8e^{2\Tfrak}}
\biggr)|t_0-t|.
\end{eqnarray*}
Thus, there exists a right triangle with altitude $b(t_0)$
and base of length
$\frac{b(t_0)}{b(t_0)+{\Cfrak\Ffrak\Lfrak}/(8e^{2\Tfrak})}$
under the graph of $b$.
Now in general for $K_1,x\geq0$ and $K_2>0$
we have
\[
4K_1\geq\frac{x^2}{x+K_2}\quad\Rightarrow \quad x\leq\sqrt{8K_1K_2}+8K_1.
\]
The claim \eqref{equationLipschitzTriangle} now follows after some
further trivial manipulations which we omit.



Finally, by combining \eqref{equationFfrakCrude}, \eqref{equationMcGuffin1},
\eqref{equationMcGuffin2} and \eqref{equationLipschitzTriangle}
we obtain bounds
\begin{eqnarray*}
[\![F_0-G_0]\!]&\leq& \Cfrak^{3/2}\Gfrak^2\Lfrak\bigl(
2\Efrak+\one_{\Cfrak^2 \Gfrak^{3/2} \Lfrak^{3/2}(\sqrt{\Efrak
}+\Efrak)\geq1}\bigr),\\{}
[\![F_0+V_0-G_0]\!]
&\leq&\Cfrak^3\Gfrak^5\Lfrak^2
\bigl(4\Efrak^2+\one_{\Cfrak^2 \Gfrak^{3/2} \Lfrak^{3/2}(\sqrt{\Efrak
}+\Efrak)\geq1}\bigr),
\end{eqnarray*}
whence the result after using Chebyshev bounds and simplifying
brutally.~%
\end{pf*}

\section{Matrix identities}\label{sectionUrMatrixIdentities}
Throughout this section, we fix a block algebra~$\SSS$.
Working in a purely algebraic setting, we build up a catalog of identities
satisfied by finite chunks of an infinite matrix with entries in $\SSS$.
The identities are chosen to illuminate the structure of random matrices
of the form \eqref{equationTypicalBlockRandom}
and
are a further contribution to our stock of tools for concentration.
All the identities\vadjust{\goodbreak} derived here are meaningful in the case $\SSS=\CC$,
in which case many of these identities are familiar from the study of
resolvents of Wigner matrices.

\subsection{An ad hoc infinite matrix formalism}
\label{subsectionInfiniteMatrices}
When we write $\Mat_{k\times\ell}(\SSS)$,
we now allow $k$ or $\ell$ or both to be infinite,
in which case we mean for the corresponding matrix indices to range over
all positive integers.
Addition, multiplication and adjoints of (possibly) infinite matrices
are defined as before,
although we never attempt to multiply such matrices unless one of them
has only
finitely many nonzero entries. For each integer $N>0$, let $\III_N$
denote the family of nonempty subsets of the set $\{1,\ldots,N\}$.
Given a finite nonempty set
$I=\{i_1<\cdots<i_k\}$ of positive integers, let $\fbold_I\in\Mat
_{k\times\infty}(\SSS)$
and ${\mathbf{e}}_I\in\Mat_\infty(\SSS)$ be defined by
\[
\fbold_I(i,j)=\sum_{\alpha=1}^k\one_{(i,j)=(\alpha,i_\alpha
)}1_\SSS \quad\mbox{and}\quad
{\mathbf{e}}_I(i,j)=\sum_{\alpha=1}^k \one_{(i,j)=(i_\alpha
,i_\alpha)}1_\SSS,
\]
respectively.
Note that
$\fbold_I\fbold_I^*=\Ibold_{|I|}\otimes1_\SSS$ and $\fbold
_I^*\fbold_I={\mathbf{e}}_I$,
where $|I|$ denotes the cardinality of $I$. Note that for all $A\in
\Mat_\infty(\SSS)$
and finite sets $I$ and $J$ of positive integers,
the finite matrix $\fbold_IA\fbold_J^*\in\Mat_{|I|\times|J|}(\SSS
)$ is the result of
striking all rows of~$A$ with indices not in $I$ and all columns of $A$
with indices not in $J$.
Thus, the matrices $\fbold_I$ allow us to pick out finite chunks of a
matrix $A\in\Mat_\infty(\SSS)$
and to use the familiar rules of matrix algebra itself to manipulate
the chunks.
For $A\in\Mat_\infty(\SSS)$ with only finitely many nonzero entries,
we define
$\trace_\SSS A=\sum_i A(i,i)$.
For such $A$, we also define
$[\![A]\!]=[\![\fbold_I A\fbold_I^*]\!]$ for any
finite set $I$ of positive integers such that
${\mathbf{e}}_IA{\mathbf{e}}_I=A$, which is independent of $I$.
For each $\zeta\in\SSS$,
let $\Ibold_\infty\otimes\zeta\in\Mat_\infty(\SSS)$
denote the infinite diagonal matrix with diagonal entries $\zeta$.

\subsection{The setup for studying matrix identities}\label
{subsectionRecipes}

\subsubsection{Data and assumption}
We fix a triple $(X,\Lambda,\Phi)$
where
\begin{itemize}
\item$X\in\Mat_\infty(\SSS)$,
\item$\Lambda\in\SSS$ and
\item$\Phi\in B(\SSS)$,
\end{itemize}
subject to the condition
%
%
\begin{equation}
\label{equationUrMatrixIdentitySetup}
\fbold_I\biggl(\frac{X}{\sqrt{N}}-\Ibold_\infty\otimes\Lambda
\biggr)\fbold_I^*\in\GL_{|I|}(\SSS)
\qquad\mbox{for $N$ and $I\in\III_N$.}
\end{equation}
Here and below $N$ is understood to range over the positive integers.
Below we will define and analyze various functions of the triple
$(X,\Lambda,\Phi)$,
calling them \textit{recipes}.

\begin{Remark}\label{RemarkCobbling}
All triples $(X,\Lambda,\Phi)$ needed to prove Theorem \ref
{TheoremMainResultBis} arise as follows.
Let $(\SSS,L,\Theta,e)$ be any SALT block design.
Let $\bigcup L(\Xi^N)\in\Mat_\infty(\SSS)_\sa$ denote the\vadjust{\goodbreak}
infinite matrix gotten by cobbling together the matrices $L(\Xi^N)\in
\Mat_N(\SSS)$ for varying $N$ using assumption \eqref{equationPepper4}.
Let $\zbold$ be as in Theorem~\ref{TheoremMainResultBis}.
Let $\tbold$ be a real random variable independent of $\sigma(\FFF
,\zbold)$
which with probability $1/2$ is concentrated at the origin and with
probability $1/2$
is standard exponential. (The motivation for using the random variable
$\tbold$ comes from Proposition~\ref{PropositionModifiedKahuna} above.)
Let $\Phi_L\in B(\SSS)$ be as in Definition~\ref{DefinitionPhiPsi}.
Then the triple
%
%
\begin{equation}\label{equationPreDeeplyRelevant}
(X,\Lambda,\Phi)=\Bigl(\bigcup L(\Xi^N),\Theta+\zbold e+\ii\tbold
1_\SSS,\Phi_L\Bigr)
\end{equation}
satisfies \eqref{equationUrMatrixIdentitySetup} with probability $1$
by Remark~\ref{RemarkYCapsuleComplement}.
\end{Remark}

\subsubsection{The first group of recipes}
For $N$ and $I\in\III_N$, we define
\begin{eqnarray*}
R^N_I&=&
\fbold_I^*\biggl(\fbold_I\biggl(\frac{X}{\sqrt{N}}-\Ibold_\infty
\otimes\Lambda\biggr)\fbold_I^*\biggr)^{-1}\fbold_I\in\Mat
_\infty(\SSS), \\
 F^N_I &=& \frac{1}{N}\trace_\SSS R_I^N\in\SSS
,\\
T_I^{N}&=&\biggl(\zeta\mapsto\frac{1}{N}
\sum_{i,j\in I}R_I^N(i,j)\zeta R_I^N(j,i)\biggr)\in B(\SSS),\\
U_I^N&=&\frac{1}{N}\sum_{i,j\in I} R_I^N(i,j)^{\otimes2}\in\SSS
^{\otimes2}.
\end{eqnarray*}
Note that $R^N_I$ is well-defined by assumption \eqref
{equationUrMatrixIdentitySetup}.
For $N$ put
\[
\III_N^{(2)}=\{(I,J)\in\III_N\times\III_N\vert J\subset I,
I\setminus J\in\III_N, |J|\leq2\}.
\]
For $N$ and $(I,J)\in\III_N^{(2)}$ put
\[
R_{I,J}^N=\fbold_JR_I^N\fbold_J^*\in\Mat_{|J|}(\SSS).
\]
The recipes in the first group do not depend on $\Phi$, whereas the
remaining recipes
we are about to define do depend on $\Phi$.

\subsubsection{Recipes of the second group}
For $N$ and $I\in\III_N$ put
\begin{eqnarray*}
E_I^{N}&=&1_\SSS+\bigl(\Lambda+\Phi(F^N_I)\bigr)F^N_I\in\SSS,\\
H_I^{N}& = &\cases{
-\bigl(\Lambda+\Phi(F_I^N)\bigr)^{-1}\in\SSS^\times,&\quad{\mbox{if $[\![E^N_I]\!]<
1/2$,}} \vspace*{2pt}\cr
0\in\SSS, &\quad {\mbox{if $[\![E^N_I]\!]\geq1/2$.}}}
\end{eqnarray*}
Note that $H^N_I$ is well-defined by
Lemma~\ref{LemmaNeumannExpansion}.
For $N$, $(I,J)\in\III_N^{(2)}$ and $j_1,j_2\in J$, we define
\begin{eqnarray*}
H_{I,J}^N&=&\Ibold_{|J|}\otimes H_{I\setminus J}^N\in\Mat_{|J|}(\SSS
),\\
\frac{Q_{I,J}^N}{\sqrt{N}}&=&-\frac{\fbold_JX\fbold_J^*}{\sqrt{N}}+
\frac{\fbold_J XR_{I\setminus J}^NX\fbold_{J}^*}{N}-\Ibold
_{|J|}\otimes\Phi(F_{I\setminus J}^N)\in\Mat_{|J|}(\SSS),\\
Q^N_{I,J,j_1,j_2}&=&\fbold_{j_1}\fbold_J^*Q^N_{I,J}\fbold_J\fbold
_{j_2}^*\in\SSS,\\
\frac{P_{I,J}^{N}}{\sqrt{N}}
&=&\biggl(A\mapsto\frac{1}{N}\trace_\SSS(R_{I\setminus J}^{N}
X\fbold_J^*
A\fbold_JXR^{N}_{I\setminus J})-T_{I\setminus J}^{N}\circ\Phi\circ
\trace_\SSS(A)\biggr)\\
&\in& B\bigl(\Mat_{|J|}(\SSS),\SSS\bigr),\\
P^N_{I,J,j_1,j_2}&=&\bigl(\zeta\mapsto P^N_{I,J}(\fbold_J\fbold
_{j_1}^*\zeta\fbold_{j_2}\fbold_J^*)\bigr)\in B(\SSS),\\
\Delta_{I,J}^N&=&
H^N_{I,J}Q^N_{I,J}+\sqrt{N}\Ibold_{|J|}\one_{[\![E^N_{I\setminus J}]\!]\geq1/2}\in\Mat_{|J|}(\SSS).
\end{eqnarray*}

\subsubsection{Abuses of notation}
We write
\[
\Delta^kR_{I,J}^N=(\Delta^N_{I,J})^kR_{I,J}^N\quad \mbox{and}\quad
\Delta R_{I,J}^N=\Delta^1R_{I,J}^N.
\]
We often write $j$ where we should more correctly write $\{j\}$,
for example, we write $Q^N_{I,j}$ instead of $Q^N_{I,\{j\}}$.
Note that
\[
R_{I,j}^N=R_I^N(j,j),\qquad H_{I,j}^N=H_{I\setminus j}^N,\qquad
Q^N_{I,j,j,j}=Q^N_{I,j},\qquad P^N_{I,j,j,j}=P^N_{I,j}.
\]
In the same spirit, we occasionally write $N$ in place of $\{1,\ldots
,N\}$.
\begin{Remark}\label{RemarkDeeplyRelevantzbold}
This is a continuation of Remark~\ref{RemarkCobbling}
and furthermore a continuation of the thread of remarks including Remarks
\ref{RemarkPrePreCobbling},~\ref{RemarkPreCobbling}
and~\ref{RemarkCobblingBis}.
Suppose now that $(\SSS,L,\Theta,e)$ is a self-adjoint linearization
of some $f\in\Mat_n({\CC\langle\Xbold\rangle})_\sa$.
Then with $\tau_{\SSS,e}$ as in Definition \ref
{DefinitionSSStracialState}, we have
%
%
\begin{equation}\label{equationDeeplyRelevant}
S_{\mu_f^N}(\zbold)=\tau_{\SSS,e}(F^N_N)\qquad \mbox{on the event
$\tbold=0$.}
\end{equation}
This representation of $S_{\mu_f^N}(\zbold)$ is just a rewrite of
\eqref{equationPreCobblingRandom}.
Next, let $(\underline{\SSS},\underline{L},\underline{\Theta
}+\Diamond_\SSS,\underline{e})$ be the ``underlined'' SALT block
design derived from $(\SSS,L,\Theta,e)$ via Lemma~\ref{LemmaForcedOnUs}.
Consider the triple
%
%
\begin{equation}\label{equationPreDeeplyRelevantUnderline}
(\underline{X},\underline{\Lambda},\underline{\Phi})=
\Bigl(\bigcup\underline{L}(\Xi^N),
\underline{\Theta}+\Diamond_\SSS+\zbold\underline{e}+\ii\tbold
1_{\underline{\SSS}},\Phi_{\underline{L}}\Bigr),
\end{equation}
which again satisfies assumption \eqref
{equationUrMatrixIdentitySetup} with probability $1$.
Let the recipes attached to the underlined triple $(\underline
{X},\underline{\Lambda},\underline{\Phi})$
be denoted with underlines in order to distinguish them from those
attached 
to the triple $(X,\Lambda,\Phi)$ defined by \eqref
{equationPreDeeplyRelevant}.
We then have
%
%
\begin{equation}\label{equationRaisonDetre1}
\partial_1\underline{F}^N_I=T^N_I
\quad\mbox{and}\quad
\partial_2\underline{F}^N_I=U^N_I
\end{equation}
by
\eqref{equationBulletTransposeRationale}, \eqref
{equationHalfTransposeRationale},
\eqref{equationUnderlineConstruction3} and \eqref
{equationUnderlineConstruction4}.
The relations \eqref{equationRaisonDetre1} joined with
the relations \eqref{equationRaisonDetre2}
will be crucial for the proof of relation \eqref{equationMainResultBis4}
of Theorem~\ref{TheoremMainResultBis}.
\end{Remark}

\begin{Remark}
The recipe $U^N_I$ does not figure in any identities stated in Section~\ref
{sectionUrMatrixIdentities}
but does become an important random variable later. We therefore
include its definition here
so that Section~\ref{subsectionRecipes} can serve as a handy catalog of
the basic random variables.
\end{Remark}

\begin{Remark}\label{RemarkLocalInverseObservation}
Note that
$R^N_I$ is the inverse of the matrix
${\mathbf{e}}_I(\frac{X}{\sqrt{N}}-\Ibold_\infty\otimes
\Lambda){\mathbf{e}}_I$ as computed in
the algebra ${\mathbf{e}}_I\Mat_\infty(\SSS){\mathbf{e}}_I$ the
identity element of which is ${\mathbf{e}}_I$.
This observation simplifies calculations below on several occasions.
\end{Remark}

\subsection{Basic identities} \label{subsectionMatrixIdentities}
We obtain block-type generalizations of matrix identities familiar
from the study of resolvents of Wigner matrices.

\begin{Lemma}\label{LemmaRescaling}
For $N$ and $I\in\III_{N}$, along with any positive integer $k$,
%
%
\begin{equation}\label{equationRescaling}
R^{N+1}_I
=R^N_I+\sum_{\nu=1}^{k-1}\biggl(\delta_NR^N_I\frac{{\mathbf{e}}_I
X{\mathbf{e}}_I^*}{\sqrt{N}}\biggr)^\nu R^N_I+
\biggl(\delta_N
R^N_I\frac{{\mathbf{e}}_I X{\mathbf{e}}_I^*}{\sqrt{N}}\biggr)
^kR^{N+1}_I,
\end{equation}
where $\delta_N=\sqrt{N}(\frac{1}{\sqrt{N}}-\frac{1}{\sqrt
{N+1}})$.
\end{Lemma}

\begin{pf}
By induction, we may assume $k=1$. Then, in view of Remark~\ref{RemarkLocalInverseObservation},
formula \eqref{equationRescaling} is merely an instance of the
resolvent identity \eqref{equationResolventIdentity}.~%
\end{pf}

\begin{Lemma}\label{LemmaConcreteVeryDry}
For $N$ and $(I,J)\in\III_N^{(2)}$,
%
%
\begin{eqnarray}\label{equationConcretePreVeryDry}
R^N_{I,J}& = &
\biggl(
\frac{\fbold_JX\fbold_J^*}{\sqrt{N}}-\Ibold_{|J|}\otimes\Lambda
-\frac{\fbold_JXR^N_{I\setminus J}X\fbold_J^*}{N}\biggr)^{-1},\\
\label{equationConcreteVeryDry}
R^N_I - R^N_{I\setminus J}
& = &
\biggl(\fbold_J^*-R^N_{I\setminus J}\frac{X}{\sqrt{N}}\fbold
_J^*\biggr)
R^N_{I,J}\biggl(\fbold_J-\fbold_J \frac{X}{\sqrt{N}}R^N_{I\setminus
J}\biggr).
\end{eqnarray}
In particular, we automatically have $R_{I,J}^{N}\in\GL_{|J|}(\SSS)$.
\end{Lemma}
\begin{pf} In Proposition~\ref{PropositionVeryDry}, let us now take
\begin{eqnarray*}
\AAA&=&{\mathbf{e}}_I\Mat_\infty(\SSS){\mathbf{e}}_I,\qquad x=
{\mathbf{e}}_I\biggl(\frac{X}{\sqrt{N}}-\Ibold_\infty\otimes
\Lambda\biggr){\mathbf{e}}_I,\qquad
\pi={\mathbf{e}}_J,\\
\pi^\perp&=&{\mathbf{e}}_{I\setminus J},\qquad
\sigma={\mathbf{e}}_I,
\end{eqnarray*}
and again let us take advantage of Remark~\ref{RemarkLocalInverseObservation}.
Rewritten in the form
\[
{\mathbf{e}}_JR^N_I{\mathbf{e}}_J=\fbold_J^*R^N_{I,J}\fbold_J=
\fbold_J^*
\bigl(\fbold_J(x-xR^N_{I\setminus J}x)\fbold_J^*
\bigr)^{-1}\fbold_J,
\]
identity \eqref{equationConcretePreVeryDry} becomes a special case of
\eqref{equationPreVeryDry}.
Similarly, rewritten in the form
\[
R^N_I-R^N_{I\setminus J}
=
({\mathbf{e}}_J-R^N_{I\setminus J}x{\mathbf{e}}_J
)R^N_I({\mathbf{e}}_J-{\mathbf{e}}_JxR^N_{I\setminus J}),
\]
identity \eqref{equationConcreteVeryDry}
becomes a specialization of \eqref{equationVeryDry}.
\end{pf}

\begin{Lemma}\label{LemmaRHQE}
For $N$ and $(I,J)\in\III_N^{(2)}$,
along with any positive integer $k$,
%
%
\begin{equation}
\label{equationMasterRIJbis}
R^N_{I,J}=H_{I,J}^N+
\sum_{\nu=1}^{k-1}
\frac{(H_{I,J}^NQ^N_{I,J})^\nu H_{I,J}^N}{N^{\nu/2}}
+\frac{\Delta^kR_{I,J}^N}{N^{k/2}}.
\end{equation}
\end{Lemma}

\begin{pf} By induction on $k$, we may assume $k=1$.
Rewrite \eqref{equationConcretePreVeryDry} in the form
%
%
\begin{equation}
\label{equationMasterRIJ}
-\bigl(\Ibold_{|J|}\otimes\bigl(\Lambda+\Phi(F^N_{I\setminus J})\bigr)\bigr)R^N_{I,J}=
\Ibold_{|J|}\otimes1_\SSS+\frac{Q^N_{I,J}R^N_{I,J}}{\sqrt{N}}.
\end{equation}
Then left-multiply by $H_{I,J}^N$ on both sides and rearrange slightly
to get the result.
\end{pf}

\subsection{More elaborate identities}\label{subsectionMoreElaborate}
We specialize and combine the basic identities.

\subsubsection{Comparison of $R$ to $H$}
For $N$, $(I,J)\in\III^{(2)}_N$ and $j_1,j_2\in J$,
we~have
%
%
\begin{eqnarray}\label{equationDoubletonCase1}
R^N_I(j_1,j_2)-\delta_{j_1j_2}H_{I\setminus J}&=&
\fbold_{j_1}\fbold_J^*\frac{\Delta R^N_{I,J}}{\sqrt{N}} \fbold_J
\fbold_{j_2}^*,\\
\label{equationDoubletonCase2}
\quad R^N_I(j_1,j_2)-\delta_{j_1j_2}H_{I\setminus J}-\frac{H_{I\setminus
J}Q^N_{I,J,j_1,j_2}H_{I\setminus J}}{\sqrt{N}}
&=&\fbold_{j_1}\fbold_J^*\frac{\Delta^2 R^N_{I,J}}{N} \fbold_J
\fbold_{j_2}^*
\end{eqnarray}
by merely rewriting \eqref{equationMasterRIJbis} in the cases $k=1$
and $k=2$, respectively,
at the level of individual matrix entries.

\subsubsection{Increments of $F$ and of $H$}
For $N$ and $(I,J)\in\III^{(2)}_N$ we have
%
%
\begin{eqnarray}
\label{equationMasterFIj}
N(F_I^N-F^N_{I\setminus J})
&=&\trace_\SSS(R^N_{I,J})+
\trace_\SSS\biggl(R_{I\setminus J}^{N}
\frac{X}{\sqrt{N}}\fbold_J^*
R^N_{I,J}\fbold_J\frac{X}{\sqrt{N}}R^{N}_{I\setminus J}\biggr)
\nonumber
\\[-8pt]
\\[-8pt]
\nonumber
&=&\biggl(\trace_\SSS+T_{I\setminus J}^N\circ\Phi\circ\trace_\SSS+
\frac{P_{I,J}^N}{\sqrt{N}}\biggr)(R^N_{I,J})
\end{eqnarray}
by applying $\trace_\SSS$ to both sides of \eqref{equationConcreteVeryDry}.
We note also the identity
%
%
\begin{eqnarray}\label{equationMasterHIJ}
H^N_I-H^N_{I\setminus J}&=&
H^N_I\one_{[\![E^N_{I\setminus J}]\!]\geq1/2}
-H_{I\setminus J}^N\one_{[\![E^N_I]\!]\geq1/2}
\nonumber
\\[-8pt]
\\[-8pt]
\nonumber
&&{}+H^N_I\Phi
(F_I^N-F_{I\setminus J}^N)H^N_{I\setminus J}
\end{eqnarray}
obtained by exploiting the resolvent identity \eqref
{equationResolventIdentity} in evident fashion.

\subsubsection{The Schwinger--Dyson error}
For $N$ and $I\in\III_N$ such that $|I|\geq2$ we have
%
%
\begin{equation}\label{equationUrSchwingerDyson}
E^N_I+\frac{|I|-N}{N}1_\SSS=
\frac{1}{N}\sum_{j\in I}
\biggl(\Phi(F^N_I-F^N_{I\setminus j})R^N_{I,j}-\frac
{Q^N_{I,j}R^N_{I,j}}{\sqrt{N}}\biggr)\vadjust{\goodbreak}
\end{equation}
after applying $\frac{1}{N}\sum_{j\in I}(\cdot)$ to both sides
of \eqref{equationMasterRIJ} in the singleton case $J=\{j\}$ and rearranging.
\begin{Remark}
Identity \eqref{equationUrSchwingerDyson} is an
approximate version of the Schwinger--Dyson equation.
Identities of this sort have long been in use for study of Wigner matrices.
\end{Remark}

\subsubsection{Refined Schwinger--Dyson error}
For $N$ and $I\in\III_N$ such that $|I|\geq2$, we also have
%
%
\begin{eqnarray}\label{equationUrSchwingerDysonCorrected}
&&E^N_I+\frac{|I|-N}{N}1_\SSS+\frac{1}{N}\sum_{j\in I}\frac
{Q^N_{I,j}H^N_{I\setminus j}}{\sqrt{N}}
\nonumber
\\[-9pt]
\\[-9pt]
\nonumber
&&\qquad=
\frac{1}{N}\sum_{j\in I}\biggl(
\Phi(F^N_I-F^N_{I\setminus j})R^N_{I,j}-\frac{Q^N_{I,j}\Delta
R^N_{I,j}}{N}\biggr),
\end{eqnarray}
by \eqref{equationMasterRIJbis} for $k=1$ in the singleton case $J=\{
j\}$
and \eqref{equationUrSchwingerDyson}, after rearrangement.
\subsubsection{Comparison of $F$ to $H$}
For $N$ and $I\in\III_N$, we have
%
%
\begin{equation}\label{equationHFComp1}
H_I^N-F^N_I=H_I^NE_I^N-F_I^N\one_{[\![E^N_I]\!]\geq1/2}
\end{equation}
by direct appeal to the definitions. One
then obtains
for $|I|\geq2$ the identity
%
%
\begin{eqnarray}\label{equationHFComp2}
&&H^N_I-F^N_I+\frac{1}{N}\sum_{j\in I}\frac{F_{I\setminus
j}^NQ^N_{I,j}H^N_{I\setminus j}}{\sqrt{N}}\nonumber\\[-1pt]
&&\qquad=\frac{N-|I|}{N}F_I^N+H_I^N(E_I^N)^2-(F_I^N+F_I^NE_I^N)\one
_{[\![E^N_I]\!]\geq1/2}
\nonumber
\\[-9pt]
\\[-9pt]
\nonumber
&&\qquad\quad{}+
\frac{1}{N}\sum_{j\in I}\biggl(
F^N_I\Phi(F^N_I-F^N_{I\setminus j})R^N_{I,j}-\frac
{F^N_IQ^N_{I,j}\Delta R^N_{I,j}}{N}\\[-1pt]
&&\hspace*{126pt}\qquad{}-\frac{(F_I^N-F^N_{I\setminus
j})Q^N_{I,j}H^N_{I\setminus j}}{\sqrt{N}}\biggr)\nonumber
\end{eqnarray}
by iterating \eqref{equationHFComp1} and combining
it with \eqref{equationUrSchwingerDysonCorrected}.

\subsubsection{Refined increment of $F$}
For $N$ and $(I,J)\in\III^{(2)}_N$ we have
%
%
\begin{eqnarray}
\label{equationMasterFIjBis}
&&F_I^N-F^N_{I\setminus J}-|J|\frac{(1_{B(\SSS)}+T^N_{I\setminus
J}\circ\Phi)
(H^N_{I\setminus J})}{N}
\nonumber
\\[-9pt]
\\[-9pt]
\nonumber
&&\qquad =
\frac{P^N_{I,J}(R_{I,J}^N)+(\trace_\SSS+T^N_{I\setminus J}\circ\Phi
\circ\trace_\SSS)
(\Delta R^N_{I,J})}{N^{3/2}}
\end{eqnarray}
by rearrangement of \eqref{equationMasterFIj}, using \eqref
{equationMasterRIJbis}
for $k=1$.

\subsubsection{Increment of $F$ with respect to $N$}
For $N$ and $I\in\III_N$,
we have
%
%
\begin{eqnarray}\label{equationFRescaling}
&&(N+1)F^{N+1}_I-NF^N_I-\frac{1}{2}\bigl(F^N_I+T^N_I(\Lambda)\bigr)\nonumber\\
&&\qquad=\frac{1}{N}\trace_\SSS\biggl(
\biggl(N\delta_N-\frac{1}{2}\biggr)
\bigl({\mathbf{e}}_I+R_I^N(\Ibold_\infty\otimes\Lambda)\bigr)
R_I^{N}\\
&&\hspace*{31pt}\qquad\quad{}
+\frac{(N\delta_N)^2}{N}
\bigl({\mathbf{e}}_I+R_I^N(\Ibold_\infty\otimes\Lambda
)\bigr)^2R_I^{N+1}\biggr)\nonumber
\end{eqnarray}
by Lemma~\ref{LemmaRescaling} in the case $k=2$
after using Remark~\ref{RemarkLocalInverseObservation},
applying $\trace_\SSS$ on both sides and rearranging.
Note that
$\frac{1}{2}-\frac{1}{2N}\leq N\delta_N\leq\frac{1}{2}$.

\subsubsection{The link recipe}
Let
\[
\Link^N=\frac{1}{2}\bigl(F^N_N+T^N_N(\Lambda)\bigr)-F^{N+1}_{N+1}
+H^{N+1}_N+T^{N+1}_N(\Phi
(H^{N+1}_N)),
\]
where here and below in similar contexts we abuse notation by writing $N$
where we should more correctly write $\{1,\ldots,N\}$.
We then have
%
%
\begin{eqnarray}\label{equationTrivialRearrangement}
&&N(F^{N+1}_{N+1}-F^N_N)-\Link^N\nonumber \\
&&\qquad= (N+1)F^{N+1}_N-NF^N_N-\tfrac{1}{2}\bigl(F^N_N+T^N_N(\Lambda)\bigr)
\\
&&\qquad\quad{}+(N+1)(F_{N+1}^{N+1}-F^{N+1}_N)-H^{N+1}_N-
T^{N+1}_N\bigl(\Phi
(H^{N+1}_N)\bigr)\nonumber
\end{eqnarray}
by mere rearrangement of terms.

\subsection{The bias identity}\label{subsectionBias}
We derive the most intricate identity used in the paper.
\subsubsection{Further refinement of the Schwinger--Dyson error} We first need an intermediate result which continues
the process of expansion begun in identity \eqref
{equationUrSchwingerDysonCorrected}.
For $N$ and $I\in\III_N$ such that $|I|\geq2$,
we have
%
%
\begin{eqnarray}\label{equationPreBias}
&&E_I^N+\frac{|I|-N}{N}1_\SSS\nonumber\hspace*{-35pt}\\
&&\quad{}+\frac{1}{N}\sum_{j\in I}\biggl(
\frac{(Q_{I,j}^NH_{I\setminus j}^N)^2-(\Phi+\Phi\circ
T^N_{I\setminus j}\circ\Phi)(R^N_{I,j})R^N_{I,j}}{N}
+\frac{(Q_{I,j}^NH_{I\setminus j}^N)^3}{N^{3/2}}\biggr)\hspace*{-35pt}\\
&&\qquad=\frac{1}{N}\sum_{j\in I}
\biggl(-\frac{Q^N_{I,j}H^N_{I\setminus j}}{\sqrt{N}}-\frac{Q_{I,j}^N
\Delta^3R_{I,j}^N}{N^2}+\frac{\Phi\circ
P^N_{I,j}(R^N_{I,j})R_{I,j}^N}{N^{3/2}}\biggr)\nonumber\hspace*{-35pt}
\end{eqnarray}
by expanding\vspace*{-1pt} the terms $\frac{Q^N_{I,j}R^N_{I,j}}{\sqrt{N}}$ in
\eqref{equationUrSchwingerDyson}
by using \eqref{equationMasterRIJbis} for $k=3$ in the singleton case
$J=\{j\}$,
and furthermore expanding the terms $\Phi(F^N_I-F^N_{I\setminus
j})R^N_{I,j}$ in \eqref{equationUrSchwingerDyson} by using \eqref
{equationMasterFIj} in the singleton case $J=\{j\}$,
after suitable rearrangement.
\subsubsection{The bias identity}
Fix $N\geq2$ and $j\in N$ arbitrarily.
To compactify notation put
\begin{eqnarray*}
\widetilde{T}^N_{j}&=&\Phi+\Phi\circ T^N_{N\setminus j}\circ\Phi,\qquad
\widetilde{P}^N_{j} = \Phi\circ P^N_{N\setminus j},\qquad
\check{R}^N_j=H^N_{N\setminus j}Q^N_{N, j}H^N_{N\setminus j},\\
\Error_j^N&=&
(Q_{N,j}^NH_{N\setminus j}^N)^2-
\widetilde{T}^N_j(H^N_{N\setminus j})H^N_{N\setminus
j}+(Q_{N,j}^NH^N_{N,j})^3/\sqrt{N},\\
\Error_j^{N,1}&=&
\bigl(\widetilde{T}^N_j(H^N_{N\setminus j})\check{R}^N_j+\widetilde
{T}^N_j(\check{R}^N_j)H^N_{N\setminus j}+\widetilde
{P}^N_j(H^N_{N\setminus j})H_{N\setminus
j}^N\bigr)/N-Q^N_{N,j}H^N_{N\setminus j},\\
\Error_j^{N,2}&=&
\widetilde{T}^N_j(H^N_{N\setminus j})\Delta^2R^N_{N,j}+
\widetilde{T}^N_j(\check{R}^N_j)\Delta R^N_{N,j}+\widetilde
{T}^N_j(\Delta^2 R^N_{N,j})R^N_{N,j}\\
&&{}
+\widetilde{P}^N_j(H^N_{N\setminus j})\Delta R^N_{N,j}+\widetilde
{P}^N_j(\Delta R^N_{N,j}) R^N_{N,j}-Q_{N,j}^N \Delta^3R_{N,j}^N.
\end{eqnarray*}
At last, we obtain the \textit{bias identity}
%
%
\begin{equation}\label{equationBias}
E_N^N+\frac{1}{N}\sum_{j\in N}\frac{\Error_{j}^N}{N}=\frac
{1}{N}\sum_{j\in N}\biggl(\frac{\Error^{N,1}_j}{\sqrt{N}}+\frac
{\Error^{N,2}_j}{N^2}
\biggr)
\end{equation}
by using \eqref{equationMasterRIJbis}
several times with $k=1,2$ in the singleton case $J=\{j\}$
to expand the terms $(\Phi+\Phi\circ T^N_{N\setminus j}\circ\Phi
)(R^N_{N,j})R^N_{N,j}$
and $\Phi\circ P^N_{N,j}(R^N_{N,j})R^N_{N,j}$ in \eqref{equationPreBias},
after suitable rearrangement.

\section{$L^p$ estimates for the block Wigner model}\label
{sectionMatrixIdentities}

We introduce a straightforward generalization of the usual Wigner
matrix model
with matrix entries in a block algebra.
Using all the tools collected in
Section~\ref{sectionConcentration} and Section \ref
{sectionUrMatrixIdentities},
we investigate how control of moments of ``randomized resolvents''
propagates to give control of moments of many related random variables.
Our main result in this section is Theorem~\ref{TheoremBias} below
which converts
identity \eqref{equationBias} above to a key approximation.

\subsection{The block Wigner model}
\label{subsectionBlockWigner}
The ad hoc infinite matrix formalism of
Section~\ref{subsectionInfiniteMatrices} will be the algebraic
framework for our discussion
of the block Wigner model.

\subsubsection{Data}
Data for the block Wigner model consist of
\begin{itemize}
\item a block algebra $\SSS$,
\item a random matrix $X\in\Mat_\infty(\SSS)_\sa$,
\item a (deterministic) linear map $\Phi\in B(\SSS)$,
\item a (deterministic) tensor $\Psi\in\SSS^{\otimes2}$,
\item a random element $\Lambda\in\SSS$ and
\item a random variable $\Gfrak\in[1,\infty)$.
\end{itemize}
%
\subsubsection{\texorpdfstring{The $\sigma$-fields $\FFF(i,j)$ and the auxiliary
random variables $\zbold$ and $\tbold$}
{The sigma-fields $\FFF(i,j)$ and the auxiliary
random variables $\zbold$ and $\tbold$}}

In addition to the data above, as a convenience of bookkeeping,
we keep for use in the present setup the same system
$\{\FFF(i,j)\}_{1\leq i\leq j<\infty}$ of independent
$\sigma$-fields mentioned\vadjust{\goodbreak} in Section~\ref{subsubsectionModelData}.
As before, let $\FFF$ denote the $\sigma$-field generated by all the
$\FFF(i,j)$.
More generally, for any set $I$ of positive integers, let $\FFF_I$
denote the $\sigma$-field
generated by the family $\{\FFF(i,j)\vert i,j\in I\}$.
We also keep the random variables $\zbold$ from Theorem \ref
{TheoremMainResultBis}
and $\tbold$ from Remark~\ref{RemarkCobbling} on hand,
and we continue to assume that $\FFF$, $\zbold$ and $\tbold$ are independent.
In Section~\ref{sectionMatrixIdentities}, the random variables
$\zbold$ and $\tbold$ intervene only through
their $\sigma$-fields but in Section~\ref{sectionEndgame} these
random variables themselves take an active part.

\subsubsection{Assumptions}
Of the sextuple $(\SSS,X,\Phi,\Psi,\Lambda,\Gfrak)$,
we assume the following:
%
%
\begin{eqnarray}
\label{equationPepper1X}\quad
&\displaystyle\sup_{i,j=1}^\infty{\Vert[\![X(i,j)]\!] \Vert}_p<\infty\qquad \mbox{for
$1\leq p<\infty$,}&\\
\label{equationPepper2X}
&\mbox{$X(i,j)$ is $\FFF(i\wedge j,i\vee j)$-measurable
and of mean zero for all $i$ and $j$,}&\\
\label{equationPepper3X}
\qquad\quad &\mbox{$\Phi=(\zeta\mapsto\Ebold X(i,j)\zeta X(j,i))$ and
$\Psi=\Ebold( X(i,j)^{\otimes2})$ for distinct $i$ and $j$,}&\\
\label{equationLambdaPepper1}
&
\mbox{$[\![\Lambda]\!]_p<\infty$\qquad for $p\in[1,\infty)$,}&\\
\label{equationLambdaPepper2}
&\mbox{$\Lambda$ is $\sigma(\zbold,\tbold)$-measurable and
$\Gfrak$ is $\sigma(\zbold)$-measurable,}&\\
\label{equationLambdaPepper3}
&\displaystyle\fbold_I\biggl(\frac{X}{\sqrt{N}}-\Ibold_\infty\otimes\Lambda
\biggr)\fbold_I^*\in\GL_{|I|}(\SSS)\qquad
\mbox{for $N$ and $I\in\III_N$.}&
\end{eqnarray}
%
For $N$ and $I\in\III_N$, we then put
\[
R^N_I=\fbold_I^*\biggl(\fbold_I\biggl(\frac{X}{\sqrt{N}}-\Ibold
_\infty\otimes\Lambda\biggr)\fbold_I^*\biggr)^{-1}\fbold_I
\in\Mat_\infty(\SSS),
\]
which is a generalized resolvent (Green's function).
Finally, we assume that
%
%
\begin{equation}\label{equationRNIBound}
\sup_N\bigvee_{I\in\III_N}{\Vert[\![R^N_I/\Gfrak]\!] \Vert
}_p<\infty.
\end{equation}
We work with a fixed instance $(\SSS,X,\Phi,\Psi,\Lambda,\Gfrak)$
of the block Wigner model over~$\SSS$
for the rest of Section~\ref{sectionMatrixIdentities}.

\begin{Remark}\label{RemarkSynthesisBis}
Fix real numbers $a$ and $b$.
Using assumptions \eqref{equationPepper1}--\eqref{equationPepper8}
it is easy to verify directly that the collection
%
%
\begin{eqnarray}\label{equationBlockWignerSpecializationBis}
&&(\SSS,X,\Phi,\Psi,\Lambda,\Gfrak)
\nonumber
\\[-8pt]
\\[-8pt]
\nonumber
&&\qquad=\biggl(\CC,\bigcup_N (a\Xi^N_1+b\Xi
^N_2),a^2+b^2,a^2-b^2,\zbold+\ii\tbold,1+1/\Im\zbold\biggr)
\end{eqnarray}
satisfies assumptions~\eqref{equationPepper1X}--\eqref{equationRNIBound}.
This specialization is merely a slight variant of the standard Wigner
matrix model.
The reader might find it helpful to concentrate on this case when
making a first pass through
the $L^p$ estimates assembled below.
\end{Remark}

\begin{Remark}\label{RemarkSynthesis}
All instances of the block Wigner model needed for the proof of Theorem
\ref{TheoremMainResultBis} arise as follows.
Let $(\SSS,L,\Theta,e)$ be a SALT block\vadjust{\goodbreak} design and let $c_0$, $c_1$,
$c_2$ and $\Tfrak$ be
the constants from Definition~\ref{DefinitionSALT}. We keep the
notation of Remark~\ref{RemarkCobbling}.
Let $\Phi_L$ and $\Psi_L$ be as in Definition~\ref{DefinitionPhiPsi}.
Then, using assumptions~\eqref{equationPepper1}--\eqref{equationPepper8}
along with Remark~\ref{RemarkYCapsuleComplement}, it is easy to
verify that the collection
%
%
\begin{eqnarray}\label{equationBlockWignerSpecialization}
&&(\SSS,X,\Phi,\Psi,\Lambda,\Gfrak)
\nonumber
\\
&&\qquad=\biggl(\SSS,\bigcup_N L(\Xi^N),\Phi_L,\Psi_L,\Theta
+\zbold e+\ii\tbold1_\SSS,\\
&&\hspace*{50pt}\qquad{}c_0\bigl(1+[\![L(\Xi)]\!]\bigr)^{c_1}(1+1/\Im
\zbold)^{c_2}\biggr)\nonumber
\end{eqnarray}
satisfies assumptions \eqref{equationPepper1X}--\eqref{equationRNIBound}.
\end{Remark}

\subsubsection{Random variables defined by recipes}
Since assumption \eqref{equationLambdaPepper3} is a verbatim
repetition of assumption
\eqref{equationUrMatrixIdentitySetup}, all the recipes of Section
\ref{subsectionRecipes}
define random variables in the present setting.
The object $R^N_I$ figuring in assumption \eqref{equationRNIBound} is
of course
a recipe. We now furthermore have
random variables $F^N_I$, $H^N_I$,
$T^N_I$, $U^N_I$, etc. at our disposal. The
compound objects $\Link^N$, $\Error^N_j$, etc. figuring in the more elaborate identities
also become random variables in the present setting.

\begin{Remark}\label{RemarkRNIBoundAmplification}
We note that \eqref{equationRNIBound} can be considerably refined in
the specialization \eqref{equationBlockWignerSpecialization} of the
block Wigner model.
Namely, for each $N$ and $I\in\III_N$ we have almost sure bounds
%
%
\begin{eqnarray}\label{equationRNIBoundAmplification}
[\![F^N_I]\!]&\leq&[\![R^N_I]\!]\leq\Gfrak\biggl(1+\biggl[\!\biggl[\frac{\fbold_I X\fbold
_I^*}{\sqrt{N}}\biggr]\!\biggr]\biggr)^{c_1},\\
\label{equationRNIBoundAmplificationBis}
[\![F^N_I]\!]&\leq&[\![R^N_I]\!]\leq\tfrac{1}{2} \qquad\mbox{on the
event $\tbold\geq\Tfrak$}.
\end{eqnarray}
One can also easily verify that
%
%
\begin{equation}\label{equationIsSDtunnel}
(G_L\dvtx \DDD_L\rightarrow\SSS,\Phi_L,\Lambda,\Tfrak,\Gfrak
)\qquad \mbox{is an SD tunnel.}
\end{equation}
Of course this SD tunnel is random since $\Lambda$ and $\Gfrak$ are random.
\end{Remark}

\subsubsection{Partially averaged random variables}
We complete our enumeration of the random variables we will be studying.
For $N$ and $I\in\III_N$, we define
\[
\overline{F}_I^N=\Gfrak\Ebold(F^N_I/\Gfrak\vert\zbold,\tbold)\in
\SSS
\quad\mbox{and}\quad
\overline{E}^N_I=1_\SSS+\bigl(\Lambda+\Phi(\overline{F}_I^N)\bigr)\overline
{F}_I^N\in\SSS.
\]
Since $[\![F^N_I]\!]/\Gfrak$ is integrable by assumption \eqref
{equationRNIBound},
in fact $\overline{F}_I^N$ and $\overline{E}_I^N$
are well-defined, almost surely.
Theorem~\ref{TheoremBias} below gives a delicate approximation to
$\overline{E}^N_N$.
(Recall our abuse of notation $N=\{1,\ldots,N\}$.)

\subsection{Basic estimates} We start gathering
consequences of \eqref{equationRNIBound}.

\subsubsection{\texorpdfstring{The norms ${|\!|\!|\cdot|\!|\!|}_{p,k}$}
{The norms |||.||| p,k}}

Given a constant $p\in[1,\infty)$, a positive integer $k$ and a
finite-dimensional-Banach-space-valued random variable $Z$
(defined on the same probability space as $\Gfrak$),
we write ${|\!|\!|Z|\!|\!|}_{p,k}={\Vert[\![Z]\!]/\Gfrak^k \Vert}_p$ to compress
notation.

\begin{Remark}
We emphasize that in Section~\ref{sectionMatrixIdentities} we make no
assumption concerning the strength of the repulsion of $\zbold$ from
the real axis. Indeed, we make no assumptions about $\zbold$ at all.
But nevertheless, looking ahead to the completion of the proof of
Theorem~\ref{TheoremMainResultBis}, we are obliged to keep track of
issues involving the repulsion strength.
To do so, we will use the parameter $k$ appearing
in the norm ${|\!|\!|\cdot|\!|\!|}_{p,k}$ and in a similar seminorm
introduced in Section~\ref{subsubsectionSemiNorm}
below.
\end{Remark}

\begin{Proposition}\label{PropositionRControl}
For each constant $p\in[1,\infty)$, we have
%
%
\begin{eqnarray}\label{equationUniform1}
\qquad\quad\sup_N\bigvee_{I\in\III_N}
{|\!|\!|R^N_I|\!|\!|}_{p,1}\vee{|\!|\!|F^N_I|\!|\!|}_{p,1}\vee{|\!|\!|H^N_I|\!|\!|}_{p,1}\vee{|\!|\!|T^N_I|\!|\!|}_{p,2}
\vee{|\!|\!|U^N_I|\!|\!|}_{p,2}&<&\infty,\\
\label{equationUniform1bis}
\sup_N\bigvee_{(I,J)\in\III_N^{(2)}}
{|\!|\!|R^N_{I,J}|\!|\!|}_{p,1}\vee{|\!|\!|H^N_{I,J}|\!|\!|}_{p,1}&<&\infty,\\
\label{equationUniform2}
\sup_N\bigvee_{(I,J)\in\III_N^{(2)}}
N{|\!|\!|F^N_I-F^N_{I\setminus J}|\!|\!|}_{p,3}&<&\infty.
\end{eqnarray}
\end{Proposition}
\begin{pf}
The claim made in \eqref{equationUniform1} for $R^N_I$ just repeats
the hypothesis \eqref{equationRNIBound} in different notation. We have
\[[\![R^N_I]\!]\geq[\![F^N_I]\!]
\vee[\![T^N_I]\!]^{1/2}\vee\frac{1}{2}[\![H^N_I]\!]\vee\frac
{1}{\sqrt{s}}[\![U^N_I]\!]^{1/2}
\]
obviously in the first two cases, by Lemma \ref
{LemmaNeumannExpansion} in the penultimate
case and Lemma~\ref{LemmaHilbertSchmidt} in the last, where $s^2$ is
the dimension of $\SSS$ over the complex numbers. Thus, \eqref
{equationUniform1} holds in general.
Clearly, we have
\[[\![R^N_{I,J}]\!]\leq[\![R^N_I]\!]
\quad\mbox{and}\quad [\![H^N_{I,J}]\!]=[\![H^N_{I\setminus J}]\!],
\]
whence \eqref{equationUniform1bis} via \eqref{equationUniform1}.
By Lemma~\ref{LemmaHilbertSchmidt} and identity \eqref{equationMasterFIj},
we have
\[
\frac{N}{|J|^2}[\![F_I^N-F^N_{I\setminus J}]\!]
\leq[\![R^N_I]\!]+\frac{s}{N}
[\![R^N_I]\!][\![R^N_{I\setminus J}]\!]^2
\sum_{(i,j)\in(I\setminus J)\times J}[\![X(i,j)]\!]^2.
\]
From this, estimate \eqref{equationUniform2} follows by
assumption \eqref{equationPepper1X}, the Minkowski inequality and
estimate \eqref{equationUniform1}.
\end{pf}

\subsubsection{\texorpdfstring{The seminorms ${|\!|\!|\cdot|\!|\!|}_{p,k,I}$}
{The seminorms |||.||| p,k,I}}\label{subsubsectionSemiNorm}
Given a constant $p\in[1,\infty)$, a positive integer $k$,
a set $I$ of positive integers\vadjust{\goodbreak}
and a finite-dimensional-Banach-space-valued random variable $Z$
defined on the same probability space as $\Gfrak$
such that ${|\!|\!|Z|\!|\!|}_{p,k}<\infty$, we define
\[
{|\!|\!|Z|\!|\!|}_{p,k,I}= {\Vert[\![\Ebold(Z/\Gfrak^k\vert\FFF
_I,\zbold,\tbold)]\!] \Vert}_p.
\]
Since the random variable $[\![Z/\Gfrak^k]\!]$ is assumed to be in
$L^p\subset L^1$,
the conditional expectation appearing on the right is well defined,
almost surely,
and moreover
\[
{|\!|\!|Z|\!|\!|}_{p,k}\geq{|\!|\!|Z|\!|\!|}_{p,k,I}\geq{|\!|\!|Z|\!|\!|}_{p,k,J}
\]
for any set $J\subset I$ by Jensen's inequality. In particular,
\[
{\Vert\Ebold(Z/\Gfrak^k\vert\zbold,\tbold) \Vert}_p={|\!|\!|Z|\!|\!|}_{p,k,\varnothing}
\]
whenever ${|\!|\!|Z|\!|\!|}_{p,k}<\infty$.

\begin{Proposition}\label{PropositionQControl}
For each constant $p\in[1,\infty)$, we have
%
%
\begin{eqnarray}
\label{equationUniform25}
\sup_N\bigvee_{(I,J)\in\III_N^{(2)}}
{|\!|\!|Q^N_{I,J}|\!|\!|}_{p,1}\vee{|\!|\!|P^N_{I,J}|\!|\!|}_{p,2}&<&\infty,\\
\label{equationUniform3}
\sup_N\bigvee_{(I,J)\in\III_N^{(2)}}
{|\!|\!|Q^N_{I,J}|\!|\!|}_{p,1,I\setminus J}\vee
{|\!|\!|P^N_{I,J}|\!|\!|}_{p,2,I\setminus J}&=&0.
\end{eqnarray}
\end{Proposition}

It is hard to overestimate the importance of this
proposition. This is the estimate
ultimately driving convergence. Our exploitation of it is of course an
imitation of the procedure of~\cite{BaiSil}.
\begin{pf*}{Proof of Proposition~\ref{PropositionQControl}}
Fix $N$, $(I,J)\in\III^{(2)}_N$ and $j_1,j_2\in J$ arbitrarily.
By definition, we have
%
%
\begin{eqnarray}\label{equationQExpansion}
Q^N_{I,J,j_1,j_2} &=& -X(j_1,j_2)+
\frac{1}{\sqrt{N}}\bigl(
\fbold_{j_1} XR_{I\setminus J}^NX\fbold_{j_2}^*-
N\delta_{j_1j_2}\Phi(F^N_{I\setminus J})\bigr),
\nonumber
\\[-8pt]
\\[-8pt]
\nonumber
\qquad P^N_{I,J,j_1,j_2} &=&
\biggl(\zeta\mapsto
\frac{1}{\sqrt{N}}\bigl(\trace_\SSS(R_{I\setminus J}^{N}
X\fbold_{j_1}^*
\zeta\fbold_{j_2}XR^{N}_{I\setminus J})
-N
\delta_{j_1j_2}T_{I\setminus J}^{N}(\Phi(\zeta))\bigr)
\biggr).
\end{eqnarray}
By \eqref{equationPepper1X} and \eqref{equationUniform1},
the random variables $[\![Q^N_{I,J,j_1,j_2}]\!]$
and $[\![P^N_{I,J,j_1,j_2}]\!]$ are integrable,
hence the conditional expectations
\[
\Ebold(Q^N_{I,J,j_1,j_2}/\Gfrak\vert\FFF_{I\setminus J},\zbold
,\tbold) \quad\mbox{and}\quad
\Ebold(P^N_{I,J,j_1,j_2}/\Gfrak^2\vert\FFF_{I\setminus J},\zbold
,\tbold)
\]
are well-defined
and vanish almost surely
by assumptions \eqref{equationPepper1X}, \eqref{equationPepper2X}
and \eqref{equationPepper3X}.
By Proposition~\ref{PropositionKeyEstimate},
Remark~\ref{RemarkTwoTypesBilinear},
estimate \eqref{equationUniform1}
and the hypotheses of the block Wigner model,
the quantities
\[
{|\!|\!|Q^N_{I,J,j_1,j_2}+X(j_1,j_2)|\!|\!|}_{p,1},\qquad
{|\!|\!|X(j_1,j_2)|\!|\!|}_{p,1}\quad \mbox{and}\quad
{|\!|\!|P^N_{I,J,j_1,j_2}|\!|\!|}_{p,2}
\]
are bounded uniformly in $N$, $I$, $J$, $j_1$ and $j_2$. Thus,
claims \eqref{equationUniform25} and \eqref{equationUniform3} hold.
\end{pf*}

\subsection{More elaborate estimates}
We combine and specialize the basic estimates.

\begin{Proposition}\label{PropositionEControl}
For each constant $p\in[1,\infty)$, we have
%
%
\begin{eqnarray}
\label{equationEControl}
\sup_N\mathop{\bigvee_{
I\in\III_N}}_{
\mathrm{s.t.}\ |I|\geq2}
\sqrt{N}{\biggl|\!\biggl|\!\biggl|E^N_I+\frac{|I|-N}{N}1_\SSS\biggr|\!\biggr|\!\biggr|}_{p,4}&<&\infty,\\
\label{equationDeltaControl}
\sup_N\mathop{\bigvee_{
(I,J)\in\III_N^{(2)}\ \mathrm{s.t.}}}_{|I|\geq N-\sqrt{N}
}
{|\!|\!|\Delta^N_{I,J}|\!|\!|}_{p,4}
&<&\infty,\\
\label{equationE2ControlBis}
\sup_N\mathop{\bigvee_{
I\in\III_N\ \mathrm{s.t.}}}_{
|I|\geq N-99
}
N{|\!|\!|E^N_I|\!|\!|}_{p,6,\varnothing}&<&\infty.
\end{eqnarray}
\end{Proposition}

\begin{pf} We take Propositions~\ref{PropositionRControl} and \ref
{PropositionQControl}
for granted at every step.
Identity~\eqref{equationUrSchwingerDyson} implies the estimate \eqref
{equationEControl}.
Estimate \eqref{equationEControl}
and the Chebyshev bound
%
%
\begin{equation}\label{equationEChebychevBound}
\one_{[\![E^N_I]\!]\geq1/2}\leq(2[\![E^N_I]\!])^c \qquad (c\geq0)
\end{equation}
imply estimate \eqref{equationDeltaControl}.
Identity \eqref{equationUrSchwingerDysonCorrected} and estimate
\eqref{equationDeltaControl}
imply the estimate
\[
\sup_N\mathop{\bigvee_{I\in\III_N}}_{
\mathrm{s.t.}\ |I|\geq N-\sqrt{N}
}
N{\biggl|\!\biggl|\!\biggl|E^N_I+\frac{|I|-N}{N}1_\SSS+\frac{1}{N}\sum_{j\in I}\frac
{Q^N_{I,j}H^N_{I\setminus j}}{\sqrt{N}}\biggr|\!\biggr|\!\biggr|}_{p,6}<\infty.
\]
Estimate \eqref{equationE2ControlBis} follows via
\eqref{equationUniform3}.
\end{pf}

\begin{Proposition}\label{PropositionFIncrementControl}
For each constant $p\in[1,\infty)$, we have
%
%
\begin{eqnarray}
\label{equationFRescalingControl}
\sup_N
N^{3/2}{\biggl|\!\biggl|\!\biggl|F^{N+1}_{N+1}-F^N_N-\frac{\Link^N}{N}\biggr|\!\biggr|\!\biggr|}_{p,7}&<&\infty
,\\
\label{equationEbarControl}
\sup_N\mathop{\bigvee_{I\in\III_N \ \mathrm{s.t.}}}_{
|I|\geq N-99} N^2{|\!|\!|\overline{E}^N_I-E^N_I|\!|\!|}_{p,14,\varnothing}&<&\infty.
\end{eqnarray}
\end{Proposition}

\begin{pf} We take Propositions~\ref{PropositionRControl}, \ref
{PropositionQControl}
and~\ref{PropositionEControl}
for granted at every step. We have
\[
\sup_N
N{\biggl|\!\biggl|\!\biggl|(N+1)F^{N+1}_N-NF^N_N-\frac{1}{2}\bigl(F^N_N+T^N_N(\Lambda )\bigr)\biggr|\!\biggr|\!\biggr|}_{p,3}<\infty
\]
by identity \eqref{equationFRescaling}
along with assumption \eqref{equationLambdaPepper1}.
The estimate
\[
\sup_N\mathop{\bigvee_{(I,J)\in\III_N^{(2)}}}_{
|I|\geq N-\sqrt{N}}
N^{3/2}{\biggl|\!\biggl|\!\biggl|F_I^N-F^N_{I\setminus J}-|J|\frac{(1_{B(\SSS
)}+T^N_{I\setminus J}\circ\Phi) (H^N_{I\setminus J})}{N}\biggr|\!\biggr|\!\biggr|}_{p,7}<\infty\vadjust{\goodbreak}
\]
follows from identity \eqref{equationMasterFIjBis}.
Estimate \eqref{equationFRescalingControl} then follows
via the definition of $\Link^N$.
From the last estimate above, it also follows that
\[
\sup_N \mathop{\sup_{I\in\III_N}}_{
|I|\geq N-99}N^2{\Vert\Var_\SSS(F^N_I/\Gfrak^7|\zbold,\tbold) \Vert
}_p<\infty
\]
via Proposition~\ref{PropositionVarianceBound}, whence estimate
\eqref{equationEbarControl}.
\end{pf}

\begin{Proposition}\label{PropositionRIncrementControl}
For each constant $p\in[1,\infty)$, we have
%
%
\begin{eqnarray}
\label{equationRIJp}
\sup_N\mathop{\bigvee_{
(I,J)\in\III_N^{(2)}\ \mathrm{s.t.}}}_{
|I|\geq N-\sqrt{N}} \bigvee_{j_1,j_2\in J}
\sqrt{N}{|\!|\!|R^N_I(j_1,j_2)-\delta_{j_1j_2}H^N_{I\setminus J}|\!|\!|}_{p,5}&<&\infty,\\
\label{equationRIJp2}
\sup_N\mathop{\bigvee_{
(I,J)\in\III_N^{(2)}\ \mathrm{s.t.}}}_{
|I|\geq N-\sqrt{N}} \bigvee_{j_1,j_2\in J}N{|\!|\!|R^N_I(j_1,j_2)-\delta
_{j_1,j_2}H^N_{I\setminus J}|\!|\!|}_{p,9,I\setminus J}&<&\infty,\\
\label{equationHIJp}
\sup_N\mathop{\bigvee_{
(I,J)\in\III_N^{(2)}\ \mathrm{s.t.}}}_{
|I|\geq N-\sqrt{N}
}
N{|\!|\!|H^N_I-H^N_{I\setminus J}|\!|\!|}_{p,9}&<&\infty,\\
\label{equationHF1comp}
\sup_N\mathop{\bigvee_{
I\in\III_N\ \mathrm{s.t.}}}_{
|I|\geq N-\sqrt{N}
}\sqrt{N}{|\!|\!|H^N_I-F^N_I|\!|\!|}_{p,5}&<&\infty,\\
\label{equationHF2comp}
\sup_N\mathop{\bigvee_{
I\in\III_N}}_{
|I|\geq N-\sqrt{N}
}N{|\!|\!|H^N_I-F^N_I|\!|\!|}_{p,9,I\setminus J}&<&\infty.
\end{eqnarray}
\end{Proposition}

\begin{pf} Taking Propositions~\ref{PropositionRControl}, \ref
{PropositionQControl}
and~\ref{PropositionEControl} for granted
and using again the Chebyshev bound \eqref{equationEChebychevBound},
one derives the estimates in question from identities~\eqref{equationDoubletonCase1},
\eqref{equationDoubletonCase2},
\eqref{equationMasterHIJ},
\eqref{equationHFComp1} and
\eqref{equationHFComp2},
respectively.
\end{pf}

\subsection{The bias theorem}\label{subsectionBiasProof}
We work out a delicate approximation to $\overline{E}_N^N$.
We use again the apparatus introduced to state and prove
Proposition~\ref{PropositionCubic},
as well as the cumulant and shuffle notation introduced in Section \ref
{subsectionUniversalCorrection}.

\subsubsection{Corrections}
For $N\geq2$ and $j=1,\ldots,N$ we define
\begin{eqnarray*}
\Corr^N_{j}
&=&\langle[\Psi,\Psi]_2,[U^N_{N\setminus j},(H^N_{N\setminus
j})^{\otimes2}]_2\rangle_4-\Phi(H^N_{N\setminus
j})H^N_{N\setminus j}\\
&&{}
+\langle\Ebold X(j,j)^{\otimes2},(H^N_{N\setminus j})^{\otimes
2}\rangle_2
-\frac{1}{\sqrt{N}}\langle\Ebold X(j,j)^{\otimes
3},(H^N_{N\setminus j})^{\otimes3}\rangle_3\\
&&{}+\frac{1}{N}\sum_{i\in N\setminus j}
\bigl\langle\Cbold^{(4)}(X(i,j)),[(R^N_{I\setminus j,i})^{\otimes
2},(H^N_{N\setminus j})^{\otimes2}]_2\bigr\rangle_4.
\end{eqnarray*}

\begin{Theorem}\label{TheoremBias}
For each constant $p\in[1,\infty)$, we have
%
%
\begin{equation}
\label{equationBiasEstimate}
\sup_{N\geq2}
N^2{\Biggl|\!\Biggl|\!\Biggl|\overline{E}^N_N+\frac{1}{N}\sum_{j=1}^N \frac{\Corr
_{j}^N}{N}\Biggr|\!\Biggr|\!\Biggr|}_{p,14,\varnothing}<\infty.
\end{equation}
\end{Theorem}

The proof of the theorem takes up the rest of Section \ref
{subsectionBiasProof}.
We need several lemmas.

\begin{Lemma}\label{LemmaErrControl}
For each constant $p\in[1,\infty)$, we have
%
%
\begin{eqnarray}
\label{equationBadErrorControl}
\sup_{N\geq2} \bigvee_{j=1}^N
{|\!|\!|\Error^N_j|\!|\!|}_{p,6}\vee{|\!|\!|\Error^{N,1}_j|\!|\!|}_{p,6}\vee
{|\!|\!|\Error^{N,2}_j|\!|\!|}_{p,14}&<&\infty,\\
\label{equationAverageErrorAway}
\sup_{N\geq2}\bigvee_{j=1}^N {|\!|\!|\Error ^{N,1}_j|\!|\!|}_{p,6,N\setminus j}&=&0.
\end{eqnarray}
\end{Lemma}

\begin{pf}
Taking Propositions~\ref{PropositionRControl},~\ref{PropositionQControl}
and~\ref{PropositionEControl} for granted, these facts can be read
off from the definitions
presented in Section~\ref{subsectionBias}.
\end{pf}

\subsubsection{Moment notation}
For any sequence ${\mathbf{i}}=i_1\cdots i_{2k}$ of positive integers
and positive integer $j$ not appearing in ${\mathbf{i}}$ put
\begin{eqnarray*}
\Mbold_j({\mathbf{i}})&=&\Ebold\Bigl[\bigl(X(j,i_1)\otimes
X(i_2,j)-\Ebold\bigl(X(j,i_1)\otimes X(i_2,j)\bigr)\bigr)
\otimes\cdots\\
&&\quad{}\otimes\bigl(X(j,i_{2k-1})\otimes X(i_{2k},j)-\Ebold
\bigl(X(j,i_{2k-1})\otimes X(i_{2k},j)\bigr)\bigr)\Bigr]\\
&\in&\SSS^{\otimes2k}.
\end{eqnarray*}

\begin{Lemma}\label{LemmaMboldProps}
For sequences ${\mathbf{i}}=i_1\cdots i_{2k}$
of positive integers, and positive integers $j$ not appearing in
${\mathbf{i}}$,
the following statements hold:
\begin{longlist}[(III)]
\item[(I)] For each fixed $k$, $[\![\Mbold_j({\mathbf{i}})]\!]$ is
bounded uniformly in
${\mathbf{i}}$ and $j$.
\item[(II)] $\Mbold_j({\mathbf{i}})$ vanishes unless
$\Pi({\mathbf{i}})\in\Part^*(2k)$.
\item[(III)] If $\Pi({\mathbf{i}})\in\Part^*_2(2k)$,
then $\Mbold_j({\mathbf{i}})$ depends only on
$\Pi({\mathbf{i}})$.\vadjust{\goodbreak}
\end{longlist}
\end{Lemma}

\begin{pf} Assumption \eqref{equationPepper1X} implies statement (I).
Assumptions \eqref{equationPepper2X} implies statement (II).
Assumptions \eqref{equationPepper2X} and \eqref{equationPepper3X}
imply statement (III).
\end{pf}

\subsubsection{Tensor products of resolvent entries}
For $N$, $I\in\III_N$ and sequences ${\mathbf
{i}}=i_1\cdots i_{2k}\in\Seq(2k,I)$ put
\[
R^N_I({\mathbf{i}})=R^N_I(i_1,i_2)\otimes\cdots\otimes
R^N_I(i_{2k-1},i_{2k})\in\SSS^{\otimes k}.
\]

\subsubsection{\texorpdfstring{The random variable $\Rub^N_{j}$}{The random variable Rub N j}}

For $N\geq2$ and $j=1,\ldots,N$ put
\begin{eqnarray*}
\Rub^N_{j}&=&\frac{1}{N^2}\mathop{\mathop{\sum_{
{\mathbf{i}}
\in\Seq(6,N\setminus j)\ \mathrm{s.t.}}}_{
\Pi({\mathbf{i}})\in\Part^*(6)\ \mathrm{and}}}_{
\Pi({\mathbf{i}})\sim\{\{1,2,3\},\{4,5,6\}\}
}
\langle\Mbold_{j}({\mathbf{i}}),[R^N_{N\setminus j}({\mathbf
{i}}),(H^N_{N\setminus j})^{\otimes3}]_3\rangle_6.
\end{eqnarray*}
Here, we employ again the notation $\sim$ for $\Gamma_3$-orbit equivalence
previously introduced in connection with the list \eqref
{equationPartitionList}.

\begin{Lemma}\label{LemmaConditionalCalc}
For $N\geq3$ and $j=1,\ldots,N$
we have
\[[\![\Gfrak^6\Ebold(\Error^N_j/\Gfrak^6\vert\FFF_{N\setminus
j},\zbold,\tbold)-\Corr^N_j-\Rub^N_{j}]\!]
\leq\frac{c}{N}[\![R^N_{N\setminus j}]\!]^6,
\]
almost surely, for a constant $c$ independent of $N$ and $j$.
\end{Lemma}

\begin{pf}
%
In the case $(I,J)=(N,\{j\})$, formula \eqref{equationQExpansion}
above simplifies to
\[
Q^N_{N,j}+X(j,j)
= \frac{1}{\sqrt{N}}\bigl(
\fbold_jXR_{N\setminus j}^NX\fbold_{j}^*-\Gfrak\Ebold
(
\fbold_{j} XR_{N\setminus j}^NX\fbold_{j}^*/\Gfrak
\vert\FFF_{N\setminus j},\zbold,\tbold)\bigr).
\]
Note that the right-hand side is independent of $X(j,j)$.
A straightforward calculation
using Lemma~\ref{LemmaMboldProps}(II) yields that for $k\in\{2,3\}$,
\begin{eqnarray*}
&&\Gfrak^{2k}\Ebold\bigl((Q^N_{N,j}H^N_{N\setminus j})^{k}/\Gfrak
^{2k}\vert\FFF_{N\setminus j},\zbold,\tbold\bigr)
-(-1)^k\langle\Ebold X(j,j)^{\otimes k},(H^N_{N\setminus j})^{\otimes
3}\rangle_3\\
&&\qquad=\frac{1}{N^{k/2}}\mathop{\sum_{
{\mathbf{i}}
\in\Seq(2k,N\setminus j)}}_{
\mathrm{s.t.}\ \Pi({\mathbf{i}})\in\Part^*(2k)}
\langle\Mbold_j({\mathbf{i}}),[R^N_{N\setminus j}({\mathbf
{i}}),(H^N_{N\setminus j})^{\otimes k}]_k\rangle_{2k}.
\end{eqnarray*}
By a calculation using Lemma~\ref{LemmaMboldProps}(II, III)
and enumeration \eqref{equationPartitionEnumeration},
with $\alpha,\beta\in N\setminus j$ arbitrarily chosen
distinct elements,
we have
\begin{eqnarray*}
&&\Gfrak^4\Ebold\bigl((Q^N_{N,j}H^N_{N\setminus j})^2/\Gfrak^4\vert\FFF
_{N\setminus j},\zbold,\tbold\bigr)
-\langle\Ebold(X(j,j)^{\otimes2}),(H^N_{N\setminus j})^{\otimes
2}\rangle_2\\
&&\qquad=\frac{1}{N}\sum_{i_1,i_2\in N\setminus j}
\langle\Mbold_{j}(\alpha\beta\alpha\beta),[R^N_{N\setminus
j}(i_1,i_2)^{\otimes2}
,(H^N_{N\setminus j})^{\otimes2}]_2\rangle_4\nonumber\\
&&\qquad\quad{}+\frac{1}{N}\sum_{i_1,i_2\in N\setminus j}
\langle\Mbold_{j}(\alpha\beta\beta\alpha),[R_{N\setminus
j}^N(i_1,i_2)\otimes R_{N\setminus j}^N(i_2,i_1),(H^N_{N\setminus
j})^{\otimes2}]_2\rangle_4\nonumber\\
&&\qquad\quad{}+\frac{1}{N}\sum_{i\in N\setminus j}\langle\Mbold
_{j}(iiii)-\Mbold_{j}(\alpha\beta\beta\alpha)-\Mbold_{j}(\alpha
\beta\alpha\beta),\\
&&\hspace*{123pt}\qquad{}
[(R^N_{N\setminus j,i})^{\otimes2},(H^N_{N\setminus j})^{\otimes
2}]_2\rangle_4\nonumber\\
&&\qquad=\langle[\Psi,\Psi]_2,
[U^N_{N\setminus j},(H^N_{N\setminus j})^{\otimes2}]_2\rangle_4+
\Phi\circ T^N_{N\setminus j}\circ\Phi(H^N_{N\setminus j})\\
&&\qquad\quad{}+\frac{1}{N}\sum_{i\in N\setminus j}\langle\Cbold
^{(4)}(X(i,j)),[(R^N_{N\setminus j,i})^{\otimes2},(H^N_{N\setminus
j})^{\otimes2}]_2\rangle_4.\nonumber
\end{eqnarray*}
It follows that
\begin{eqnarray*}
&&\Gfrak^6\Ebold(\Error^N_j/\Gfrak^6\vert\FFF_{N\setminus
j},\zbold,\tbold)-\Corr^N_j-
\Rub^N_j\\
&&\qquad=\frac{1}{N^2}\mathop{\mathop{\sum_{
{\mathbf{i}}
\in\Seq(6,N\setminus j)\ \mathrm{s.t.}}}_
{\Pi({\mathbf{i}})\in\Part^*(6)\ \mathrm{and}}}_{
\Pi({\mathbf{i}})\not\sim\{\{1,2,3\},\{4,5,6\}\}
}
\langle\Mbold_j({\mathbf{i}}),[R^N_{N\setminus j}({\mathbf
{i}}),(H^N_{N\setminus j})^{\otimes3}]_3\rangle_6,
\end{eqnarray*}
whence the result by Proposition~\ref{PropositionCubic}
and Lemma~\ref{LemmaMboldProps}(I, III).
\end{pf}

\begin{Lemma}\label{LemmaNub}
Fix $p\in[1,\infty)$ arbitrarily. For $N\geq3$ and distinct
$j,j_1,j_2\in N$,
the quantity
\[
N{|\!|\!|(H_{N\setminus j}^N)^{\otimes3} \otimes R^N_{N\setminus
j}(j_1j_1j_2j_2j_1j_2)|\!|\!|}_{p,14,N\setminus\{
j,j_1,j_2\}}
\]
is bounded uniformly in $N$, $j$, $j_1$ and $j_2$.
\end{Lemma}

\begin{pf}
Put $J=\{j,j_1,j_2\}$.
The quantity
\[
N{|\!|\!|(H^N_{N\setminus J})^{\otimes5}\otimes R^N_{N\setminus
j}(j_1,j_2)|\!|\!|}_{p,14,N\setminus J}
\]
is bounded uniformly in $N$, $j$, $j_1$ and $j_2$
by \eqref{equationUniform1} and \eqref{equationRIJp2}.
The quantity
\[
N{|\!|\!|(H_{N\setminus j}^N)^{\otimes3}\otimes R^N_{N\setminus
j}(j_1j_1j_2j_2j_1j_2) -(H^N_{N\setminus J})^{\otimes5}\otimes
R^N_{N\setminus j}(j_1,j_2) |\!|\!|}_{p,14}
\]
is bounded uniformly in $N$, $j$, $j_1$ and $j_2$
by \eqref{equationRIJp} and \eqref{equationHIJp}.
\end{pf}

\subsubsection{\texorpdfstring{Completion of the proof of Theorem \protect\ref{TheoremBias}}
{Completion of the proof of Theorem 5}}
\label{subsubsectionProofOfTheoremBias}
We have
\[
\sup_{N\geq2}N^2{|\!|\!|\overline{E}^N_N- E^N_N|\!|\!|}_{p,14,\varnothing
}<\infty
\]
by estimate \eqref{equationEbarControl}. We have
\[
\sup_{N\geq2}
N^2{\biggl|\!\biggl|\!\biggl|E^N_N+\frac{1}{N}\sum_{j=1}^N \frac{\Error_j^N}{N}\biggr|\!\biggr|\!\biggr|}_{p,14,\varnothing}<\infty
\]
by the bias identity
and Lemma~\ref{LemmaErrControl}.
We have
\[
\sup_{N\geq2}\bigvee_{j=1}^N
N{|\!|\!|\Error^N_j-\Corr_j^N-\Rub^N_j|\!|\!|}_{p,14}<\infty
\]
by Proposition~\ref{PropositionRControl} and
Lemma~\ref{LemmaConditionalCalc}.
Finally, we have
\[
\sup_{N\geq2}\bigvee_{j=1}^N N{|\!|\!|\Rub^N_j|\!|\!|}_{p,14,\varnothing
}<\infty
\]
by Lemma~\ref{LemmaMboldProps}(I) and Lemma~\ref{LemmaNub}, which
finishes the proof.


\section{Concluding arguments}\label{sectionEndgame}
We finish the proof of Theorem~\ref{TheoremMainResultBis}.
\subsection{Setup for the concluding arguments}
Throughout Section~\ref{sectionEndgame}, we fix an instance
%
%
\begin{equation}\label{equationBasicSALTBlockDesign}
(\SSS,L,\Theta,e),\qquad c_0,\qquad c_1,\qquad c_2 ,\qquad\Tfrak
\end{equation}
of Definition~\ref{DefinitionSALT}
and we work with the corresponding instance
\begin{eqnarray*}
&&(\SSS,X,\Phi,\Psi,\Lambda,\Gfrak)\\
&&\qquad=\biggl(\SSS,\bigcup_N L(\Xi^N),\Phi_L,\Psi_L,\Theta+\zbold
e+\ii\tbold1_\SSS,c_0\bigl(1+[\![L(\Xi)]\!]\bigr)^{c_1}(1+1/\Im\zbold
)^{c_2}\biggr)
\end{eqnarray*}
of the block Wigner model exhibited in Remark~\ref{RemarkSynthesis}.
We emphasize that we must consider a general example of
a SALT block design because, at various stages below, we have to
consider both a SALT block design
arising as a self-adjoint linearization, that is, via Proposition \ref
{PropositionSALTnouveau},
and also a SALT block design arising by the underline construction,
that is,
via Lemma~\ref{LemmaForcedOnUs}.

Given any $\sigma(\FFF,\zbold,\tbold)$-measurable
Banach-space-valued integrable random variable
$Z$, we define $Z\vert_{\tbold=0}$ to be any $\sigma(\FFF,\zbold
)$-measurable
random variable which on the event $\tbold=0$ equals $Z$ almost surely.
For example, we have $\Lambda\vert_{\tbold=0}=\Theta+\zbold e$.
Since the latter random variable intervenes frequently below, we
will write $\Lambda_0=\Lambda\vert_{\tbold=0}$ to compress notation.

In a similar vein, given $Z$ as above along with a $\sigma$-field
$\GGG$
on which $\zbold$ is measurable, we abuse notation by writing $\Ebold
(Z\vert\GGG)=
\Gfrak^k\Ebold(Z/\Gfrak^k\vert\GGG)$ when there exists some
positive integer $k$
large enough so that $Z/\Gfrak^k$ is integrable
and hence the conditional expectation $\Ebold(Z/\Gfrak^k\vert\GGG)$
is well defined.

We will employ the abbreviated notation
\[
G\dvtx \DDD\rightarrow\SSS
\]
in place of the more heavily subscripted notation
$G_L\dvtx \DDD_L\rightarrow\SSS$. In a similar spirit, we write
\[
G'=\DD[G],\qquad \check{G}=\bigl((G^{-1})^{\otimes2}-\Psi\bigr)^{-1}, \qquad\Bias
^N=\Bias^N_L.
\]
Note also that for every $p\in[1,\infty)$ the bounds
%
%
\begin{equation}\label{equationNoGintegrabilityIssues}
\qquad{|\!|\!|G(\Lambda)|\!|\!|}_{p,1}\vee
{|\!|\!|G'(\Lambda)|\!|\!|}_{p,2}\vee{|\!|\!|\check{G}(\Lambda)|\!|\!|}_{p,2}\vee
\sup_N{|\!|\!|\Bias^N(\Lambda)|\!|\!|}_{p,5}<\infty
\end{equation}
hold, as one checks by means of Remark~\ref{RemarkCapsules}. We also have
%
%
\begin{equation}\label{equationNoGintegrabilityIssuesBis}
{|\!|\!|G(\Lambda)^{-1}|\!|\!|}_{p,1}<\infty
\end{equation}
by the SD equation $1_\SSS+(\Lambda+\Phi(G(\Lambda)))G(\Lambda)=0$ and
assumption \eqref{equationLambdaPepper1}.
\subsection{\texorpdfstring{Application of Proposition \protect\ref{PropositionModifiedKahuna}}
{Application of Proposition 18}}

We claim that
%
%
\begin{eqnarray}\label{equationGoodbyeGothic1}
&&\sup_N\mathop{\bigvee_{
I\in\III_N\ \mathrm{s.t.}}}_{
|I|\geq N-99
}\sqrt{N}{|\!|\!|F^N_I\vert_{\tbold=0}-G(\Lambda_0)|\!|\!|}_{p,99}<\infty,
\\\label{equationGoodbyeGothic2}
&&\sup_N\mathop{\bigvee_{
I\in\III_N\ \mathrm{s.t.}}}_{
|I|\geq N-99
}N|\!|\!|F^N_I\vert_{\tbold=0}+G'(\Lambda_0;(E^N_I\vert_{\tbold
=0})G(\Lambda_0)^{-1})
\nonumber
\\[-8pt]
\\[-8pt]
\nonumber
&&\hspace*{172pt}\qquad{}-G(\Lambda_0)|\!|\!|_{p,99}<\infty,\\
\label{equationGoodbyeGothic3}
&&\sup_N N^2{|\!|\!|\overline{F}^N_N\vert_{\tbold=0}+G'(\Lambda
_0;(\overline{E}^N_N\vert_{\tbold=0})G(\Lambda_0)^{-1})-G(\Lambda
_0)|\!|\!|}_{p,99}<\infty.
\end{eqnarray}
To prove the claim, we introduce several further random variables. Put
\[
\Cfrak=99e^{2\Tfrak}(1+[\![\Phi]\!]+[\![\Theta]\!]+|\zbold|).
\]
For $N$ and $I\in\III_N$ put
\[
\Lfrak^N_I=\Gfrak^2\biggl(1+\biggl[\!\biggl[\frac{\fbold_IX\fbold _I^*}{\sqrt{N}}\biggr]\!\biggr]\biggr)^{2c_1},\qquad
\Efrak^N_I=\Ebold( [\![E^N_I]\!]\vert\FFF,\zbold).
\]
Also for $N$ put
\[
\overline{\Lfrak}^N=\Gfrak^2\Ebold\biggl(1+\biggl[\!\biggl[\frac{\fbold _NX\fbold
_N^*}{\sqrt{N}}\biggr]\!\biggr]\biggr)^{2c_1},\qquad
\overline{\Efrak}^N = \Ebold( [\![\overline {E}^N_N]\!]\vert
\zbold).
\]
By Proposition~\ref{PropositionModifiedKahuna} applied conditionally,
with help from Remarks~\ref{RemarkYCapsuleComplement} and
\ref{RemarkRNIBoundAmplification} to check hypotheses, we have
\begin{eqnarray*}
&&\sqrt{N}[\![F^N_I\vert_{\tbold=0}-G(\Lambda_0)]\!]
\leq\sqrt{N}(\Cfrak\Gfrak\Lfrak^N_I)^6\bigl(\Efrak^N_I+(\Efrak
^N_I)^2\bigr),\\
&&N[\![F^N_I\vert_{\tbold=0}+G'(\Lambda_0;(E^N_I\vert_{\tbold
=0})G(\Lambda_0)^{-1})-G(\Lambda_0)]\!]\\
&&\qquad\leq N(\Cfrak\Gfrak\Lfrak^N_I)^{12}\bigl((\Efrak^N_I)^2+(\Efrak
^N_I)^4\bigr),\\
&&N^2[\![\overline{F}^N_N\vert_{\tbold=0}+G'(\Lambda _0;(\overline
{E}^N_N\vert_{\tbold=0})G(\Lambda_0)^{-1})-G(\Lambda _0)]\!]\\
&&\qquad\leq N^2(\Cfrak\Gfrak\overline{\Lfrak}^N)^{12}\bigl((\overline{\Efrak
}^N)^2+(\overline{\Efrak}^N)^4\bigr).
\end{eqnarray*}
Now fix $p\in[1,\infty)$ arbitrarily.
The right sides above can be bounded in the norm ${|\!|\!|\cdot |\!|\!|}_{p,k}$ for suitably chosen $k$, as follows.
Firstly, $\Cfrak$ has moments of all orders.
Secondly, we are in effect allowed to ignore factors of $\Gfrak$
on the right sides above at the expense of increasing $k$.
Thirdly, we have
\[
\sup_N\bigvee_{I\in\III_N}
{|\!|\!|\Lfrak^N_I|\!|\!|}_{p,2}<\infty,\qquad
\sup_N {|\!|\!|\overline{\Lfrak}^N|\!|\!|}_{p,2}<\infty
\]
by assumption \eqref{equationPepper2}.
Fourthly, we have
\[
\sup_N\mathop{\bigvee_{
I\in\III_N \mathrm{s.t.}}}_{
|I|\geq N-\sqrt{N}
}\sqrt{N}{|\!|\!|\Efrak^N_I|\!|\!|}_{p,4}<\infty,\qquad
\sup_NN{|\!|\!|\overline{\Efrak}^N|\!|\!|}_{p,14}<\infty,
\]
via \eqref{equationEControl}, \eqref{equationE2ControlBis}, \eqref
{equationEbarControl}
and Jensen's inequality.
The claims \eqref{equationGoodbyeGothic1}, \eqref
{equationGoodbyeGothic2} and \eqref{equationGoodbyeGothic3}
are proved.

From \eqref{equationGoodbyeGothic2}, we then deduce
%
%
\begin{equation}\label{equationGoodbyeGothic4}
\sup_N\mathop{\bigvee_{
I\in\III_N \ \mathrm{s.t.}}}_{|I|\geq N-99
}N{|\!|\!|F^N_I\vert_{\tbold=0}-G(\Lambda_0)|\!|\!|}_{p,99,\varnothing
}<\infty
\end{equation}
via \eqref{equationE2ControlBis}, \eqref
{equationNoGintegrabilityIssues} and \eqref
{equationNoGintegrabilityIssuesBis}.

\subsection{\texorpdfstring{Proof of statement \protect\eqref{equationMainResultBis1}
of Theorem \protect\ref{TheoremMainResultBis}}
{Proof of statement (22) of Theorem 4}}

In this paragraph, we assume that $(\SSS,L,\Theta,e)$ is a
self-adjoint linearization of
some $f\in\Mat_n(\CC\langle\Xbold\rangle)$. Then (recall) we have formulas
%
%
\begin{equation}\label{equationInitialTrope}
\tau_{\SSS,e}(G(\Lambda_0))=S_{\mu_f}(\zbold)
\quad\mbox{and}\quad \tau_{\SSS,e}(F^N_N\vert_{\tbold=0}
)=S_{\mu_f^N}(\zbold)
\end{equation}
by Remarks~\ref{RemarkCobblingBis} and \ref
{RemarkDeeplyRelevantzbold}, respectively.
Thus, we have
\[[\![S_{\mu_f^N}(\zbold)-S_{\mu_f}(\zbold)]\!]\leq
[\![F^N_N\vert_{\tbold=0}-G(\Lambda_0)]\!].
\]
Now fix $p\in[1,\infty)$ arbitrarily.
By \eqref{equationGoodbyeGothic1}, it follows that
\[
\sup_N\sqrt{N}{|\!|\!|S_{\mu_f^N}(\zbold)-S_{\mu_f}(\zbold )|\!|\!|}_{2p,99}<\infty.
\]
Now this last bound holds no matter what strength of repulsion of
$\zbold$ from the real axis we choose.
Let us now choose the repulsion strength strong enough so that ${\Vert
\Gfrak^{99} \Vert}_{2p}<\infty$.
Then we reach the desired conclusion \eqref{equationMainResultBis1}.

The preceding proof explains by example how bounds in the norm ${|\!|\!|\cdot|\!|\!|}_{p,k}$
with $k$ independent of $p$
translate to bounds in the norm ${\Vert[\![\cdot]\!] \Vert}_p$
provided that the strength of repulsion of $\zbold$
from the real axis is sufficiently strong, depending on $p$.
In the remainder of the proof of Theorem~\ref{TheoremMainResultBis},
we will omit similar details of translation.

\subsection{\texorpdfstring{Easy consequences of \protect\eqref{equationGoodbyeGothic1} and
\protect\eqref{equationGoodbyeGothic4}}
{Easy consequences of (176) and (179)}}
Estimates \eqref{equationGoodbyeGothic1} and \eqref{equationGoodbyeGothic4}
along with Propositions~\ref{PropositionQControl} and \ref
{PropositionRIncrementControl} yield the following bounds:
%
%
\begin{eqnarray}\label{equationCheap1}
\sup_N\mathop{\sup_{
I\in\III_N}}_{
|I|\geq N-99
}\sqrt{N}{|\!|\!|H^N_I\vert_{\tbold=0}-G(\Lambda_0)|\!|\!|}_{p,99}&<&\infty
,\\
\label{equationCheap2}
\sup_N\mathop{\sup_{
I\in\III_N}}_{
|I|\geq N-99
}N{|\!|\!|H^N_I\vert_{\tbold=0}-G(\Lambda_0)|\!|\!|}_{p,99,\varnothing
}&<&\infty,\\
\label{equationCheap3}
\sup_N\mathop{\sup_{
(I,J)\in\III_N^{(2)}}}_{
|I|\geq N-99
}
\bigvee_{j_1,j_2\in J}\sqrt{N}{|\!|\!|R^N_I(j_1,j_2)\vert_{\tbold
=0}-\delta_{j_1,j_2}G(\Lambda_0)|\!|\!|}_{p,99}&<&\infty,\\
\label{equationCheap4}
\sup_N\mathop{\sup_{
(I,J)\in\III_N^{(2)}}}_{
|I|\geq N-99
}
\bigvee_{j_1,j_2\in J}N{|\!|\!|R^N_I(j_1,j_2)\vert_{\tbold=0}-\delta
_{j_1,j_2}G(\Lambda_0)|\!|\!|}_{p,99,\varnothing}&<&\infty.
\end{eqnarray}

\subsection{Bootstrapping: Application of the secondary trick}
Let
%
%
\begin{equation}\label{equationBasicSALTBlockDesignBis}
(\underline{\SSS},\underline{L},\underline{\Theta}+\Diamond_\SSS
,\underline{e}),\qquad
\underline{c}_0,\qquad\underline{c}_1,\qquad\underline{c}_2,\qquad\underline{\Tfrak}
\end{equation}
be the instance of Definition~\ref{DefinitionSALT}
obtained by applying the underline construction to the instance \eqref
{equationBasicSALTBlockDesign}.
Consider as well the corresponding instance
\begin{eqnarray*}
&&(\underline{\SSS},\underline{X},\underline{\Phi},\underline
{\Psi},\underline{\Lambda},\underline{\Gfrak})\\
&&\qquad=\biggl(\underline{\SSS},\bigcup_N \underline{L}(\Xi^N),\Phi
_{\underline{L}},\Psi_{\underline{L}},\underline{\Theta}+\Diamond
_\SSS+\zbold\underline{e}+\ii\tbold1_{\underline{\SSS}},\\
&&\hspace*{76pt}\qquad{}
\underline{c}_0\bigl(1+[\![\underline{L}(\Xi)]\!]\bigr)^{\underline
{c}_1}(1+1/\Im\zbold)^{\underline{c}_2}\biggr)\nonumber
\end{eqnarray*}
of the block Wigner model constructed in Remark~\ref{RemarkSynthesis}.
By Lemma~\ref{LemmaForcedOnUs}, we can take
$\underline{c}_0=3c_0^2$, $\underline{c}_1=2c_1$ and $\underline{c}_2=2c_2$.
It follows that we can take $\underline{\Gfrak}=\frac{3}{2}\Gfrak^2$.
By~\eqref{equationRaisonDetre2} and \eqref{equationRaisonDetre1}
in combination with \eqref{equationGoodbyeGothic1} and \eqref
{equationGoodbyeGothic4},
we thus obtain bounds
%
%
\begin{eqnarray}\label{equationUnderhanded1}
\sup_N\mathop{\sup_{
I\in\III_N}}_{
|I|\geq N-99}\sqrt{N}{|\!|\!|T^N_I\vert_{\tbold=0}-G'(\Lambda _0)|\!|\!|}_{p,199}&<&\infty,\\
\label{equationUnderhanded3}
\sup_N\mathop{\sup_{
I\in\III_N}}_{
|I|\geq N-99
}\sqrt{N}{|\!|\!|U^N_I\vert_{\tbold=0}-\check{G}(\Lambda _0)|\!|\!|}_{p,199}&<&\infty,\\
\label{equationUnderhanded4}
\sup_N\mathop{\sup_{
I\in\III_N}}_{
|I|\geq N-99}N{|\!|\!|U^N_I\vert_{\tbold=0}-\check{G}(\Lambda _0)|\!|\!|}_{p,199,\varnothing}&<&\infty.
\end{eqnarray}
We can dispense now with the underlined SALT block design \eqref
{equationBasicSALTBlockDesignBis}
for the rest of the proof. We just needed it to get the estimates
immediately above.

\subsection{\texorpdfstring{Proof of statement \protect\eqref{equationMainResultBis2}
of Theorem \protect\ref{TheoremMainResultBis}}
{Proof of statement (23) of Theorem 4}}
Using again \eqref{equationInitialTrope},
we see that it is enough to prove for every $p\in[1,\infty)$ that
\[
\sup_N N^{3/2}{|\!|\!|F^{N+1}_{N+1}\vert_{\tbold=0}-F^N_N\vert
_{\tbold=0}|\!|\!|}_{p,999}<\infty.
\]
In turn, by estimate \eqref{equationFRescalingControl}, it is enough
to prove that
\[
\sup_N N^{1/2}{|\!|\!|\Link^N\vert_{\tbold=0}|\!|\!|}_{p,999}<\infty.
\]
But the latter follows in a straightforward way from
\eqref{equationSchwingerDysonCorrectionSpecial}, \eqref
{equationGoodbyeGothic1}, \eqref{equationCheap1} and~\eqref{equationUnderhanded1}.

\subsection{The last estimate}
We pause to explain in general terms how we are going to estimate
the seminorm ${|\!|\!|\cdot|\!|\!|}_{p,k,\varnothing}$
applied to the difference between a tensor product
of random variables of the form $U^N_I\vert_{\tbold=0}$, $H^N_I\vert
_{\tbold=0}$ and
$R^N_I(j,j)\vert_{\tbold=0}$ on the one hand
and a tensor product of random variables of the form
$G(\Lambda_0)$ and $\check{G}(\Lambda_0)$ on the other.
It is worthwhile to have a relatively abstract discussion of the method now
so that we can skip an unpleasant proliferation of indices below.

Let $A_1,\ldots,A_m\in\SSS$ be random and $\sigma(\FFF,\zbold)$-measurable.
Let $B_1,\ldots,B_m\in\SSS$ be random and $\sigma(\zbold)$-measurable.
Let $k_1,\ldots,k_m$ be positive integers and put $k=k_1+\cdots+k_m$.
Assume that for every $p\in[1,\infty)$ we have
\[
\bigvee_{i=1}^m{|\!|\!|A_i|\!|\!|}_{p,k_i}\vee\bigvee_{i=1}^m {|\!|\!|B_i|\!|\!|}_{p,k_i}<\infty.
\]
Now put
\[
A^{(0)}_i=B_i,\qquad A_i^{(1)}=A_i-B_i-\Ebold\bigl((A_i-B_i)\vert\zbold\bigr),
\qquad
A^{(2)}_i=\Ebold\bigl((A_i-B_i)\vert\zbold\bigr),
\]
noting that
\[
A_i=A_i^{(0)}+A_i^{(1)}+A_i^{(2)}.
\]
We then have for every $p\in[1,\infty)$ that
%
%
\begin{equation}\label{equationLastSmallTrick}
{|\!|\!|A_1\otimes\cdots A_m-B_1\otimes\cdots\otimes B_m|\!|\!|}_{p,k,\varnothing}
\leq
\mathop{\sum_{
(\nu_1,\ldots,\nu_k)\in\{0,1,2\}^m}}_{
\nu_1+\cdots+\nu_m\geq2
}
\prod_{i=1}^m{\bigl|\!\bigl|\!\bigl|A^{(\nu_i)}_i\bigr|\!\bigr|\!\bigr|}_{mp,k_i}\hspace*{-35pt}
\end{equation}
after taking into account the most obvious cancellations and applying
the H\"{o}lder inequality.
Roughly speaking, \eqref{equationLastSmallTrick} is advantageous because
in the intended application, we have $A_i^{(\nu)}=O(\frac{1}{N^{\nu/2}})$.

\subsection{\texorpdfstring{Proof of statements \protect\eqref{equationMainResultBis3}
and \protect\eqref{equationMainResultBis4}
of Theorem \protect\ref{TheoremMainResultBis}}
{Proof of statements (24) and (25) of Theorem 4}}
By Remark~\ref{RemarkCobblingTer}, the bound \eqref
{equationNoGintegrabilityIssues}
and yet another application of \eqref{equationInitialTrope},
it suffices to\vadjust{\goodbreak} prove that
\[
\sup_N N^2{\biggl|\!\biggl|\!\biggl|\overline{F}^N_N\vert_{\tbold=0}-\frac{\Bias
^N(\Lambda_0)}{N}-G(\Lambda_0)\biggr|\!\biggr|\!\biggr|}_{p,9999}<\infty.
\]
Using Theorem~\ref{TheoremBias}, \eqref{equationGoodbyeGothic3}
and \eqref{equationInitialTrope} above
it suffices to prove
\[
\sup_N\bigvee_{j=1}^N
N{|\!|\!|\widehat{\Bias}_L^N(\Lambda_0)- \Corr^N_j\vert_{\tbold=0}|\!|\!|}_{p,999}<\infty.
\]
Finally, this last bound is obtained by using the general observation
\eqref{equationLastSmallTrick} in conjunction with assumption \eqref
{equationPepper1}
and the estimates
\eqref{equationCheap1}, \eqref{equationCheap2}, \eqref{equationCheap3},
\eqref{equationCheap4}, \eqref{equationUnderhanded3} and \eqref
{equationUnderhanded4} above. The proof of Theorem \ref
{TheoremMainResultBis} is complete.

\section*{Acknowledgments}
I thank K. Dykema for teaching me Lemma~\ref{LemmaDykemaHelp}
and its proof, as well as explaining its application to free
semicircular variables.
I thank O. Zeitouni for suggesting the use of
an interpolation argument (as formalized in Lemma~\ref{LemmaGiddyUp})
to reduce Theorem~\ref{TheoremMainResult}
to Theorem~\ref{TheoremMainResultBis}.
I thank both of my co-authors A.~Guionnet and O.~Zeitouni for
teaching me much of direct relevance to this paper
in the process of writing the book~\cite{AGZ}. Especially Guionnet's
forceful advocacy of the Schwinger--Dyson
equation was influential.

%


\printaddresses

\end{document}